\documentclass[a4paper, final, 12pt]{article}
\usepackage{amsmath, amssymb, latexsym, amscd, amsthm,amsfonts,amstext}
\usepackage[mathscr]{eucal}
\usepackage{graphicx}\graphicspath{{figure/}}
\usepackage{subcaption}
\usepackage{float}
\usepackage{color}
\usepackage{hyperref}
\usepackage[utf8]{inputenc}
\usepackage[english]{babel}

\numberwithin{equation}{section}
\newtheorem{test}{Test}[section]
\newcommand{\dint}{\displaystyle\int\limits}

\begin{document}

\title{ Two 3-Dimensional Coefficient Inverse Problems for Wave-like
Equations With a Single Location of the Point Source\\
}
\author{Michael V. Klibanov, Jingzhi Li and Tian Niu}
\date{}
\maketitle

\begingroup
\footnotesize
\tableofcontents
\endgroup
\clearpage

\section{Introduction}

\label{sec:1}

In this paper two long standing open questions are addressed for the first
time. Those questions are about Lipschitz stability estimates and uniqueness
theorems for two Coefficient Inverse Problems (CIPs) for 3-Dimensional
wave-like Partial Differential Equations (PDEs) with formally determined
input data. The unknown coefficients depend on all three spatial variables.
The key novelty here is that the input data for these CIPs are generated by
only a single location of the point source at the moment of time $\left\{
t=0\right\} $. In the first CIP (CIP1) the unknown coefficient is the
potential of a wave-like PDE. The second CIP (CIP2) is concerned with a far
more challenging case when the unknown coefficient is involved in the
principal part of the hyperbolic operator. 

In addition to Lipschitz stability and uniqueness results, numerical methods
are developed for both CIPs. Their convergence analysis establishes the
global convergence property of both numerical techniques. The latter is
unlike the local convergence of conventional numerical procedures for CIPs.
The method for CIP1 is tested computationally. Its extensive numerical
testing for complicated media structures consistently demonstrates a very
accurate performance of this method for noisy input data.

\textbf{Definition 1.1}. \emph{A numerical method for a CIP is called
\textquotedblleft globally convergent" if there is a theorem, which
guarantees that iterations of this method deliver points in a sufficiently
small neighborhood of the true solution of this CIP without any advanced
knowledge of this neighborhood. On the other hand, a numerical method for a
CIP is called \textquotedblleft locally convergent" if a theorem guarantees
its convergence to that true solution only if the starting point of its
iterations is located in a sufficiently small neighborhood of that solution. 
}

Below $\mathbf{x}=\left( x,y,z\right) \in \mathbb{R}^{3}$ denotes spatial
points$,$ $T>0$ is a number, and $t\in \left( -\infty ,T\right) $ is time.
To our best knowledge, the only similar results for the case of a single
point source are H\"{o}lder stability and uniqueness theorems of the paper 
\cite{KR} for CIP2. However, it is assumed in \cite{KR} that a system of two
coupled PDEs derived from the original wave-like equation is written in
finite differences with respect to two spatial variables $x$ and $y$. On the
other hand, we use in this paper a significantly lighter assumption that a
certain associated PDE is written in finite differences with respect to $t$.
A high accuracy of our computational results confirms the adequacy of this
assumption. In addition, the Lipschitz stability estimates of the current
paper are stronger than the H\"{o}lder stability estimate of \cite{KR}. A
numerical method is not a part of \cite{KR}.

A CIP for a PDE is the problem of the determination of an unknown
coefficient(s) of that PDE from measurements of boundary data. Measured
boundary data are called \textquotedblleft input data". Since solution of a
PDE depends nonlinearly on its coefficients, then any CIP is nonlinear. Let $%
p$ be the number of free variables in the input data for a CIP\ and $n$ be
the number of free variables in the unknown coefficient. The input data are
called \textquotedblleft formally determined" ones if%
\begin{equation}
p=n.  \label{1.1}
\end{equation}%
Throughout this paper $p=n=3.$

As to condition (\ref{1.1}), we refer now to uniqueness theorems of
Berezansky for some CIPs for the Schr\"{o}dinger equation in $\mathbb{R}^{3}$
with the unknown potential \cite{Ber}. However, condition (\ref{1.1}) is
violated in \cite{Ber} via $p=5>n=3$. In addition, it is assumed in \cite%
{Ber} that the unknown potential is a piecewise analytic function. On the
other hand, it is desirable to impose basically only the assumption that the
unknown coefficient of a governing PDE belongs to one of main function
spaces, such as, e.g. $C^{k}.$ The latter is the case of the current paper.
Stability estimates were not obtained in \cite{Ber}, and a numerical method
was not developed.

The case, which is both the most economical one for practical measurements
and the most challenging one for analytical studies, is the case when the
input data for a CIP are generated by a single measurement event, which is
considered here. Currently all stability and uniqueness theorems for CIPs
with single measurement data and with condition (\ref{1.1}) at $n\geq 2$ are
obtained only by the method of the paper \cite{BukhKlib}. That paper is the
first one, where Carleman estimates were introduced in the field of Inverse
Problems. This idea was later explored in many publications. Since the
current paper is not a survey, we refer now for some follow up works only to 
\cite{Bell,Boul,Hamr,Isakov,Klib92,Ksurvey,KL,KR,MFG,Liu} and publications
cited therein. However, the technique of \cite{BukhKlib} works only under
the assumption that a certain initial condition does not vanish in the
entire domain of interest, which is not the case of the current paper. 

As stated above, we develop here two new numerical methods for our two CIPs
and test computationally one of them. They are versions of the so-called
convexification numerical method, which has the global convergence property
in terms of Definition 1.1. The convexification was originated in two
publications \cite{Klib95,Klib97} with the goal to avoid the well known
phenomenon of multiple local minima and ravines of conventional non-convex
least squares cost functionals of CIPs, see, e.g. \cite%
{Beilina,Beilina2,Chavent} for those functionals. The paper \cite{Scales}
contains a numerical example of multiple local minima. Minimization methods
for those functionals can be trapped by local minima and ravines.

While the initial publications \cite{Klib95,Klib97} about the
convexification method were purely analytical ones, more recently a number
of works was published, in which the convergence analysis is combined with
numerical studies for a variety of versions of the convexification method
for many CIPs, see, e.g. \cite{Bak,Khoa,KL,Kind,TFD,MFG,moving} and
references cited therein for some samples of those publications. We also
refer to \cite{QY} for the most recent work on the convexification method
for the case of the Ginzburg-Landau equation.

The convexification method is the numerical development of the idea of \cite%
{BukhKlib}. This method is a concept rather than a ready-to-use instruction.
The latter means that each new CIP requires the development of its own
version of the convexification method. Given an appropriate CIP with
formally determined input data, as in (\ref{1.1}), convexification
constructs a weighted Tikhonov-like functional $J_{\lambda }$ for it, where $%
\lambda \geq 1$ is a parameter. The weight is the Carleman Weight Function
(CWF). This function is involved as the weight in the Carleman estimate for
the corresponding PDE operator. The functional $J_{\lambda }$ is minimized
on a certain convex bounded set $G\left( \rho \right) \subset H,$ where $%
\rho >0$ is the diameter of $G\left( \rho \right) $ and $H$ is an
appropriate Hilbert space.

The key theorem of the convexification method states that for an appropriate
choice of $\lambda $ the functional $J_{\lambda }$ is strongly convex on $%
G\left( \rho \right) $ and has a unique minimizer on this set. Furthermore,
the distance between that minimizer and the true solution of that CIP is
estimated via the level of noise $\sigma $ in the input data. That minimizer
is called \textquotedblleft regularized solution" in the regularization
theory, and an estimate of its accuracy is an important requirement of this
theory \cite{T}. Finally, it is proven that if $\sigma \rightarrow 0,$ then
the gradient descent method of the minimization of $J_{\lambda }$ converges
to the true solution of the CIP if starting from an arbitrary point of a
set, which is similar with $G\left( \rho /3\right) .$ Furthermore, an
explicit convergence rate of the gradient descent method is derived. Since a
smallness condition is not imposed on $\rho $, then this is the global
convergence in terms of Definition 1.1.

\section{Statements of Two Coefficient Inverse Problems}

\label{sec:2}

In this paper, the governing forward problem is 
\begin{equation}
m\left( \mathbf{x}\right) u_{tt}=\Delta u+q\left( \mathbf{x}\right) u+\delta
\left( t\right) \delta \left( \mathbf{x}\right) ,\text{ }\mathbf{x}\in 
\mathbb{R}^{3},\text{ }t\in \left( -\infty ,T\right) ,  \label{2.1}
\end{equation}%
\begin{equation}
u\mid _{t\leq 0}=0.  \label{2.2}
\end{equation}%
Thus, (\ref{2.1}) and (\ref{2.2}) imply that the wave propagation process
starts at the moment of time $\left\{ t=0\right\} $ at the single location
of the point source at $\left\{ \mathbf{x}=0\right\} .$ As to equation (\ref%
{2.1}), in acoustics $1/\sqrt{m\left( \mathbf{x}\right) }$ is the speed of
the propagation of acoustic waves \cite[page 206]{CK}. Consider now the case
when $m\left( \mathbf{x}\right) \equiv 1.$ Let the function $p\left( \mathbf{%
x}\right) $ be the density of the medium, where acoustical waves propagate.
Let 
\begin{equation}
\left. 
\begin{array}{c}
p\mathbf{\in }C^{2}\left( \mathbb{R}^{3}\right) , \\ 
p\left( \mathbf{x}\right) \geq 1\text{ in }\mathbb{R}^{3}, \\ 
p\left( \mathbf{x}\right) =1\text{ for }x\in \mathbb{R}^{3}\diagdown \Omega .%
\end{array}%
\right. \text{ }  \label{2.20}
\end{equation}%
Then in (\ref{2.1}) \cite{Ber} 
\begin{equation}
q\left( \mathbf{x}\right) =-\sqrt{p\left( \mathbf{x}\right) }\Delta \left( 
\frac{1}{\sqrt{p\left( \mathbf{x}\right) }}\right) .  \label{2.03}
\end{equation}%
Therefore, CIPs considered here have applications in acoustics.

CIPs of this paper are complicated in their own rights. Hence, to simplify
the presentation as well as to simplify numerical studies in section 10, we
assume that our domain of interest $\Omega \subset \mathbb{R}^{3}$ is just a
rectangular prism,%
\begin{equation}
\Omega =\left\{ \mathbf{x}=\left( x,y,z\right) :-A<x,y<A,a<z<B\right\} ,
\label{2.3}
\end{equation}%
where $a,A,B>0$ are some numbers and $a<B$. By (\ref{2.3}) 
\begin{equation}
0\notin \overline{\Omega }.  \label{2.4}
\end{equation}%
Hence, by (\ref{2.4}), the point source in (\ref{2.1}) at the moment of time 
$\left\{ t=0\right\} $ in the position $\left\{ \mathbf{x}=0\right\} $ is
located outside of the closed domain $\overline{\Omega }.$ Denote%
\begin{equation}
Q_{T}=\Omega \times \left( 0,T\right) ,\text{ }S_{T}=\partial \Omega \times
\left( 0,T\right) ,  \label{2.5}
\end{equation}%
where $\partial \Omega $ is the boundary of the domain $\Omega .$ Obviously, 
$\partial \Omega =$ $\Gamma _{0}\cup \Gamma _{1}\cup \Gamma _{2},$ where  
\begin{equation}
\left. 
\begin{array}{c}
\Gamma _{0}=\left\{ \mathbf{x}=\left( x,y,z\right) :\left\vert x\right\vert
,\left\vert y\right\vert <A,\text{ }z=a\right\} , \\ 
\Gamma _{1}=\left\{ \mathbf{x}=\left( x,y,z\right) :\left\vert x\right\vert
,\left\vert y\right\vert <A,\text{ }z=B\right\} , \\ 
\Gamma _{2}=\left\{ \mathbf{x}=\left( x,y,z\right) :x,y=\pm A,z\in \left(
a,B\right) \right\} , \\ 
\text{ }\Gamma _{i,T}=\Gamma _{i}\times \left( 0,T\right) ,\text{ }i=0,1,2.%
\end{array}%
\right.   \label{2.6}
\end{equation}%
We assume everywhere below that coefficients of equation (\ref{2.1}) satisfy
the following conditions: 
\begin{equation}
q\in C^{4}\left( \mathbb{R}^{3}\right) ,  \label{2.7}
\end{equation}%
\begin{equation}
q\left( \mathbf{x}\right) =0\text{ for }\mathbf{x}\notin \Omega ,
\label{2.8}
\end{equation}%
\begin{equation}
q\left( \mathbf{x}\right) \geq 0\text{ in }\mathbb{R}^{3},  \label{2.9}
\end{equation}%
\begin{equation}
m\left( \mathbf{x}\right) \in C^{16}\left( \mathbb{R}^{3}\right) ,
\label{2.10}
\end{equation}%
\begin{equation}
m(\mathbf{x})\geq 1,~\mathbf{x}\in \mathbb{R}^{3},  \label{2.11}
\end{equation}%
\begin{equation}
m(\mathbf{x})=1\text{ for }\mathbf{x}\notin \Omega .  \label{2.12}
\end{equation}%
We need the regularity condition (\ref{2.7}) in section 4: to apply Lemma 1
of \cite{KR1}. The regularity condition (\ref{2.10}) follows from the
construction of the solution of Cauchy problem (\ref{2.1}), (\ref{2.2}) in
Theorem 4.1 of the book \cite{Rom}.

\textbf{Remark 2.1}. \emph{Note that traditionally extra regularity
conditions are of a secondary concern in the theory of CIPs, see, e.g. \cite%
{Nov,Rom,Rom1}.}

\textbf{Coefficient Inverse Problem 1 (CIP1)}. In equation (\ref{2.1}), let
the coefficient $m\left( \mathbf{x}\right) \equiv 1.$ Let the coefficient $%
q\left( \mathbf{x}\right) $ satisfies conditions (\ref{2.7})-(\ref{2.9}).
Suppose that the function $q\left( \mathbf{x}\right) $ is unknown for $%
\mathbf{x}\in \Omega .$ On the other hand, assume that the following
functions $g_{0}\left( \mathbf{x},t\right) $ and $g_{1}\left( \mathbf{x}%
,t\right) $ are known:%
\begin{equation}
u\mid _{S_{T}}=g_{0}\left( \mathbf{x},t\right) ,\text{ }  \label{2.13}
\end{equation}%
\begin{equation}
u_{z}\mid _{\Gamma _{0,T}}=g_{1}\left( \mathbf{x},t\right) .  \label{2.14}
\end{equation}%
Find the function $q\left( \mathbf{x}\right) $ for $\mathbf{x}\in \Omega .$

\textbf{Coefficient Inverse Problem 2 (CIP2).} In equation (\ref{2.1}), let
the coefficient $q\left( \mathbf{x}\right) \equiv 0$ and let the coefficient 
$m\left( \mathbf{x}\right) $ be unknown for $\mathbf{x}\in \Omega .$. In
addition, assume that conditions (\ref{2.10})-(\ref{2.12}) hold. Assume that
functions $g_{0}\left( \mathbf{x},t\right) $ and $g_{1}\left( \mathbf{x}%
,t\right) $ in (\ref{2.13}) and (\ref{2.14}) are known. Find the function $%
m\left( \mathbf{x}\right) $ for $\mathbf{x}\in \Omega .$

We now come back to the above mentioned application of CIP1 in acoustics.
Suppose that the function $q\left( \mathbf{x}\right) $ in (\ref{2.03}) is
found. Assume that the function $p\left( \mathbf{x}\right) $ in (\ref{2.03})
satisfies conditions (\ref{2.20}). Consider the function $\eta \left( 
\mathbf{x}\right) =\left( p\left( \mathbf{x}\right) \right) ^{-1/2}.$ Then,
using (\ref{2.20}) and (\ref{2.03}), we obtain 
\begin{equation}
\left. 
\begin{array}{c}
\Delta \eta +q\left( \mathbf{x}\right) \eta =0, \\ 
\eta \mid _{\partial \Omega }=1,\text{ }\partial _{\nu }\eta \mid _{\partial
\Omega }=0,%
\end{array}%
\right.   \label{2.114}
\end{equation}%
where $\partial _{\nu }$ is the normal derivative at $\partial \Omega .$
Thus, we have here both Dirichlet and Neumann for the linear elliptic
equation in the first line of (\ref{2.114}) with the known coefficient $%
q\left( \mathbf{x}\right) $. Uniqueness of the solution $\eta \in
C^{2}\left( \overline{\Omega }\right) $ of problem (\ref{2.114}) is well
known, see, e.g. \cite[\S 1 of Chapter 4]{LRS} and \cite[section 2.6.2]{KL}.
Furthermore, this solution, if it exists, can be found numerically using via
a modification of the version of \cite{KT} of the so-called
\textquotedblleft quasi-reversibility" method. Convergence of this method is
established in \cite{KT}. For brevity, we do not discuss further here the
issue of this paragraph. 

\section{Finite Differences With Respect to $t$}

\label{sec:3}

Consider a partition of the interval $\left[ 0,T\right] $ with the grid step
size $h$, 
\begin{equation}
0=t_{0}<t_{1}<...<t_{k}=T,\text{ }t_{i}-t_{i-1}=h,\text{ }k>2.  \label{3.1}
\end{equation}%
We now define $t-$finite differences of the function $f\left( t\right) \in
C^{2}\left[ 0,T\right] $ on the discrete set (\ref{3.1}). Using Taylor
formula, it is easy to verify that formulas below have the $O\left( h\right) 
$ accuracy as $h\rightarrow 0^{+},$ although we fix $h$ below. We set for
all points $t_{j}$ of the grid, except of the point $t_{k}=T$ 
\begin{equation}
\partial _{t}^{h}f\left( t_{j}\right) =\frac{f\left( t_{j+1}\right) -f\left(
t_{j}\right) }{h},\text{ }j=0,...,k-1.  \label{3.2}
\end{equation}%
In addition, for $t=t_{k}=T$ we set as in \cite{TFD} 
\begin{equation}
\partial _{t}^{h}f\left( t_{k}\right) =\frac{3f\left( t_{k}\right) -4f\left(
t_{k-1}\right) +f\left( t_{k-2}\right) }{2h}.  \label{3.3}
\end{equation}%
Everywhere below in this paper, for any function $f\left( t\right) \in C^{2}%
\left[ 0,T\right] $ we denote $f^{h}\left( t\right) $ the following vector
function%
\begin{equation}
f^{h}\left( t\right) =\left( f\left( t_{0}\right) ,f\left( t_{1}\right)
,...,f\left( t_{k}\right) \right) ^{T}.  \label{3.4}
\end{equation}

We also define the discrete Volterra integral of the function $f\left(
t\right) $ as in \cite{TFD}:%
\begin{equation}
\left. 
\begin{array}{c}
\left( \displaystyle\int \limits_{0}^{t_{i}}f\left( s\right) ds\right) _{h}=h%
\displaystyle\sum \limits_{j=0}^{i}f\left( t_{j}\right) ,\text{ }i=1,...,k,
\\ 
\left( \displaystyle\int \limits_{0}^{t_{0}}f\left( s\right) ds\right)
_{h}=0.%
\end{array}%
\right.   \label{3.5}
\end{equation}%
Let $s\geq 1$ be an integer. Below we need the following notations:%
\begin{equation}
\left. 
\begin{array}{c}
Q_{T}^{h}=\left\{ \left( \mathbf{x},t\right) :\mathbf{x}\in \Omega ,t\in
\left\{ 0=t_{0},t_{1},...,t_{k}=T\right\} \right\} , \\ 
H_{h}^{s}\left( Q_{T}^{h}\right) =\left\{ u^{h}\left( \mathbf{x},t\right)
:\left\Vert u^{h}\right\Vert _{H_{h}^{s}\left( Q_{T}^{h}\right) }^{2}=\left( %
\displaystyle\int \limits_{0}^{T}\left\Vert u^{h}\left( \mathbf{x},\alpha
\right) \right\Vert _{H^{s}\left( \Omega \right) }^{2}d\alpha \right)
_{h}\right\} <\infty , \\ 
L_{2,h}\left( \Gamma _{0,T}^{h}\right) =\left\{ u^{h}\left( \mathbf{x}%
,t\right) :\left\Vert u^{h}\right\Vert _{L_{2,h}\left( \Gamma
_{0,T}^{h}\right) }^{2}=\left( \displaystyle\int \limits_{0}^{T}\left\Vert
u^{h}\left( \mathbf{x},\alpha \right) \right\Vert _{L_{2}\left( \Gamma
_{0}\right) }^{2}d\alpha \right) _{h}\right\} <\infty , \\ 
C^{s}\left( \overline{Q}_{T_{1}}^{h}\right) =\left\{ u^{h}\left( \mathbf{x}%
,t\right) :\left\Vert u^{h}\right\Vert _{C^{s}\left( \overline{Q}%
_{T_{1}}^{h}\right) }=\max_{j\in \left[ 0,k\right] }\left\Vert u^{h}\left( 
\mathbf{x},t_{j}\right) \right\Vert _{C^{s}\left( \overline{\Omega }\right)
}<\infty \right\} .%
\end{array}%
\right.   \label{3.7}
\end{equation}

\section{Transformation Procedure for CIP1}

\label{sec:4}

We set in this section $m\left( \mathbf{x}\right) \equiv 1.$ Then the cone $%
\left\{ t=\left\vert \mathbf{x}\right\vert \right\} \subset \mathbb{R}%
^{3}\times \left( 0,\infty \right) $ is the characteristic cone for problem (%
\ref{2.1}), (\ref{2.2}). It follows from formulas (7.5)- (7.7) of Chapter 7
of \cite{LRS} that the solution $u\left( \mathbf{x},t\right) $ of problem (%
\ref{2.1}), (\ref{2.2}) has the form%
\begin{equation}
u\left( \mathbf{x},t\right) =\frac{\delta \left( t-\left\vert \mathbf{x}%
\right\vert \right) }{4\pi \left\vert \mathbf{x}\right\vert }+H\left(
t-\left\vert \mathbf{x}\right\vert \right) \widehat{u}\left( \mathbf{x}%
,t\right) ,  \label{4.1}
\end{equation}%
where $H\left( x\right) $ is the Heaviside function,%
\begin{equation*}
H\left( x\right) =\left\{ 
\begin{array}{c}
1\text{ if }x>0, \\ 
0,\text{ if }x<0.%
\end{array}%
\right. 
\end{equation*}%
Hence, 
\begin{equation}
u\left( \mathbf{x},t\right) =0\text{ for }t\in \left( -\infty ,\left\vert 
\mathbf{x}\right\vert \right) .  \label{4.2}
\end{equation}%
Furthermore, using (\ref{2.3}), (\ref{2.4}), (\ref{2.7})-(\ref{2.9}),
formulas (7.5)- (7.7) of Chapter 7 of \cite{LRS}, Lemma 1 and Corollary 1 of 
\cite{KR1} as well as the proof of that lemma, we obtain 
\begin{equation}
\widehat{u}\left( \mathbf{x},t\right) \geq 0\text{ for }t\geq \left\vert 
\mathbf{x}\right\vert ,  \label{4.3}
\end{equation}%
\begin{equation}
\partial _{t}^{k}\widehat{u},\widehat{u}_{x_{i}},\widehat{u}_{x_{i}x_{j}}\in
C\left( \overline{Q}_{T}\cap \left\{ t\geq \left\vert \mathbf{x}\right\vert
\right\} \right) ,\text{ }k=0,1,2;\text{ }i,j=1,2,3.  \label{4.4}
\end{equation}%
In addition, it was proven in \cite[\S 1 of Chapter 7]{LRS} that%
\begin{equation}
\lim_{t\rightarrow \left\vert \mathbf{x}\right\vert ^{+}}\widehat{u}\left( 
\mathbf{x},t\right) =\frac{1}{8\pi \left\vert \mathbf{x}\right\vert }%
\displaystyle\int \limits_{L\left( 0,\mathbf{x}\right) }q\left( \mathbf{\xi }%
\right) ds,  \label{4.5}
\end{equation}%
where $L\left( 0,\mathbf{x}\right) $ is the interval of the straight line
connecting points $\left\{ \mathbf{0}\right\} $ and $\left\{ \mathbf{x}%
\right\} $ and $ds$ is the element of the length of this line.

\subsection{Integral differential equation }

\label{sec:4.1}

Consider the function 
\begin{equation}
v\left( \mathbf{x},t\right) =\displaystyle\int\limits\limits_{-\infty
}^{t}u\left( \mathbf{x},\tau \right) d\tau .  \label{4.50}
\end{equation}%
Integrating equation (\ref{2.1}) with respect to $t$ and using (\ref{2.2}),
we obtain%
\begin{equation}
v_{tt}=\Delta v+q\left( \mathbf{x}\right) v+\delta \left( \mathbf{x}\right) 
\text{ in }D_{T}^{4}=\mathbb{R}^{3}\times \left( 0,T\right) .  \label{4.6}
\end{equation}%
In addition, (\ref{2.3}) and (\ref{4.1})-(\ref{4.5}) imply%
\begin{equation}
v\left( \mathbf{x},t\right) =\frac{H\left( t-\left\vert \mathbf{x}%
\right\vert \right) }{4\pi \left\vert \mathbf{x}\right\vert }+H\left(
t-\left\vert \mathbf{x}\right\vert \right) \overline{v}\left( \mathbf{x}%
,t\right) ,  \label{4.7}
\end{equation}%
\begin{equation}
\overline{v}\left( \mathbf{x},t\right) =\frac{t-\left\vert \mathbf{x}%
\right\vert }{2\left\vert \mathbf{x}\right\vert }\displaystyle%
\int\limits\limits_{L\left( 0,\mathbf{x}\right) }q\left( \mathbf{\xi }%
\right) ds+O\left( \left( t-\left\vert \mathbf{x}\right\vert \right)
^{2}\right) ,  \label{4.8}
\end{equation}%
\begin{equation}
\left\vert \frac{O\left( \left( t-\left\vert \mathbf{x}\right\vert \right)
^{2}\right) }{\left( t-\left\vert \mathbf{x}\right\vert \right) ^{2}}%
\right\vert \leq C\text{ in }\overline{Q}_{T}\cap \left\{ t\geq \left\vert 
\mathbf{x}\right\vert \right\} ,  \label{4.9}
\end{equation}%
\begin{equation}
\overline{v}\left( \mathbf{x},t\right) \geq 0\text{ in }\overline{Q}_{T}\cap
\left\{ t\geq \left\vert \mathbf{x}\right\vert \right\} ,  \label{4.10}
\end{equation}%
\begin{equation}
\left. 
\begin{array}{c}
\partial _{t}^{k}\overline{v},\partial _{t}^{s}\overline{v}_{x_{i}},\partial
_{t}^{s}\overline{v}_{x_{i}x_{j}}\in C\left( \overline{Q}_{T}\cap \left\{
t\geq \left\vert \mathbf{x}\right\vert \right\} \right) , \\ 
k=0,1,2;\text{ }s=0,1;\text{ }i,j=1,2,3,%
\end{array}%
\right.  \label{4.11}
\end{equation}%
where $C>0$ is a certain number.

Since $v\left( \mathbf{x},t\right) =0$ for $t<\left\vert \mathbf{x}%
\right\vert ,$ we change variables in equation (\ref{4.6}) via introduction
of a new function $w\left( \mathbf{x},t\right) ,$%
\begin{equation}
w\left( \mathbf{x},t\right) =v\left( \mathbf{x},t+\left\vert \mathbf{x}%
\right\vert \right) .  \label{4.12}
\end{equation}%
Let%
\begin{equation}
c_{1}=\max_{\mathbf{x}\in \overline{\Omega }}\left\vert \mathbf{x}%
\right\vert .  \label{4.120}
\end{equation}%
We assume below that 
\begin{equation}
T_{1}=T-c_{1}>0.  \label{4.13}
\end{equation}%
Using (\ref{2.3}), (\ref{2.4}), (\ref{2.6}) and (\ref{4.13}), denote%
\begin{equation}
\left. 
\begin{array}{c}
Q_{T_{1}}=\Omega \times \left( 0,T_{1}\right) ,\text{ }S_{T_{1}}=\partial
\Omega \times \left( 0,T_{1}\right) , \\ 
\Gamma _{i,T_{1}}=\Gamma _{i}\times \left( 0,T_{1}\right) ,\text{ }i=0,1,2.%
\end{array}%
\right.   \label{4.14}
\end{equation}%
By (\ref{2.4}), (\ref{4.50}) and (\ref{4.12})-(\ref{4.14}) equation (\ref%
{4.6}) becomes%
\begin{equation}
\Delta w-2\displaystyle\sum \limits_{i=1}^{3}\frac{x_{i}}{\left\vert \mathbf{%
x}\right\vert }w_{tx_{i}}-\frac{2}{\left\vert \mathbf{x}\right\vert }%
w_{t}=-q\left( \mathbf{x}\right) w\text{ in }Q_{T_{1}}.  \label{4.15}
\end{equation}%
Initial and boundary conditions for the function $w$ are:%
\begin{equation}
w\left( \mathbf{x},0\right) =\frac{1}{4\pi \left\vert \mathbf{x}\right\vert }%
,  \label{4.16}
\end{equation}%
\begin{equation}
w\mid _{S_{T_{1}}}=\widehat{g}_{0}\left( \mathbf{x},t\right) ,\text{ }%
w_{z}\mid _{\Gamma _{0,T_{1}}}=\widehat{g}_{1}\left( \mathbf{x},t\right) ,
\label{4.17}
\end{equation}%
where functions $\widehat{g}_{0}\left( \mathbf{x},t\right) $ and $\widehat{g}%
_{1}\left( \mathbf{x},t\right) $ are known and%
\begin{equation}
\text{ }\widehat{g}_{0}\in C^{2}\left( \overline{\Gamma _{i,T_{1}}}\right) ,%
\text{ }i=0,1,2;\text{ }\widehat{g}_{1}\in C^{2}\left( \overline{\Gamma
_{0,T_{1}}}\right) .  \label{4.170}
\end{equation}

Using (\ref{4.7}), (\ref{4.10}) and (\ref{4.12}), we obtain%
\begin{equation}
w\left( \mathbf{x},t\right) \geq \frac{1}{4\pi \left\vert \mathbf{x}%
\right\vert }\geq \frac{1}{4\pi c_{1}}\text{ in }Q_{T_{1}}.  \label{4.18}
\end{equation}%
In particular, (\ref{4.17}) and (\ref{4.18}) imply%
\begin{equation}
\widehat{g}_{0}\left( \mathbf{x},t\right) \geq \frac{1}{4\pi c_{1}}.
\label{4.180}
\end{equation}%
Hence, we can consider the function 
\begin{equation}
p\left( \mathbf{x},t\right) =\ln w\left( \mathbf{x},t\right) .  \label{4.19}
\end{equation}%
Substituting (\ref{4.19}) in (\ref{4.15}), we obtain for $\left( \mathbf{x}%
,t\right) \in Q_{T_{1}}$ 
\begin{equation}
\Delta p+\left( \nabla p\right) ^{2}-2\displaystyle\sum \limits_{i=1}^{3}%
\frac{x_{i}}{\left\vert \mathbf{x}\right\vert }p_{tx_{i}}-2p_{t}\displaystyle%
\sum \limits_{i=1}^{3}\frac{x_{i}}{\left\vert \mathbf{x}\right\vert }%
p_{x_{i}}-\frac{2p_{t}}{\left\vert \mathbf{x}\right\vert }=-q\left( \mathbf{x%
}\right) .  \label{4.20}
\end{equation}%
Differentiate both sides of (\ref{4.20}) with respect to $t$ and denote%
\begin{equation}
r\left( \mathbf{x},t\right) =\frac{\partial p}{\partial t}\left( \mathbf{x}%
,t\right) .  \label{4.21}
\end{equation}%
By (\ref{4.16}) and (\ref{4.21}) 
\begin{equation}
p\left( \mathbf{x},t\right) =\displaystyle\int \limits_{0}^{t}r\left( 
\mathbf{x},\tau \right) d\tau -\ln \left( 4\pi \left\vert \mathbf{x}%
\right\vert \right) \text{ in }Q_{T_{1}}.  \label{4.22}
\end{equation}%
Since $\partial _{t}q\left( \mathbf{x}\right) \equiv 0,$ then, using (\ref%
{4.22}), we obtain after the differentiation of (\ref{4.21})%
\begin{equation}
\Delta r+2\nabla r\left( \displaystyle\int \limits_{0}^{t}\nabla r\left( 
\mathbf{x},\tau \right) d\tau -\nabla \ln \left( \left\vert \mathbf{x}%
\right\vert \right) \right) -2\displaystyle\sum \limits_{i=1}^{3}\frac{x_{i}%
}{\left\vert \mathbf{x}\right\vert }r_{tx_{i}}-  \label{4.23}
\end{equation}%
\begin{equation*}
-2r\displaystyle\sum \limits_{i=1}^{3}\frac{x_{i}}{\left\vert \mathbf{x}%
\right\vert }r_{x_{i}}-2r_{t}\displaystyle\sum \limits_{i=1}^{3}\frac{x_{i}}{%
\left\vert \mathbf{x}\right\vert }\displaystyle\int \limits%
_{0}^{t}r_{x_{i}}\left( \mathbf{x},\tau \right) d\tau -2r_{t}\left( \frac{1}{%
\left\vert \mathbf{x}\right\vert }-1\right) =0\text{ in }Q_{T_{1}}.
\end{equation*}%
We have used here the following formula%
\begin{equation*}
-2r_{t}\displaystyle\sum \limits_{i=1}^{3}\frac{x_{i}}{\left\vert \mathbf{x}%
\right\vert }\left( -\partial _{x_{i}}\ln \left( \left\vert \mathbf{x}%
\right\vert \right) \right) =2r_{t}.
\end{equation*}%
In addition, using (\ref{4.17}), (\ref{4.19}) and (\ref{4.21}), we obtain
the following Dirichlet and Neumann boundary conditions for the function $%
r\left( \mathbf{x},t\right) $%
\begin{equation}
r\mid _{S_{T_{1}}}=\overline{g}_{0}\left( \mathbf{x},t\right) ,r_{z}\mid
_{\Gamma _{0,T_{1}}}=\overline{g}_{1}\left( \mathbf{x},t\right) ,
\label{4.24}
\end{equation}%
where 
\begin{equation}
\overline{g}_{0}\left( \mathbf{x},t\right) =\frac{\partial _{t}\widehat{g}%
_{0}}{\widehat{g}_{0}}\left( \mathbf{x},t\right) ,\text{ }\overline{g}%
_{1}\left( \mathbf{x},t\right) =\frac{\partial _{t}\widehat{g}_{1}}{\widehat{%
g}_{0}}\left( \mathbf{x},t\right) -\frac{\left( \partial _{t}\widehat{g}%
_{0}\right) \widehat{g}_{1}}{\widehat{g}_{0}^{2}}\left( \mathbf{x},t\right) .
\label{4.240}
\end{equation}%
As to the denominators in formulas (\ref{4.240}), we refer to (\ref{4.180}).

\subsection{Finite differences with respect to $t$ for problem (\protect\ref%
{4.23})-(\protect\ref{4.240})}

\label{sec:4.2}

We now need to work with problem (\ref{4.23}), (\ref{4.24}). Replace in (\ref%
{3.1}) $T$ with $T_{1}$ and use formulas (\ref{3.2})-(\ref{3.5}). We rewrite
equation (\ref{4.23}) in terms of finite differences with respect to $t$ and
come up this way with a system of coupled elliptic equations with respect to
the functions $r\left( \mathbf{x},t_{j}\right) .$ As it is always the case
in the theory of finite differences (see, e.g. \cite{Sam}), we assume that
this system holds. Our system is:%
\begin{equation*}
L_{j}\left( r^{h}\right) \left( \mathbf{x},t_{j}\right) =\Delta r\left( 
\mathbf{x},t_{j}\right) +2\nabla r\left( \mathbf{x},t_{j}\right) \left[
\left( \displaystyle\int \limits_{0}^{t_{j}}\nabla r^{h}\left( \mathbf{x}%
,\tau \right) d\tau \right) _{h}-\nabla \ln \left( \left\vert \mathbf{x}%
\right\vert \right) \right] -
\end{equation*}%
\begin{equation}
-2\displaystyle\sum \limits_{i=1}^{3}\frac{x_{i}}{\left\vert \mathbf{x}%
\right\vert }\partial _{t}^{h}r_{x_{i}}\left( \mathbf{x},t_{j}\right) -
\label{4.25}
\end{equation}%
\begin{equation*}
-2r\left( \mathbf{x},t_{j}\right) \sum_{i=1}^{3}\frac{x_{i}}{\left\vert 
\mathbf{x}\right\vert }r_{x_{i}}\left( \mathbf{x},t_{j}\right) -2\partial
_{t}^{h}r\left( \mathbf{x},t_{j}\right) \sum_{i=1}^{3}\frac{x_{i}}{%
\left\vert \mathbf{x}\right\vert }\left(
\dint\limits_{0}^{t_{j}}r_{x_{i}}^{h}\left( \mathbf{x},\tau \right) d\tau
\right) _{h}-
\end{equation*}%
\begin{equation*}
-2\partial _{t}^{h}r\left( \mathbf{x},t_{j}\right) \left( \frac{1}{%
\left\vert \mathbf{x}\right\vert }-1\right) =0,\text{ }\mathbf{x}\in \Omega ,%
\text{ }j=0,...,k.
\end{equation*}%
In addition, using (\ref{2.5}), (\ref{2.6}), (\ref{3.4}) and (\ref{4.24}) we
obtain the following boundary conditions for the semi-discrete function $%
r^{h}\left( \mathbf{x},t\right) $:%
\begin{equation}
r^{h}\left( \mathbf{x},t\right) \mid _{S_{T_{1}}^{h}}=\overline{g}%
_{0}^{h}\left( \mathbf{x},t\right) ,\text{ }r_{z}^{h}\left( \mathbf{x}%
,t\right) \mid _{\Gamma _{0,T_{1}}^{h}}=\overline{g}_{1}^{h}\left( \mathbf{x}%
,t\right) ,  \label{4.26}
\end{equation}%
\begin{equation}
\left. 
\begin{array}{c}
S_{T_{1}}^{h}=\left\{ \mathbf{x}\in \partial \Omega
,t=0=t_{0},t_{1},...,t_{k}\right\} ,\text{ } \\ 
\Gamma _{0,T_{1}}^{h}=\left\{ \mathbf{x}\in \Gamma
_{0},t=0=t_{0},t_{1},...,t_{k}\right\} .%
\end{array}%
\right.   \label{4.026}
\end{equation}%
By (\ref{4.170}) 
\begin{equation}
\text{ }\overline{g}_{0}\left( \mathbf{x},t_{j}\right) \in C^{2}\left(
\Gamma _{i}\right) ,\text{ }\overline{g}_{1}\left( \mathbf{x},t_{j}\right)
\in C^{2}\left( \Gamma _{0}\right) ,\text{ }i=0,1,2;\text{ }j=0,...,k.
\label{4.260}
\end{equation}%
Suppose that problem (\ref{4.25})-(\ref{4.260}) is solved. Then, using (\ref%
{4.20}) and (\ref{4.22}), we set%
\begin{equation}
p\left( \mathbf{x},t_{j}\right) =\left( \displaystyle\int \limits%
_{0}^{t_{j}}r^{h}\left( \mathbf{x},\tau \right) d\tau \right) _{h}-\ln
\left( 4\pi \left\vert \mathbf{x}\right\vert \right) ,\text{ }j=0,...,k,
\label{4.27}
\end{equation}%
\begin{equation*}
q_{h}\left( \mathbf{x}\right) =-\frac{1}{T_{1}}\left( \displaystyle\int %
\limits_{0}^{T_{1}}\left( \Delta p^{h}+\left( \nabla p^{h}\right) ^{2}-2%
\displaystyle\sum \limits_{i=1}^{3}\left( x_{i}/\left\vert \mathbf{x}%
\right\vert \right) \partial _{t}^{h}p_{x_{i}}^{h}\right) \left( \mathbf{x}%
,\tau \right) d\tau \right) _{h}+
\end{equation*}%
\begin{equation}
+\frac{1}{T_{1}}\left( \displaystyle\int \limits_{0}^{T_{1}}\left( 2\partial
_{t}^{h}p\displaystyle\sum \limits_{i=1}^{3}\left( x_{i}/\left\vert \mathbf{x%
}\right\vert \right) p_{x_{i}}^{h}+2\frac{\partial _{t}^{h}p}{\left\vert 
\mathbf{x}\right\vert }\right) \left( \mathbf{x},\tau \right) d\tau \right)
_{h}.  \label{4.28}
\end{equation}

\textbf{Remark 4.1. }\emph{Thus, we consider the function }$q_{h}\left( 
\mathbf{x}\right) $\emph{\ in (\ref{4.28}) as the approximation of the
unknown coefficient }$q\left( \mathbf{x}\right) $\emph{\ generated by our
above finite difference approximation. }

The transformation procedure for CIP1 is completed.

\section{Transformation Procedure for CIP2 }

\label{sec:5}

\subsection{A system of two coupled PDEs}

\label{sec:5.1}

We now consider Cauchy problem (\ref{2.1}), (\ref{2.2}) in the case%
\begin{equation}
q\left( \mathbf{x}\right) \equiv 0.  \label{5.1}
\end{equation}%
The function $m\left( \mathbf{x}\right) $ generates the following conformal
Riemannian metric%
\begin{equation*}
d\tau =\sqrt{m\left( \mathbf{x}\right) }\left\vert d\mathbf{x}\right\vert ,%
\text{ }\left\vert d\mathbf{x}\right\vert =\sqrt{\left( dx\right)
^{2}+\left( dy\right) ^{2}+\left( dz\right) ^{2}}.\text{ }
\end{equation*}%
This metric generates geodesic lines $\Gamma \left( \mathbf{x},\mathbf{x}%
_{0}\right) $ connecting points $\mathbf{x\in }\mathbb{R}^{3}$ and $\mathbf{x%
}_{0}\mathbf{\in }\mathbb{R}^{3}.$ Let $\tau \left( \mathbf{x},\mathbf{x}%
_{0}\right) $ be the travel time along $\Gamma \left( \mathbf{x},\mathbf{x}%
_{0}\right) $ from the point $\mathbf{x}_{0}$ to the point $\mathbf{x}.$
Then this function satisfies the eikonal equation \cite[\S 1 of Chapter 3]%
{Rom},%
\begin{equation}
\left. 
\begin{array}{c}
\left\vert \nabla _{\mathbf{x}}\tau \left( \mathbf{x},\mathbf{x}_{0}\right)
\right\vert ^{2}=m\left( \mathbf{x}\right) , \\ 
\tau \left( \mathbf{x},\mathbf{x}_{0}\right) =O\left( \left\vert \mathbf{x}-%
\mathbf{x}_{0}\right\vert \right) \text{ as }\left\vert \mathbf{x}-\mathbf{x}%
_{0}\right\vert \rightarrow 0.%
\end{array}%
\right.   \label{5.2}
\end{equation}%
Below we use the following assumption without further comments:

\textbf{Assumption} (regularity of geodesic lines). For any two points $%
\mathbf{x\in }\mathbb{R}^{3}$ and $\mathbf{x}_{0}\mathbf{\in }\mathbb{R}^{3}$
there exists unique geodesic line $\Gamma \left( \mathbf{x},\mathbf{x}%
_{0}\right) $ connecting them.

Since we use only a single point source $\left\{ \mathbf{x}_{0}=0\right\} ,$
then we denote for brevity 
\begin{equation}
\tau \left( \mathbf{x}\right) =\tau \left( \mathbf{x},0\right) .
\label{5.02}
\end{equation}%
Note that due to (\ref{5.2}) and (\ref{5.02}) there exists a number $%
\overline{C}>0$ such that%
\begin{equation}
\left\Vert m\right\Vert _{C^{2}\left( \overline{\Omega }\right) }\leq 
\overline{C}\left\Vert \tau \left( \mathbf{x}\right) \right\Vert
_{C^{3}\left( \overline{\Omega }\right) }.  \label{5.03}
\end{equation}

Theorem 5.1 follows immediately from Theorem 4.1 of \cite{Rom}, Theorem 3.1
of \cite{KR} and (\ref{5.03}).

\textbf{Theorem 5.1.} \emph{Assume that conditions (\ref{2.10})-(\ref{2.12})
and (\ref{5.1}) hold. Then there exists unique solution of problem (\ref{2.1}%
), (\ref{2.2}) of the following form:}%
\begin{equation}
u\left( \mathbf{x},t\right) =A\left( \mathbf{x}\right) \delta \left( t-\tau
\left( \mathbf{x}\right) \right) +H\left( t-\tau \left( \mathbf{x}\right)
\right) \widehat{u}\left( \mathbf{x},t\right) ,~(\mathbf{x},t)\in D_{T}^{4},
\label{5.3}
\end{equation}%
\begin{equation}
\tau \left( \mathbf{x}\right) \in C^{16}\left( \mathbb{R}^{3}\setminus
0\right) ,  \label{5.4}
\end{equation}%
\begin{equation}
A\left( \mathbf{x}\right) \geq A_{0}>0\text{ in }\overline{\Omega },~\text{ }%
A\in C^{14}\left( \mathbb{R}^{3}\setminus 0\right) ,  \label{5.5}
\end{equation}%
\begin{equation}
\widehat{u}\in C^{3}\left( \left\{ (\mathbf{x},t):\mathbf{x}\in \mathbb{R}%
^{3}\diagdown \left\{ 0\right\} ,t\in \left[ \tau \left( \mathbf{x}\right) ,T%
\right] \text{ }\right\} \right) .  \label{5.6}
\end{equation}%
\emph{Equation }$t=\tau (\mathbf{x})$\emph{\ defines the characteristic cone
in }$D_{T}^{4}$\emph{. In (\ref{5.5}) }%
\begin{equation}
A_{0}=A_{0}\left( \left\Vert \tau (\mathbf{x})\right\Vert _{C^{3}\left( 
\overline{\Omega }\right) },\Omega \right) >0  \label{5.60}
\end{equation}%
\emph{is a number depending only on listed parameters.}

Below we use some notations of section 4 for different functions. It will be
always clear from the context which exactly functions we work with. Define
the function $v\left( \mathbf{x},t\right) $ as%
\begin{equation*}
v\left( \mathbf{x},t\right) =\dint\limits_{-\infty }^{t}\left(
\dint\limits_{-\infty }^{s}u\left( \mathbf{x},\rho \right) d\rho \right)
ds=\dint\limits_{-\infty }^{t}\left( t-\rho \right) u\left( \mathbf{x},\rho
\right) d\rho .
\end{equation*}%
Then (\ref{5.3}) implies%
\begin{equation}
m\left( \mathbf{x}\right) v_{tt}=\Delta v\text{ in }Q_{T},  \label{5.7}
\end{equation}%
\begin{equation}
v\left( \mathbf{x},t\right) =H\left( t-\tau \left( \mathbf{x}\right) \right) %
\left[ A\left( \mathbf{x}\right) \left( t-\tau \left( \mathbf{x}\right)
\right) +\left( t-\tau \left( \mathbf{x}\right) \right) ^{2}\widehat{v}%
\left( \mathbf{x},t\right) \right] ,~\text{in }Q_{T},  \label{5.8}
\end{equation}%
\begin{equation}
\widehat{v}\mathbf{\in }C^{3}\left( \overline{Q}_{T}\right) .  \label{5.9}
\end{equation}%
It follows from (\ref{5.5}) and (\ref{5.8}) that%
\begin{equation}
\lim_{t\rightarrow +0}v\left( \mathbf{x},t+\tau \left( \mathbf{x}\right)
\right) =0,  \label{5.10}
\end{equation}%
\begin{equation}
\lim_{t\rightarrow +0}v_{t}\left( \mathbf{x},t+\tau \left( \mathbf{x}\right)
\right) =A\left( \mathbf{x}\right) \geq A_{0}>0\text{ in }\overline{\Omega }.
\label{5.11}
\end{equation}%
It was proven in \cite[formula 4.7]{KR} that%
\begin{equation}
\max_{\overline{\Omega }}\tau \left( \mathbf{x}\right) \leq
c_{2}=c_{2}\left( \Omega \right) ,  \label{5.12}
\end{equation}%
where the number $c_{2}\left( \Omega \right) >0$ depends only on the domain $%
\Omega .$ We assume below that 
\begin{equation}
T_{2}=T-c_{2}>0.  \label{5.13}
\end{equation}%
Using (\ref{2.3}), (\ref{5.12}) and (\ref{5.13}), denote similarly with (\ref%
{4.14}) 
\begin{equation}
\left. 
\begin{array}{c}
Q_{T_{2}}=\Omega \times \left( 0,T_{2}\right) ,\text{ }S_{T_{2}}=\partial
\Omega \times \left( 0,T_{2}\right) , \\ 
\Gamma _{i,T_{2}}=\Gamma _{i}\times \left( 0,T_{2}\right) ,\text{ }i=0,1,2.%
\end{array}%
\right.   \label{5.14}
\end{equation}

Consider the function 
\begin{equation}
p\left( \mathbf{x},t\right) =v\left( \mathbf{x},t+\tau \left( \mathbf{x}%
\right) \right) .  \label{5.15}
\end{equation}%
Then, using (\ref{5.2}) and (\ref{5.7}), we obtain%
\begin{equation}
\Delta p-2\nabla _{\mathbf{x}}p_{t}\cdot \nabla \tau -p_{t}\Delta \tau =0%
\text{ in }Q_{T_{2}}.  \label{5.16}
\end{equation}%
It follows from (\ref{5.6}), (\ref{5.8}), (\ref{5.10}), (\ref{5.11}) and (%
\ref{5.15}) that 
\begin{equation}
p\left( \mathbf{x},t\right) =A\left( \mathbf{x}\right) t+\widehat{p}\left( 
\mathbf{x},t\right) t^{2},\text{ }\widehat{p}\in C^{3}\left( \overline{Q}%
_{T_{2}}\right) .  \label{5.17}
\end{equation}%
Introduce the function $w\left( \mathbf{x},t\right) ,$%
\begin{equation}
w\left( \mathbf{x},t\right) =p_{t}\left( \mathbf{x},t\right) .  \label{5.18}
\end{equation}%
Then (\ref{5.16}) leads to%
\begin{equation}
\Delta w-2\nabla _{\mathbf{x}}w_{t}\cdot \nabla \tau -w_{t}\Delta \tau =0%
\text{ in }Q_{T_{2}}.  \label{5.19}
\end{equation}%
In addition, using setting in (\ref{5.16}) $t=0$ and using  (\ref{5.11}), (%
\ref{5.17}) and (\ref{5.18}), we obtain%
\begin{equation}
\Delta \tau +2\frac{\nabla w\left( \mathbf{x},0\right) }{w\left( \mathbf{x}%
,0\right) }\nabla \tau =0\text{ in }\Omega .  \label{5.20}
\end{equation}%
Combining (\ref{5.19}) with (\ref{5.20}), we finally obtain the following
system of coupled nonlinear PDEs:%
\begin{equation}
\Delta w-2\nabla _{\mathbf{x}}w_{t}\cdot \nabla \tau +2w_{t}\frac{\nabla
w\left( \mathbf{x},0\right) }{w\left( \mathbf{x},0\right) }\nabla \tau =0%
\text{ in }Q_{T_{2}},  \label{5.21}
\end{equation}%
\begin{equation}
\Delta \tau +2\frac{\nabla w\left( \mathbf{x},0\right) }{w\left( \mathbf{x}%
,0\right) }\nabla \tau =0\text{ in }\Omega .  \label{5.210}
\end{equation}%
System (\ref{5.21}), (\ref{5.210}) was first derived in \cite{KR}. Dirichlet
and Neumann boundary conditions for functions $w\left( \mathbf{x},t\right) $
and $\tau \left( \mathbf{x}\right) $ are:%
\begin{equation}
w\mid _{S_{T_{2}}}=s_{0}\left( \mathbf{x},t\right) ,\text{ }w_{z}\mid
_{\Gamma _{0,T_{2}}}=s_{1}\left( \mathbf{x},t\right) ,  \label{5.22}
\end{equation}%
\begin{equation}
\tau \mid _{\partial \Omega }=\tau _{0}\left( \mathbf{x}\right) ,\text{ }%
\tau _{z}\mid _{\Gamma _{0}}=\tau _{1}\left( \mathbf{x}\right) ,
\label{5.23}
\end{equation}%
\begin{equation}
\left. 
\begin{array}{c}
\text{ }s_{0}\in C^{2}\left( \overline{\Gamma _{i,T_{2}}}\right) ,s_{1}\in
C^{2}\left( \overline{\Gamma _{0,T_{2}}}\right) , \\ 
\tau _{0}\in C^{2}\left( \overline{\Gamma _{i}}\right) ,\tau _{1}\in
C^{2}\left( \overline{\Gamma _{0}}\right) , \\ 
\text{ }i=0,1,2.%
\end{array}%
\right.   \label{5.24}
\end{equation}%
Functions  $s_{0},s_{1},\tau _{0},\tau _{1}$ in (\ref{5.22})-(\ref{5.24})
are known since they are derived from given functions (\ref{2.13}) and (\ref%
{2.14}), also, see (\ref{5.14}). Furthermore, using (\ref{5.5}), (\ref{5.17}%
) and (\ref{5.18}), we obtain%
\begin{equation}
0<\frac{1}{w\left( \mathbf{x},0\right) }=\frac{1}{A\left( \mathbf{x}\right) }%
\leq \frac{1}{A_{0}}.  \label{5.25}
\end{equation}

\subsection{Finite differences with respect to $t$ for problem (\protect\ref%
{5.21})-(\protect\ref{5.24})}

\label{sec:5.2}

Replace in (\ref{3.1}) $T$ with $T_{2}$ and, just as in section 4.2, use
formulas (\ref{3.2})-(\ref{3.4}). Rewrite the terms $\nabla _{\mathbf{x}%
}w_{t},w_{t}$ in equation (\ref{5.21}) via finite differences with respect
to $t$. Then equations (\ref{5.21}) and (\ref{5.210}) are turned in a system
of coupled nonlinear elliptic equations with respect to functions $w\left( 
\mathbf{x},t_{j}\right) $ and $\tau \left( \mathbf{x}\right) ,$%
\begin{equation}
\left. 
\begin{array}{c}
L_{j}^{\left( 1\right) }\left( w^{h},\tau \right) \left( \mathbf{x}%
,t_{j}\right) =\Delta w\left( \mathbf{x},t_{j}\right) -2\nabla \partial
_{t}^{h}w\left( \mathbf{x},t_{j}\right) \cdot \nabla \tau + \\ 
+2\partial _{t}^{h}w\left( \mathbf{x},t_{j}\right) \left( \nabla w\left( 
\mathbf{x},t_{0}\right) /w\left( \mathbf{x},t_{0}\right) \right) \nabla \tau
=0, \\ 
\text{ }\mathbf{x}\in \Omega ,\text{ }j=0,...,k,%
\end{array}%
\right.   \label{5.26}
\end{equation}%
\begin{equation}
L^{\left( 2\right) }\left( w\left( \mathbf{x},t_{0}\right) ,\tau \right) 
\mathbf{=}\Delta \tau +2\frac{\nabla w\left( \mathbf{x},t_{0}\right) }{%
w\left( \mathbf{x},t_{0}\right) }\nabla \tau =0,\text{ }\mathbf{x}\in \Omega
.  \label{5.27}
\end{equation}%
Dirichlet and Neumann boundary conditions for this elliptic system are
derived from (\ref{5.22}) and (\ref{5.23})%
\begin{equation}
w^{h}\mid _{S_{T_{2}}^{h}}=s_{0}^{h}\left( \mathbf{x},t\right) ,\text{ }%
w_{z}^{h}\mid _{\Gamma _{0,T_{2}}^{h}}=s_{1}^{h}\left( \mathbf{x},t\right) ,
\label{5.28}
\end{equation}%
\begin{equation}
\tau \mid _{\partial \Omega }=\tau _{0}\left( \mathbf{x}\right) ,\text{ }%
\tau _{z}\mid _{\Gamma _{0}}=\tau _{1}\left( \mathbf{x}\right) ,
\label{5.29}
\end{equation}%
\begin{equation}
\left. 
\begin{array}{c}
S_{T_{2}}^{h}=\left\{ \mathbf{x}\in \partial \Omega
,t=0=t_{0},t_{1},...,t_{k}\right\} ,\text{ } \\ 
\Gamma _{0,T_{2}}^{h}=\left\{ \mathbf{x}\in \Gamma
_{0},t=0=t_{0},t_{1},...,t_{k}\right\} .%
\end{array}%
\right.   \label{5.30}
\end{equation}%
By (\ref{5.24})%
\begin{equation}
\left. 
\begin{array}{c}
s_{0}\left( \mathbf{x},t_{j}\right) ,\tau _{0}\in C^{2}\left( \overline{%
\Gamma _{i}}\right) ;\text{ }s_{1}\left( \mathbf{x},t_{j}\right) ,\tau
_{1}\in C^{2}\left( \overline{\Gamma _{0}}\right) , \\ 
\text{ }j=0,...,k;\text{ }i=0,1,2.%
\end{array}%
\right.   \label{5.31}
\end{equation}%
Assume that problem (\ref{5.26})-(\ref{5.31}) is solved, i.e. assume that
the vector function $\left( w^{h}\left( \mathbf{x},t\right) ,\tau \left( 
\mathbf{x}\right) \right) $ is found. Then we find the target coefficient $%
m\left( \mathbf{x}\right) $ via (\ref{5.2}) and (\ref{5.02}),%
\begin{equation}
m\left( \mathbf{x}\right) =\left\vert \nabla \tau \left( \mathbf{x}\right)
\right\vert ^{2}.  \label{5.32}
\end{equation}%
The transformation procedure for CIP2 is completed.

\section{Formulations of Some Theorems }

\label{sec:6}

Proofs of Theorems 6.2 and 6.3 of this section can be found in section 7.
Also, see Remark 6.2 in subsection 6.1.

\subsection{A Carleman estimate for the Laplace operator }

\label{sec:6.1}

First, we formulate a Carleman estimate for the Laplace operator in the
rectangular prism $\Omega $ defined in (\ref{2.3}). While usually Carleman
estimates are formulated for the case of zero boundary conditions, we
formulate Theorem 6.1 for the case of non-zero Dirichlet and Neumann
boundary conditions at $\Gamma _{0}$ since we need to estimate solutions of
our CIPs via non-zero boundary data.

Let 
\begin{equation}
d>B  \label{6.1}
\end{equation}%
be a number and $\lambda \geq 1$ be a parameter. We introduce the Carleman
Weight Function (CWF) for the Laplace operator in the domain $\Omega $ as%
\begin{equation}
\varphi _{\lambda }\left( z\right) =e^{\lambda \left( z-d\right) ^{2}}.
\label{6.2}
\end{equation}%
Hence, by (\ref{2.4}), (\ref{6.1}) and (\ref{6.2}) the function $\varphi
_{\lambda }\left( z\right) $ is monotonically decreasing for $z\in \left(
a,B\right) ,$ and also%
\begin{equation}
\max_{\left[ a,B\right] }\varphi _{\lambda }\left( z\right) =\varphi
_{\lambda }\left( a\right) =e^{\lambda \left( a-d\right) ^{2}},\text{ }\min_{%
\left[ a,B\right] }\varphi _{\lambda }\left( z\right) =\varphi _{\lambda
}\left( B\right) =e^{\lambda \left( B-d\right) ^{2}}.  \label{6.3}
\end{equation}%
Denote%
\begin{equation}
H_{0}^{2}\left( \Omega \right) =\left\{ u\in H^{2}\left( \Omega \right)
:u\mid _{\partial \Omega \diagdown \Gamma _{0}}=0\right\} .  \label{6.4}
\end{equation}

\textbf{Theorem 6.1} (Carleman estimate with non-zero boundary terms at $%
\Gamma _{0}$). \emph{There exists a sufficiently large number }$\lambda
_{0}=\lambda _{0}\left( d,\Omega \right) \geq 1$\emph{\ and a number }$%
C=C\left( d,\Omega \right) >0$\emph{\ such that the following Carleman
estimate for the Laplace operator holds}%
\begin{equation}
\left. 
\begin{array}{c}
\displaystyle\int \limits_{\Omega }\left( \Delta u\right) ^{2}\varphi
_{\lambda }^{2}\left( z\right) d\mathbf{x}\geq C\lambda \displaystyle\int %
\limits_{\Omega }\left\vert \nabla u\right\vert ^{2}\varphi _{\lambda
}^{2}\left( z\right) d\mathbf{x+}C\lambda ^{3}\displaystyle\int \limits%
_{\Omega }u^{2}\varphi _{\lambda }^{2}\left( z\right) d\mathbf{x-} \\ 
-C\lambda ^{3}e^{2\lambda \left( a-d\right) ^{2}}\left( \left\Vert u\left(
x,y,0\right) \right\Vert _{L_{2}\left( \Gamma _{0}\right) }^{2}+\left\Vert
u_{z}\left( x,y,0\right) \right\Vert _{L_{2}\left( \Gamma _{0}\right)
}^{2}\right) , \\ 
\forall u\in H_{0}^{2}\left( \Omega \right) ,\text{ }\forall \lambda \geq
\lambda _{0}.%
\end{array}%
\right.   \label{6.40}
\end{equation}%
\emph{Both numbers }$\lambda _{0}$\emph{\ and }$C$\emph{\ depend only on }$d$%
\emph{\ and }$\Omega .$

\textbf{Remark 6.1. }\emph{A similar Carleman estimate is also valid for a
more general elliptic operator of the second order with variable
coefficients and for a more general domain \ }$\Omega ,$\emph{\ see, e.g. 
\cite[Theorem 2.4.1]{KL}, \cite[\S 1 of Chapter 4]{LRS}. The CWF in this
case depends on two large parameters rather than on just one in (\ref{6.2}).
However, the dependence on two large parameters causes a significant
complication is the numerical implementation, to which our section 10 is
devoted. The latter is the reason of our use of a simplified Carleman
estimate.}

\textbf{Remark 6.2. }\emph{A close analog of Theorem 6.1 was proven in \cite[%
Theorem 2.4.1]{Khoa} and \cite[Theorem 10.3.1]{KL} under the condition that
boundary terms in the second line of (\ref{6.40}) equal zero. Incorporation
of non-zero boundary terms requires only a slight change of proofs of these
references. Hence, we omit the proof of Theorem 6.1 \cite{Khoa,KL}.}

\subsection{Coefficient Inverse Problem 1}

\label{sec:6.2}

Let $M_{1}>0$ be a number. Consider the set $D\left( M_{1}\right) $ of
vector functions $r^{h}\left( \mathbf{x},t\right) ,$%
\begin{equation}
D\left( M_{1}\right) =\left\{ 
\begin{array}{c}
r^{h}\left( \mathbf{x},t\right) :r\left( \mathbf{x},t_{j}\right) \in
C^{2}\left( \overline{\Omega }\right) ,\left\Vert r\left( \mathbf{x}%
,t_{j}\right) \right\Vert _{C^{1}\left( \overline{\Omega }\right) }\leq
M_{1}, \\ 
\text{ }j=0,...,k%
\end{array}%
\right\} .  \label{6.5}
\end{equation}%
Let $h_{0}\in \left( 0,1\right) $ be an arbitrary number. We set for our
grid step size $h$ of the above finite differences with respect to $t$%
\begin{equation}
h\geq h_{0}.  \label{6.6}
\end{equation}

\textbf{Theorem 6.2 }(Lipschitz stability and uniqueness for CIP1). \emph{%
Assume that there exist two vector functions }$r_{1}^{h},r_{2}^{h}\in
D\left( M_{1}\right) $\emph{\ satisfying the system of equations (\ref{4.25}%
) with two different pairs of boundary conditions (\ref{4.26}),}%
\begin{equation}
r_{i}^{h}\left( \mathbf{x},t\right) \mid _{S_{T_{1}}^{h}}=\overline{g}%
_{0,i}^{h}\left( \mathbf{x},t\right) ,\text{ }\partial _{z}r_{i}^{h}\left( 
\mathbf{x},t\right) \mid _{\Gamma _{0,T_{1}}^{h}}=\overline{g}%
_{1,i}^{h}\left( \mathbf{x},t\right) ,\text{ }i=1,2.  \label{6.7}
\end{equation}%
\emph{Keeping in mind (\ref{4.026}) and (\ref{6.4}), let} 
\begin{equation}
\overline{g}_{0,1}^{h}\left( \mathbf{x},t\right) =\overline{g}%
_{0,2}^{h}\left( \mathbf{x},t\right) \text{ for }\left( \mathbf{x},t\right)
\in S_{T_{1}}^{h}\diagdown \Gamma _{0,T_{1}}^{h}.  \label{6.8}
\end{equation}%
\emph{Let }$c_{1}$\emph{\ be the number defined in (\ref{4.120}) and let (%
\ref{6.6}) be valid. Then there exists a number }$C_{1}=C_{1}\left(
c_{1},h_{0},M_{1},T_{1},\Omega \right) >0$\emph{\ such that the following
Lipschitz stability estimate holds}%
\begin{equation}
\left\Vert r_{1}^{h}-r_{2}^{h}\right\Vert _{H_{h}^{1}\left(
Q_{T_{1}}^{h}\right) }\leq C_{1}\left( \left\Vert \overline{g}_{0,1}^{h}-%
\overline{g}_{0,2}^{h}\right\Vert _{L_{2,h}\left( \Gamma
_{0,T_{1}}^{h}\right) }+\left\Vert \overline{g}_{1,1}^{h}-\overline{g}%
_{1,2}^{h}\right\Vert _{L_{2,h}\left( \Gamma _{0,T_{1}}^{h}\right) }\right) ,
\label{6.9}
\end{equation}%
\emph{where }$H_{h}^{1}\left( Q_{T_{1}}^{h}\right) $\emph{\ and }$%
L_{2,h}\left( \Gamma _{0,T_{1}}^{h}\right) $\emph{\ are spaces defined in (%
\ref{3.7}), where }$T$\emph{\ is replaced with }$T_{1}$\emph{.} \emph{%
Furthermore, let the functions }$q_{h,1}\left( \mathbf{x}\right) $\emph{\
and }$q_{h,2}\left( \mathbf{x}\right) $\emph{\ be the approximations (\ref%
{4.28}) of two functions }$q_{1}\left( \mathbf{x}\right) $\emph{\ and }$%
q_{2}\left( \mathbf{x}\right) $\emph{\ (see Remark 4.1), in which case
vector functions }$r_{1}^{h}$\emph{\ and }$r_{2}^{h}$\emph{\ respectively
are used in (\ref{4.28}). Then the following Lipschitz stability estimate
holds} 
\begin{equation}
\left\Vert q_{h,1}-q_{h,2}\right\Vert _{L_{2}\left( \Omega \right) }\leq
C_{1}\left( \left\Vert \overline{g}_{0,1}^{h}-\overline{g}%
_{0,2}^{h}\right\Vert _{L_{2,h}\left( \Gamma _{0,T_{1}}^{h}\right)
}+\left\Vert \overline{g}_{1,1}^{h}-\overline{g}_{1,2}^{h}\right\Vert
_{L_{2,h}\left( \Gamma _{0,T_{1}}^{h}\right) }\right) .  \label{6.10}
\end{equation}%
\emph{In particular, if }$\overline{g}_{0,1}^{h}=\overline{g}_{0,2}^{h}$%
\emph{\ in }$S_{T_{1}}^{h}$\emph{\ and }$\overline{g}_{1,1}^{h}=\overline{g}%
_{1,2}^{h}$\emph{\ in }$\Gamma _{0,T_{1}}^{h},$\emph{\ then} $%
r_{1}^{h}\left( \mathbf{x},t\right) =r_{2}^{h}\left( \mathbf{x},t\right) $
in $Q_{T_{1}}^{h}$ \emph{and} $q_{h,1}\left( \mathbf{x}\right)
=q_{h,2}\left( \mathbf{x}\right) $ \emph{in} $\Omega $\emph{, i.e.
uniqueness holds. The above number }$C_{1}=C_{1}\left(
c_{1},h_{0},M_{1},T_{1},\Omega \right) >0$\emph{\ depends only on listed
parameters.}

Below $C_{1}>0$ denotes different numbers depending only on parameters
indicated in Theorem 6.2.

\subsection{Coefficient Inverse Problem 2}

\label{sec:6.3}

Let $M_{2}>0$ and $M_{3}>0$ be two numbers. Consider the set $D\left(
M_{2},M_{3}\right) $ of vector functions $\left( w^{h}\left( \mathbf{x}%
,t\right) ,\tau \left( \mathbf{x}\right) \right) ,$%
\begin{equation}
D\left( M_{2},M_{3}\right) =\left\{ 
\begin{array}{c}
\left( w^{h}\left( \mathbf{x},t\right) ,\tau \left( \mathbf{x}\right)
\right) :w\left( \mathbf{x},t_{j}\right) \in C^{2}\left( \overline{\Omega }%
\right) ,\text{ }\tau \left( \mathbf{x}\right) \in C^{3}\left( \overline{%
\Omega }\right) , \\ 
\left\Vert w\left( \mathbf{x},t_{j}\right) \right\Vert _{C^{1}\left( 
\overline{\Omega }\right) }\leq M_{2},\text{ }\left\Vert \tau \left( \mathbf{%
x}\right) \right\Vert _{C^{3}\left( \overline{\Omega }\right) }\leq M_{3},
\\ 
\text{ }j=0,...,k, \\ 
w\left( \mathbf{x},t_{0}\right) \geq A_{0}\left( M_{3},\Omega \right) >0%
\end{array}%
\right\} ,  \label{6.11}
\end{equation}%
where $A_{0}\left( M_{3},\Omega \right) $ is the number of Theorem 5.1. As
to the last line of (\ref{6.11}), we refer to (\ref{5.25}).

\textbf{Theorem 6.3 }(Lipschitz stability and uniqueness for CIP2). \emph{%
Let (\ref{6.6}) holds. Assume that there exist two pairs of vector functions 
}%
\begin{equation}
\left( w_{1}^{h},\tau _{1}\right) ,\left( w_{2}^{h},\tau _{2}\right) \in
D\left( M_{2},M_{3}\right)   \label{6.110}
\end{equation}%
\emph{satisfying the system of equations (\ref{5.26}) and (\ref{5.27}) with
two different sets of boundary conditions (\ref{5.28}), (\ref{5.29})
satisfying (\ref{5.31}), }%
\begin{equation}
w_{i}^{h}\mid _{S_{T_{2}}^{h}}=s_{0,i}^{h}\left( \mathbf{x},t\right) ,\text{ 
}\partial _{z}w_{i}^{h}\mid _{\Gamma _{0,T_{2}}^{h}}=s_{1,i}^{h}\left( 
\mathbf{x},t\right) ,\text{ }i=1,2,  \label{6.12}
\end{equation}%
\begin{equation}
\tau _{i}\mid _{\partial \Omega }=\tau _{0,i}\left( \mathbf{x}\right) ,\text{
}\partial _{z}\tau _{i}\mid _{\Gamma _{0}}=\tau _{1,i}\left( \mathbf{x}%
\right) ,\text{ }i=1,2.  \label{6.13}
\end{equation}%
\emph{\ Keeping in mind (\ref{2.6}), (\ref{4.026}) and (\ref{6.4}), let} 
\begin{equation}
\left. 
\begin{array}{c}
s_{0,1}^{h}\left( \mathbf{x},t\right) =s_{0,2}^{h}\left( \mathbf{x},t\right) 
\text{ \emph{for} }\left( \mathbf{x},t\right) \in S_{T_{2}}^{h}\diagdown
\Gamma _{0,T_{2}}^{h}, \\ 
\tau _{0,1}\left( \mathbf{x}\right) =\tau _{0,2}\left( \mathbf{x}\right) 
\text{ \emph{for} }\mathbf{x}\in \Gamma _{1}\cup \Gamma _{2}.%
\end{array}%
\right.   \label{6.14}
\end{equation}%
\emph{Furthermore, let by (\ref{5.32})} 
\begin{equation}
m_{i}\left( \mathbf{x}\right) =\left\vert \nabla \tau _{i}\left( \mathbf{x}%
\right) \right\vert ^{2},\text{ }i=1,2.  \label{6.15}
\end{equation}%
\emph{Then the following two Lipschitz stability estimates hold}%
\begin{equation}
\left. 
\begin{array}{c}
\left\Vert w_{1}^{h}-w_{2}^{h}\right\Vert _{H_{h}^{1}\left(
Q_{T_{2}}^{h}\right) }+\left\Vert \tau _{1}-\tau _{2}\right\Vert
_{H^{1}\left( \Omega \right) }\leq  \\ 
\leq C_{2}\left( \left\Vert s_{0,1}^{h}-s_{0,2}^{h}\right\Vert
_{L_{2,h}\left( \Gamma _{0,T_{2}}^{h}\right) }+\left\Vert
s_{1,1}^{h}-s_{1,2}^{h}\right\Vert _{L_{2,h}\left( \Gamma
_{0,T_{2}}^{h}\right) }\right) + \\ 
+C_{2}\left( \left\Vert \tau _{0,1}-\tau _{0,2}\right\Vert _{L_{2,h}\left(
\Gamma _{0}\right) }+\left\Vert \tau _{1,1}-\tau _{1,2}\right\Vert
_{L_{2}\left( \Gamma _{0}\right) }\right) ,%
\end{array}%
\right.   \label{6.16}
\end{equation}%
\begin{equation}
\left. 
\begin{array}{c}
\left\Vert m_{1}-m_{2}\right\Vert _{L_{2}\left( \Omega \right) }\leq
C_{2}\left( \left\Vert s_{0,1}^{h}-s_{0,2}^{h}\right\Vert _{L_{2,h}\left(
\Gamma _{0,T_{2}}^{h}\right) }+\left\Vert s_{1,1}^{h}-s_{1,2}^{h}\right\Vert
_{L_{2,h}\left( \Gamma _{0,T_{2}}^{h}\right) }\right) + \\ 
+C_{2}\left( \left\Vert \tau _{0,1}-\tau _{0,2}\right\Vert _{L_{2,h}\left(
\Gamma _{0}\right) }+\left\Vert \tau _{1,1}-\tau _{1,2}\right\Vert
_{L_{2}\left( \Gamma _{0}\right) }\right) .%
\end{array}%
\right.   \label{6.160}
\end{equation}%
\emph{In (\ref{6.16}) and (\ref{6.160}) }$H_{h}^{1}\left(
Q_{T_{2}}^{h}\right) $\emph{\ and }$L_{2,h}\left( \Gamma
_{0,T_{2}}^{h}\right) $\emph{\ are spaces defined in (\ref{3.7}), where }$T$%
\emph{\ is replaced with }$T_{2},$\emph{\ and also the number }$%
C_{2}=C_{2}\left( h_{0},M_{2},M_{3},T_{2},\Omega \right) >0$\emph{\ depends
only on listed parameters.} \emph{In particular, setting in (\ref{6.16}) and
(\ref{6.160}) }$s_{0,1}^{h}=s_{0,2}^{h}$\emph{\ in }$\Gamma _{0,T_{2}}^{h}$%
\emph{\ and }$\tau _{1,1}=\tau _{1,2}$\emph{\ in }$\Gamma _{0},$\emph{\ we
obtain }$w_{1}^{h}=w_{2}^{h}$\emph{\ in }$Q_{T_{2}}^{h},$\emph{\ }$\tau
_{1}=\tau _{2}$\emph{\ in }$\Omega $\emph{\ and }$m_{1}=m_{2}$\emph{\ in }$%
\Omega ,$\emph{\ i.e. we obtain uniqueness. }

Below $C_{2}>0$ denotes different numbers dependent on the same parameters
as ones listed above.

\section{Proofs of Theorems 6.2 and 6.3}

\label{sec:7}

As to the proof of Theorem 6.1, we refer to Remark 6.2.

\subsection{Proof of Theorem 6.2}

\label{sec:7.1}

Denote%
\begin{equation}
\left. 
\begin{array}{c}
\emph{\ }\widetilde{r}^{h}\left( \mathbf{x},t\right) =r_{1}^{h}\left( 
\mathbf{x},t\right) -r_{2}^{h}\left( \mathbf{x},t\right) ,\text{ } \\ 
\widetilde{g}_{0}^{h}\left( \mathbf{x},t\right) =\overline{g}%
_{0,1}^{h}\left( \mathbf{x},t\right) -\overline{g}_{0,2}^{h}\left( \mathbf{x}%
,t\right) ,\text{ }\widetilde{g}_{1}^{h}\left( \mathbf{x},t\right) =%
\overline{g}_{1,1}^{h}\left( \mathbf{x},t\right) -\overline{g}%
_{1,2}^{h}\left( \mathbf{x},t\right) .%
\end{array}%
\right.   \label{7.9}
\end{equation}%
Substitute in (\ref{4.25}) $r_{1}^{h}\left( \mathbf{x},t\right) $ first and
then $r_{2}^{h}\left( \mathbf{x},t\right) .$ Next, subtract equation for $%
r_{2}^{h}\left( \mathbf{x},t\right) $ from equation for $r_{1}^{h}\left( 
\mathbf{x},t\right) $ and use notation in the first line of (\ref{7.9}). We
obtain%
\begin{equation*}
\Delta \widetilde{r}\left( \mathbf{x},t_{j}\right) +2\nabla \widetilde{r}%
\left( \mathbf{x},t_{j}\right) \left[ \left( \displaystyle\int \limits%
_{0}^{t_{j}}\nabla r_{1}^{h}\left( \mathbf{x},\tau \right) d\tau \right)
_{h}-\nabla \ln \left( \left\vert \mathbf{x}\right\vert \right) \right] +
\end{equation*}%
\begin{equation}
+2\nabla r_{2}\left( \mathbf{x},t_{j}\right) \left( \displaystyle\int \limits%
_{0}^{t_{j}}\nabla \widetilde{r}^{h}\left( \mathbf{x},\tau \right) d\tau
\right) _{h}-2\displaystyle\sum \limits_{i=1}^{3}\frac{x_{i}}{\left\vert 
\mathbf{x}\right\vert }\partial _{t}^{h}\widetilde{r}_{x_{i}}\left( \mathbf{x%
},t_{j}\right) -  \label{7.10}
\end{equation}%
\begin{equation*}
-2\widetilde{r}\left( \mathbf{x},t_{j}\right) \displaystyle\sum \limits%
_{i=1}^{3}\frac{x_{i}}{\left\vert \mathbf{x}\right\vert }\partial
_{x_{i}}r_{1}\left( \mathbf{x},t_{j}\right) -2r_{2}\left( \mathbf{x}%
,t_{j}\right) \displaystyle\sum \limits_{i=1}^{3}\frac{x_{i}}{\left\vert 
\mathbf{x}\right\vert }\partial _{x_{i}}\widetilde{r}\left( \mathbf{x}%
,t_{j}\right) -
\end{equation*}%
\begin{equation*}
-2\partial _{t}^{h}\widetilde{r}\left( \mathbf{x},t_{j}\right) \displaystyle%
\sum \limits_{i=1}^{3}\frac{x_{i}}{\left\vert \mathbf{x}\right\vert }\left( %
\displaystyle\int \limits_{0}^{t_{j}}\partial _{x_{i}}r_{1}^{h}\left( 
\mathbf{x},\tau \right) d\tau \right) _{h}-
\end{equation*}%
\begin{equation*}
-2\partial _{t}^{h}r_{2}\left( \mathbf{x},t_{j}\right) \displaystyle\sum %
\limits_{i=1}^{3}\frac{x_{i}}{\left\vert \mathbf{x}\right\vert }\left( %
\displaystyle\int \limits_{0}^{t_{j}}\partial _{x_{i}}\widetilde{r}%
^{h}\left( \mathbf{x},\tau \right) d\tau \right) _{h}-
\end{equation*}%
\begin{equation*}
-2\partial _{t}^{h}\widetilde{r}^{h}\left( \mathbf{x},t_{j}\right) \left( 
\frac{1}{\left\vert \mathbf{x}\right\vert }-1\right) =0,\text{ }\mathbf{x}%
\in \Omega ,\text{ }j=0,...,k.
\end{equation*}%
By (\ref{3.2}), (\ref{3.3}) and (\ref{6.6})%
\begin{equation}
\left. 
\begin{array}{c}
\left\vert \partial _{t}^{h}f\left( t_{j}\right) \right\vert \leq \left(
\left\vert f\left( t_{j+1}\right) \right\vert +\left\vert f\left(
t_{j}\right) \right\vert \right) /h_{0},\text{ }j=0,...,k-1, \\ 
\left\vert \partial _{t}^{h}f\left( t_{k}\right) \right\vert \text{ }\leq
2\left( \left\vert f\left( t_{k}\right) \right\vert +\left\vert f\left(
t_{k-1}\right) \right\vert +\left\vert f\left( t_{k-2}\right) \right\vert
\right) /h_{0}, \\ 
\forall f\in C^{2}\left[ 0,T_{1}\right] .%
\end{array}%
\right.   \label{7.11}
\end{equation}%
Hence, using (\ref{3.5}), (\ref{6.5}), (\ref{7.10}) and (\ref{7.11}), we
obtain the following inequality%
\begin{equation}
\left. 
\begin{array}{c}
\left\vert \Delta \widetilde{r}\left( \mathbf{x},t_{j}\right) \right\vert
\leq C_{1}\displaystyle\sum \limits_{i=0}^{k}\left( \left\vert \nabla 
\widetilde{r}\left( \mathbf{x},t_{i}\right) \right\vert +\left\vert 
\widetilde{r}\left( \mathbf{x},t_{i}\right) \right\vert \right) ,\text{ } \\ 
\mathbf{x}\in \Omega ,\text{ }j=0,...,k.%
\end{array}%
\right.   \label{7.12}
\end{equation}%
In addition, by (\ref{6.7}), (\ref{6.8}) and (\ref{7.9})%
\begin{equation}
\left. 
\begin{array}{c}
\widetilde{r}^{h}\mid _{S_{T_{1}}^{h}\diagdown \Gamma _{0,T_{1}}^{h}}=0\text{%
,} \\ 
\text{ }\widetilde{r}^{h}\mid _{\Gamma _{0,T_{1}}^{h}}=\widetilde{g}%
_{0}^{h}\left( \mathbf{x},t\right) ,\text{ }\partial _{z}\widetilde{r}%
^{h}\mid _{\Gamma _{0,T_{1}}^{h}}=\widetilde{g}_{1}^{h}\left( \mathbf{x}%
,t\right) .%
\end{array}%
\right. \text{ }  \label{7.13}
\end{equation}%
Square both sides of (\ref{7.12}). Then multiply by the function $\varphi
_{\lambda }^{2}\left( z\right) $ defined in (\ref{6.2}), use Cauchy-Schwarz
inequality and integrate over the domain $\Omega $. We obtain for  $j=0,...,k
$%
\begin{equation}
\displaystyle\int \limits_{\Omega }\left( \Delta \widetilde{r}\left( \mathbf{%
x},t_{j}\right) \right) ^{2}\varphi _{\lambda }^{2}\left( z\right) d\mathbf{%
x\leq }C_{1}\displaystyle\sum \limits_{i=0}^{k}\displaystyle\int \limits%
_{\Omega }\left( \left\vert \nabla \widetilde{r}\left( \mathbf{x}%
,t_{i}\right) \right\vert ^{2}+\widetilde{r}^{2}\left( \mathbf{x}%
,t_{i}\right) \right) \varphi _{\lambda }^{2}\left( z\right) d\mathbf{x}.
\label{7.14}
\end{equation}%
Applying Carleman estimate (\ref{6.40}) of Theorem 6.1 to the left hand side
of (\ref{7.14}), using (\ref{3.5}), (\ref{3.7}), (\ref{6.2}), (\ref{6.3})
and (\ref{7.13}) and summing up resulting inequalities with respect to $%
j=0,...,k$, we obtain%
\begin{equation*}
\left. 
\begin{array}{c}
C_{1}\lambda ^{3}e^{2\lambda \left( a-d\right) ^{2}}\left( \left\Vert 
\widetilde{g}_{0}^{h}\right\Vert _{L_{2,h}\left( \Gamma
_{0,T_{1}}^{h}\right) }^{2}+\left\Vert \widetilde{g}_{1}^{h}\right\Vert
_{L_{2,h}\left( \Gamma _{0,T_{1}}^{h}\right) }^{2}\right) + \\ 
+C_{1}\displaystyle\sum \limits_{i=0}^{k}\displaystyle\int \limits_{\Omega
}\left( \left\vert \nabla \widetilde{r}\left( \mathbf{x},t_{i}\right)
\right\vert ^{2}+\widetilde{r}^{2}\left( \mathbf{x},t_{i}\right) \right)
\varphi _{\lambda }^{2}\left( z\right) d\mathbf{x}\geq  \\ 
\geq \lambda \displaystyle\sum \limits_{i=0}^{k}\displaystyle\int \limits%
_{\Omega }\left\vert \nabla \widetilde{r}\left( \mathbf{x},t_{i}\right)
\right\vert ^{2}\varphi _{\lambda }^{2}\left( z\right) d\mathbf{x+}\lambda
^{3}\displaystyle\sum \limits_{i=0}^{k}\displaystyle\int \limits_{\Omega }%
\widetilde{r}^{2}\left( \mathbf{x},t_{i}\right) \varphi _{\lambda
}^{2}\left( z\right) d\mathbf{x,} \\ 
\forall \lambda \geq \lambda _{0}.%
\end{array}%
\right. 
\end{equation*}%
Choosing here $\lambda _{1}\geq \lambda _{0}$ so large that $\lambda
_{1}>2C_{1},$ we obtain%
\begin{equation}
\left. 
\begin{array}{c}
\displaystyle\sum \limits_{i=0}^{k}\displaystyle\int \limits_{\Omega }\left(
\left\vert \nabla \widetilde{r}\left( \mathbf{x},t_{i}\right) \right\vert
^{2}+\widetilde{r}^{2}\left( \mathbf{x},t_{i}\right) \right) \varphi
_{\lambda _{1}}^{2}\left( z\right) d\mathbf{x\leq } \\ 
\leq C_{1}\lambda _{1}^{2}e^{2\lambda _{1}\left( a-d\right) ^{2}}\left(
\left\Vert \widetilde{g}_{0}^{h}\right\Vert _{L_{2,h}\left( \Gamma
_{0,T_{1}}^{h}\right) }^{2}+\left\Vert \widetilde{g}_{1}^{h}\right\Vert
_{L_{2,h}\left( \Gamma _{0,T_{1}}^{h}\right) }^{2}\right) .%
\end{array}%
\right.   \label{7.15}
\end{equation}%
Since by (\ref{6.2}) and (\ref{6.3})   
\begin{equation}
\varphi _{\lambda }^{2}\left( z\right) \geq e^{2\lambda \left( B-d\right)
^{2}}\text{ in }\Omega ,\text{ }\forall \lambda >0,  \label{7.150}
\end{equation}%
then 
\begin{equation*}
\left. 
\begin{array}{c}
\displaystyle\sum \limits_{i=0}^{k}\displaystyle\int \limits_{\Omega }\left(
\left\vert \nabla \widetilde{r}\left( \mathbf{x},t_{i}\right) \right\vert
^{2}+\widetilde{r}^{2}\left( \mathbf{x},t_{i}\right) \right) \varphi
_{\lambda _{1}}^{2}\left( z\right) d\mathbf{x}\geq  \\ 
\geq e^{2\lambda _{1}\left( B-d\right) ^{2}}\displaystyle\sum \limits%
_{i=0}^{k}\displaystyle\int \limits_{\Omega }\left( \left\vert \nabla 
\widetilde{r}\left( \mathbf{x},t_{i}\right) \right\vert ^{2}+\widetilde{r}%
^{2}\left( \mathbf{x},t_{i}\right) \right) \varphi _{\lambda _{1}}^{2}\left(
z\right) d\mathbf{x.}%
\end{array}%
\right. 
\end{equation*}%
Hence, dividing both sides of (\ref{7.15}) by $e^{2\lambda _{1}\left(
B-d\right) ^{2}},$ we obtain 
\begin{equation}
\left. 
\begin{array}{c}
\displaystyle\sum \limits_{i=0}^{k}\displaystyle\int \limits_{\Omega }\left(
\left\vert \nabla \widetilde{r}\left( \mathbf{x},t_{i}\right) \right\vert
^{2}+\widetilde{r}^{2}\left( \mathbf{x},t_{i}\right) \right) d\mathbf{x\leq }
\\ 
\leq C_{1}e^{2\lambda _{1}\left( a-d\right) ^{2}}\left( \left\Vert 
\widetilde{g}_{0}^{h}\right\Vert _{L_{2,h}\left( \Gamma
_{0,T_{1}}^{h}\right) }^{2}+\left\Vert \widetilde{g}_{1}^{h}\right\Vert
_{L_{2,h}\left( \Gamma _{0,T_{1}}^{h}\right) }^{2}\right) .%
\end{array}%
\right.   \label{7.16}
\end{equation}%
It follows from (\ref{3.7}) and (\ref{7.9}) that (\ref{7.16}) implies the
first target estimate (\ref{6.9}) of this theorem. Next, by (\ref{6.9}) and (%
\ref{7.12})%
\begin{equation}
\left. 
\begin{array}{c}
\displaystyle\sum \limits_{j=0}^{k}\displaystyle\int \limits_{\Omega }\left(
\Delta \widetilde{r}\left( \mathbf{x},t_{j}\right) \right) ^{2}\leq  \\ 
\leq C_{1}\left( \left\Vert \overline{g}_{0,1}^{h}-\overline{g}%
_{0,2}^{h}\right\Vert _{L_{2,h}\left( \Gamma _{0,T_{1}}^{h}\right)
}+\left\Vert \overline{g}_{1,1}^{h}-\overline{g}_{1,2}^{h}\right\Vert
_{L_{2,h}\left( \Gamma _{0,T_{1}}^{h}\right) }\right) .%
\end{array}%
\right.   \label{7.17}
\end{equation}%
Thus, (\ref{4.28}), (\ref{6.5}), (\ref{6.9}), (\ref{7.11}) and (\ref{7.17})
imply the second target estimate (\ref{6.10}) of this theorem. $\square $

\subsection{Proof of Theorem 6.3}

\label{sec:7.2}

Denote%
\begin{equation}
\widetilde{w}^{h}\left( \mathbf{x},t\right) =w_{1}^{h}\left( \mathbf{x}%
,t\right) -w_{2}^{h}\left( \mathbf{x},t\right) ,\text{ }\widetilde{\tau }%
\left( \mathbf{x}\right) =\tau _{1}\left( \mathbf{x}\right) -\tau _{2}\left( 
\mathbf{x}\right) ,  \label{7.18}
\end{equation}%
\begin{equation}
\left. 
\begin{array}{c}
\widetilde{s}_{0}^{h}\left( \mathbf{x},t\right) =s_{0,1}^{h}\left( \mathbf{x}%
,t\right) -s_{0,2}^{h}\left( \mathbf{x},t\right) ,\text{ }\widetilde{s}%
_{1}^{h}\left( \mathbf{x},t\right) =s_{1,1}^{h}\left( \mathbf{x},t\right)
-s_{1,2}^{h}\left( \mathbf{x},t\right) ,\emph{\ } \\ 
\widetilde{\tau }_{0}\left( \mathbf{x}\right) =\tau _{0,1}\left( \mathbf{x}%
\right) -\tau _{0,2}\left( \mathbf{x}\right) ,\text{ }\widetilde{\tau }%
_{1}\left( \mathbf{x}\right) =\tau _{1,1}\left( \mathbf{x}\right) -\tau
_{1,2}\left( \mathbf{x}\right) .%
\end{array}%
\right.   \label{7.19}
\end{equation}%
By (\ref{6.14}) and (\ref{7.19})%
\begin{equation}
\left. 
\begin{array}{c}
\widetilde{s}_{0}^{h}\left( \mathbf{x},t\right) =0\text{ for }\left( \mathbf{%
x},t\right) \in S_{T_{2}}^{h}\diagdown \Gamma _{0,T_{2}}^{h}, \\ 
\widetilde{\tau }_{0}\left( \mathbf{x}\right) =0\text{ for }\mathbf{x}\in
\Gamma _{1}\cup \Gamma _{2}.%
\end{array}%
\right.   \label{7.20}
\end{equation}

First, substitute in equations (\ref{5.26}), (\ref{5.27}) the pair $\left(
w_{1}^{h}\left( \mathbf{x},t\right) ,\tau _{1}\left( \mathbf{x}\right)
\right) $ and next the pair $\left( w_{2}^{h}\left( \mathbf{x},t\right)
,\tau _{2}\left( \mathbf{x}\right) \right) .$ Then subtract equations for
the second pair from equations for the first pair. Using (\ref{7.18}), we
obtain%
\begin{equation*}
\Delta \widetilde{w}\left( \mathbf{x},t_{j}\right) -2\nabla \partial _{t}^{h}%
\widetilde{w}\left( \mathbf{x},t_{j}\right) \nabla \tau _{1}\left( \mathbf{x}%
\right) -2\nabla _{\mathbf{x}}\partial _{t}^{h}w_{2}\left( \mathbf{x}%
,t_{j}\right) \cdot \nabla \widetilde{\tau }\left( \mathbf{x}\right) +
\end{equation*}%
\begin{equation*}
+2\partial _{t}^{h}\widetilde{w}\left( \mathbf{x},t_{j}\right) \frac{\nabla
w_{1}\left( \mathbf{x},t_{0}\right) }{w_{1}\left( \mathbf{x},t_{0}\right) }%
\nabla \tau _{1}+
\end{equation*}%
\begin{equation}
+2\partial _{t}^{h}w_{2}\left( \mathbf{x},t_{j}\right) \frac{\nabla 
\widetilde{w}\left( \mathbf{x},t_{0}\right) w_{2}\left( \mathbf{x}%
,t_{0}\right) +\nabla w_{2}\left( \mathbf{x},t_{0}\right) \widetilde{w}%
\left( \mathbf{x},t_{0}\right) }{w_{1}\left( \mathbf{x},t_{0}\right)
w_{2}\left( \mathbf{x},t_{0}\right) }\nabla \tau _{1}\left( \mathbf{x}%
\right) +  \label{7.21}
\end{equation}%
\begin{equation*}
+2\partial _{t}^{h}w_{2}\left( \mathbf{x},t_{j}\right) \frac{\nabla
w_{2}\left( \mathbf{x},t_{0}\right) }{w_{2}\left( \mathbf{x},t_{0}\right) }%
\nabla \widetilde{\tau }\left( \mathbf{x}\right) =0,\text{ }\mathbf{x}\in
\Omega ,\text{ }j=0,...k,
\end{equation*}%
\begin{equation}
\Delta \widetilde{\tau }\left( \mathbf{x}\right) +2\frac{\nabla \widetilde{w}%
\left( \mathbf{x},t_{0}\right) w_{2}\left( \mathbf{x},t_{0}\right) +\nabla
w_{2}\left( \mathbf{x},t_{0}\right) \widetilde{w}\left( \mathbf{x}%
,t_{0}\right) }{w_{1}\left( \mathbf{x},t_{0}\right) w_{2}\left( \mathbf{x}%
,t_{0}\right) }\nabla \tau _{1}\left( \mathbf{x}\right) +  \label{7.22}
\end{equation}%
\begin{equation*}
+2\frac{\nabla w_{2}\left( \mathbf{x},t_{0}\right) }{w_{2}\left( \mathbf{x}%
,t_{0}\right) }\nabla \widetilde{\tau }\left( \mathbf{x}\right) =0,\text{ }%
\mathbf{x}\in \Omega .
\end{equation*}%
In addition, by (\ref{6.12})-(\ref{6.14}) and (\ref{7.18})-(\ref{7.20}) we
have the following boundary conditions for the vector function $\widetilde{w}%
^{h}\left( \mathbf{x},t\right) $ and the function $\widetilde{\tau }\left( 
\mathbf{x}\right) $%
\begin{equation}
\left. 
\begin{array}{c}
\widetilde{w}^{h}\left( \mathbf{x},t\right) \mid _{S_{T_{2}}^{h}\diagdown
\Gamma _{0,T_{2}}^{h}}=0,\text{ } \\ 
\text{ }\widetilde{w}^{h}\left( \mathbf{x},t\right) \mid _{\Gamma
_{0,T_{2}}^{h}}=\widetilde{s}_{0}^{h}\left( \mathbf{x},t\right) ,\text{ }%
\partial _{z}\widetilde{w}^{h}\left( \mathbf{x},t\right) \mid _{\Gamma
_{0,T_{2}}^{h}}=\widetilde{s}_{1}^{h}\left( \mathbf{x},t\right) ,%
\end{array}%
\right.   \label{7.23}
\end{equation}%
\begin{equation}
\left. 
\begin{array}{c}
\widetilde{\tau }\left( \mathbf{x}\right) \mid _{\Gamma _{1}\cup \Gamma
_{2}}=0, \\ 
\widetilde{\tau }\left( \mathbf{x}\right) \mid _{\Gamma _{0}}=\widetilde{%
\tau }_{0}\left( \mathbf{x}\right) ,\text{ }\partial _{z}\widetilde{\tau }%
\left( \mathbf{x}\right) \mid _{\Gamma _{0}}=\widetilde{\tau }_{1}\left( 
\mathbf{x}\right) .%
\end{array}%
\right.   \label{7.24}
\end{equation}

Using (\ref{6.11}), (\ref{6.110}), (\ref{7.11}), (\ref{7.21}) and (\ref{7.22}%
), we obtain%
\begin{equation}
\left. 
\begin{array}{c}
\left\vert \Delta \widetilde{w}\left( \mathbf{x},t_{j}\right) \right\vert
\leq C_{2}\displaystyle\sum \limits_{i=0}^{k}\left( \left\vert \nabla 
\widetilde{w}\left( \mathbf{x},t_{i}\right) \right\vert +\left\vert 
\widetilde{w}\left( \mathbf{x},t_{i}\right) \right\vert \right)
+C_{2}\left\vert \nabla \widetilde{\tau }\left( \mathbf{x}\right)
\right\vert , \\ 
\left\vert \Delta \widetilde{\tau }\left( \mathbf{x}\right) \right\vert \leq
C_{2}\left( \left\vert \nabla \widetilde{w}\left( \mathbf{x},t_{0}\right)
\right\vert +\left\vert \widetilde{w}\left( \mathbf{x},t_{0}\right)
\right\vert \right) +C_{2}\left\vert \nabla \widetilde{\tau }\left( \mathbf{x%
}\right) \right\vert , \\ 
\mathbf{x}\in \Omega ,\text{ }j=0,...,k.%
\end{array}%
\right.   \label{7.25}
\end{equation}%
Square each of two inequalities (\ref{7.25}), then sum them up, multiply by
the function $\varphi _{\lambda }^{2}\left( z\right) ,$ integrate over the
domain $\Omega $ and use Cauchy-Schwarz inequality. We obtain%
\begin{equation*}
\left. 
\begin{array}{c}
\displaystyle\int \limits_{\Omega }\left( \Delta \widetilde{w}\left( \mathbf{%
x},t_{j}\right) \right) ^{2}\varphi _{\lambda }^{2}\left( z\right) d\mathbf{x%
}+\displaystyle\int \limits_{\Omega }\left( \Delta \widetilde{\tau }\left( 
\mathbf{x}\right) \right) ^{2}\varphi _{\lambda }^{2}\left( z\right) d%
\mathbf{x}\leq  \\ 
\leq C_{2}\displaystyle\sum \limits_{i=0}^{k}\displaystyle\int \limits%
_{\Omega }\left( \left\vert \nabla \widetilde{w}\left( \mathbf{x}%
,t_{i}\right) \right\vert ^{2}+\widetilde{w}^{2}\left( \mathbf{x}%
,t_{i}\right) \right) \varphi _{\lambda }^{2}\left( z\right) d\mathbf{x}%
+C_{2}\displaystyle\int \limits_{\Omega }\left\vert \nabla \widetilde{\tau }%
\left( \mathbf{x}\right) \right\vert ^{2}\varphi _{\lambda }^{2}\left(
z\right) d\mathbf{x.}%
\end{array}%
\right. 
\end{equation*}%
Applying to the left hand side of this inequality Carleman estimate (\ref%
{6.40}) of Theorem 6.1, using boundary conditions (\ref{7.23}) and (\ref%
{7.24}), summing up resulting inequalities and using (\ref{3.5}), (\ref{3.7}%
), (\ref{6.2}) and (\ref{6.3}), we obtain%
\begin{equation}
\left. 
\begin{array}{c}
C_{2}\lambda ^{3}e^{2\lambda \left( a-d\right) ^{2}}\left( \left\Vert 
\widetilde{s}_{0}^{h}\right\Vert _{L_{2,h}\left( \Gamma
_{0,T_{2}}^{h}\right) }^{2}+\left\Vert \widetilde{s}_{1}^{h}\right\Vert
_{L_{2,h}\left( \Gamma _{0,T_{2}}^{h}\right) }^{2}\right) + \\ 
+C_{2}\lambda ^{3}e^{2\lambda \left( a-d\right) ^{2}}\left( \left\Vert 
\widetilde{\tau }_{0}\right\Vert _{L_{2}\left( \Gamma _{0}\right)
}^{2}+\left\Vert \widetilde{\tau }_{1}\right\Vert _{L_{2}\left( \Gamma
_{0}\right) }^{2}\right) + \\ 
+C_{2}\displaystyle\sum \limits_{i=0}^{k}\displaystyle\int \limits_{\Omega
}\left( \left\vert \nabla \widetilde{w}\left( \mathbf{x},t_{i}\right)
\right\vert ^{2}+\widetilde{w}^{2}\left( \mathbf{x},t_{i}\right) \right)
\varphi _{\lambda }^{2}\left( z\right) d\mathbf{x}+ \\ 
+C_{2}\displaystyle\int \limits_{\Omega }\left\vert \nabla \widetilde{\tau }%
\left( \mathbf{x}\right) \right\vert ^{2}\varphi _{\lambda }^{2}\left(
z\right) d\mathbf{x}\geq  \\ 
\geq \lambda \displaystyle\sum \limits_{i=0}^{k}\displaystyle\int \limits%
_{\Omega }\left\vert \nabla \widetilde{w}\left( \mathbf{x},t_{i}\right)
\right\vert ^{2}\varphi _{\lambda }^{2}\left( z\right) d\mathbf{x+}\lambda
^{3}\displaystyle\sum \limits_{i=0}^{k}\displaystyle\int \limits_{\Omega }%
\widetilde{w}^{2}\left( \mathbf{x},t_{i}\right) \varphi _{\lambda
}^{2}\left( z\right) d\mathbf{x+} \\ 
+\lambda \displaystyle\int \limits_{\Omega }\left\vert \nabla \widetilde{%
\tau }\left( \mathbf{x}\right) \right\vert ^{2}\varphi _{\lambda }^{2}\left(
z\right) d\mathbf{x}+\lambda ^{3}\displaystyle\int \limits_{\Omega }%
\widetilde{\tau }^{2}\left( \mathbf{x}\right) \varphi _{\lambda }^{2}\left(
z\right) d\mathbf{x,} \\ 
\forall \lambda \geq \lambda _{0}.%
\end{array}%
\right.   \label{7.26}
\end{equation}%
Select $\lambda _{2}\geq \lambda _{0}$ so large that $\lambda _{2}>2C_{2}.$
Then (\ref{7.26}) implies%
\begin{equation}
\left. 
\begin{array}{c}
\displaystyle\sum \limits_{i=0}^{k}\displaystyle\int \limits_{\Omega }\left(
\left\vert \nabla \widetilde{w}\left( \mathbf{x},t_{i}\right) \right\vert
^{2}+\widetilde{w}^{2}\left( \mathbf{x},t_{i}\right) \right) \varphi
_{\lambda _{2}}^{2}\left( z\right) d\mathbf{x}+ \\ 
+\displaystyle\int \limits_{\Omega }\left( \left\vert \nabla \widetilde{\tau 
}\left( \mathbf{x}\right) \right\vert ^{2}+\widetilde{\tau }^{2}\left( 
\mathbf{x}\right) \right) \varphi _{\lambda _{2}}^{2}\left( z\right) d%
\mathbf{x\leq } \\ 
\leq C_{2}\lambda _{2}^{2}e^{2\lambda _{2}\left( a-d\right) ^{2}}\left(
\left\Vert \widetilde{s}_{0}^{h}\right\Vert _{L_{2,h}\left( \Gamma
_{0,T_{2}}^{h}\right) }^{2}+\left\Vert \widetilde{s}_{1}^{h}\right\Vert
_{L_{2,h}\left( \Gamma _{0,T_{2}}^{h}\right) }^{2}\right) + \\ 
+C_{2}\lambda _{2}^{2}e^{2\lambda _{2}\left( a-d\right) ^{2}}\left(
\left\Vert \widetilde{\tau }_{0}\right\Vert _{L_{2}\left( \Gamma _{0}\right)
}^{2}+\left\Vert \widetilde{\tau }_{1}\right\Vert _{L_{2}\left( \Gamma
_{0}\right) }^{2}\right) .%
\end{array}%
\right.   \label{7.27}
\end{equation}%
Use again (\ref{7.150}), although replacing in it $\lambda _{1}$ with $%
\lambda _{2}.$ Hence, dividing (\ref{7.27}) by $e^{2\lambda _{2}\left(
B-d\right) ^{2}}$ and using (\ref{3.5}) and (\ref{3.7}), we obtain (\ref%
{6.16}), which is the first target estimate of this theorem. 

Next, by (\ref{6.11}) and (\ref{6.15}) 
\begin{equation}
\left\vert m_{1}\left( \mathbf{x}\right) -m_{2}\left( \mathbf{x}\right)
\right\vert \leq C_{2}\left\vert \nabla \tau _{1}\left( \mathbf{x}\right)
-\nabla \tau _{2}\left( \mathbf{x}\right) \right\vert =C_{2}\left\vert
\nabla \widetilde{\tau }\left( \mathbf{x}\right) \right\vert .  \label{7.29}
\end{equation}%
Thus, it follows from (\ref{3.5}), (\ref{3.7}), (\ref{6.16}), (\ref{7.18})
and (\ref{7.19}) that (\ref{7.27}) and (\ref{7.29}) imply the second target
estimate (\ref{6.160}) of this theorem. $\square $

\section{Convexification Numerical Method for CIP1}

\label{sec:8}

\subsection{Convexification functional}

\label{sec:8.1}

Let $N\geq 2$ be an integer and $B$ be a Banach space with the norm $%
\left\Vert \cdot \right\Vert _{B}.$ Denote 
\begin{equation}
\left. 
\begin{array}{c}
B_{N}=B\times B\times ...\times B,\text{ }N\text{ times, } \\ 
\left\Vert u\right\Vert _{B_{N}}^{2}=\displaystyle\sum \limits%
_{i=1}^{N}\left\Vert u_{i}\right\Vert _{B}^{2},\text{ }\forall u=\left(
u_{1},...,u_{N}\right) ^{T}\in B_{N}.%
\end{array}%
\right.   \label{8.1}
\end{equation}%
It follows from (\ref{3.4}), (\ref{4.25}), (\ref{7.11}) and (\ref{8.1}) that
we need to ensure in the convexification method that the vector function 
\begin{equation}
r^{h}\left( \mathbf{x},t\right) \in C_{k}^{2}\left( \overline{\Omega }%
\right) ,\text{ }t\in \left\{ 0=t_{0},t_{1},...,t_{k}=T_{1}\right\} ,
\label{8.2}
\end{equation}%
By embedding theorem%
\begin{equation}
H_{k}^{4}\left( \Omega \right) \subset C_{k}^{2}\left( \overline{\Omega }%
\right) ,\text{ }\left\Vert f\right\Vert _{C_{k}^{1}\left( \overline{\Omega }%
\right) }\leq \overline{C}\left\Vert f\right\Vert _{H_{k}^{4}\left( \Omega
\right) },\text{ }\forall f\in H_{k}^{4}\left( \Omega \right) ,  \label{8.3}
\end{equation}%
where the number $\overline{C}=\overline{C}\left( \Omega ,k\right) >0$
depends only on listed parameters. Hence, we define the set $B\left(
R\right) $ as follows: Let $R>0$ be an arbitrary number. Then%
\begin{equation}
B\left( R\right) =\left\{ 
\begin{array}{c}
r^{h}\left( \mathbf{x},t\right) \in H_{h,k}^{4}\left( Q_{T_{1}}^{h}\right)
:\left\Vert r^{h}\right\Vert _{H_{h,k}^{4}\left( Q_{T_{1}}^{h}\right) }\leq
R, \\ 
r^{h}\left( \mathbf{x},t\right) \mid _{S_{T_{1}}^{h}}=\overline{g}%
_{0}^{h}\left( \mathbf{x},t\right) ,\text{ }r_{z}^{h}\left( \mathbf{x}%
,t\right) \mid _{\Gamma _{0,T_{1}}^{h}}=\overline{g}_{1}^{h}\left( \mathbf{x}%
,t\right) 
\end{array}%
\right\} .  \label{8.4}
\end{equation}%
Hence, $B\left( R\right) $ is a closed set in the space $H_{h,k}^{4}\left(
Q_{T_{1}}^{h}\right) ,$ also, see Remark 2.1. By (\ref{3.7}), (\ref{8.3})
and (\ref{8.4})%
\begin{equation}
B\left( R\right) \subset C^{2}\left( \overline{Q}_{T_{1}}^{h}\right) ,\text{ 
}\left\Vert r^{h}\right\Vert _{C^{2}\left( \overline{Q}_{T_{1}}^{h}\right)
}\leq \overline{C}R,\text{ }\forall r^{h}\in B\left( R\right) .  \label{8.40}
\end{equation}%
The boundary conditions in (\ref{8.4}) are the same as the ones as in (\ref%
{4.26}). Thus, (\ref{8.3}) and (\ref{8.4}) imply that (\ref{8.2}) holds for
all vector functions, which belong to the set $B\left( R\right) .$

Recall that the CWF $\varphi _{\lambda }\left( z\right) $ was defined in (%
\ref{6.2}). Let $L_{j}\left( r^{h}\right) \left( \mathbf{x},t_{j}\right) ,$ $%
j=0,...,k$ be the nonlinear operators defined in (\ref{4.25}). We define the
convexification functional $J_{\lambda ,\gamma }\left( r^{h}\right) $ as:%
\begin{equation}
\left. 
\begin{array}{c}
J_{\lambda ,\gamma }\left( r^{h}\right) :B\left( R\right) \rightarrow 
\mathbb{R}, \\ 
J_{\lambda ,\gamma }\left( r^{h}\right) =e^{-2\lambda \left( a-d\right) ^{2}}%
\displaystyle\sum \limits_{j=0}^{k}\displaystyle\int \limits_{\Omega }\left[
L_{j}\left( r^{h}\right) \left( \mathbf{x},t_{j}\right) \right] ^{2}\varphi
_{\lambda }^{2}\left( z\right) d\mathbf{x}+\gamma \left\Vert
r^{h}\right\Vert _{H_{h,k}^{4}\left( Q_{T_{1}}^{h}\right) }^{2},%
\end{array}%
\right.   \label{8.5}
\end{equation}%
where $\gamma \in \left( 0,1\right) $ is the regularization parameter, and $%
e^{-2\lambda \left( a-d\right) ^{2}}$ is the balancing factor between two
terms in the right hand side of (\ref{8.5}) since by (\ref{6.3}) $%
e^{-2\lambda \left( a-d\right) ^{2}}\max_{\left[ a,B\right] }\varphi
_{\lambda }^{2}\left( z\right) =1.$Below in this section we are concerned
with the following problem:

\textbf{Minimization Problem.} \emph{Minimize the functional }$J_{\lambda
,\gamma }\left( r^{h}\right) $\emph{\ on the set} $B\left( R\right) .$

\subsection{Global strong convexity of the functional $J_{\protect\lambda ,%
\protect\gamma }\left( r^{h}\right) $}

\label{sec:8.2}

In this subsection we prove strong convexity of the functional $J_{\lambda
,\gamma }\left( r^{h}\right) $ on the set $B\left( R\right) $ defined in (%
\ref{8.5}). Since a smallness condition is not imposed on the number $R$,
then it is natural to call this \textquotedblleft global strong convexity",
similarly with Definition 1.1. Define the subspace $H_{h,k}^{4,0}\left(
Q_{T_{1}}^{h}\right) $ of the space $H_{h,k}^{4}\left( Q_{T_{1}}^{h}\right) $
as%
\begin{equation}
H_{h,k}^{4,0}\left( Q_{T_{1}}^{h}\right) =\left\{ u^{h}\in H_{h,k}^{4}\left(
Q_{T_{1}}^{h}\right) :u^{h}\mid _{S_{T_{1}}^{h}}=\partial _{z}u^{h}\mid
_{\Gamma _{0,T_{1}}^{h}}=0\right\} .  \label{8.6}
\end{equation}

\textbf{Theorem 8.1.} \emph{Assume that condition (\ref{6.6}) holds. The
functional }$J_{\lambda ,\gamma }\left( r^{h}\right) $\emph{\ has the
following properties:}

\emph{1. For any }$\lambda ,\gamma >0$\emph{\ and for any }$r^{h}\in B\left(
R\right) $\emph{\ there exists Fr\'{e}chet derivative }

$J_{\lambda ,\gamma }^{\prime }\left( r^{h}\right) \in H_{h,k}^{4,0}\left(
Q_{T_{1}}^{h}\right) $\emph{\ of the functional }$J_{\lambda ,\gamma }\left(
r^{h}\right) $\emph{\ at the point }$r^{h}.$\emph{\ This derivative is
Lipschitz continuous on }$B\left( R\right) .$\emph{\ In other words, there
exists a number }$K=K\left( h_{0},R,\lambda ,\gamma ,T,c_{1},d,\Omega
\right) >0$\emph{\ depending only on listed parameters such that }%
\begin{equation}
\left\Vert J_{\lambda ,\gamma }^{\prime }\left( r_{2}^{h}\right) -J_{\lambda
,\gamma }^{\prime }\left( r_{1}^{h}\right) \right\Vert _{H_{h,k}^{4}\left(
Q_{T_{1}}^{h}\right) }\leq K\left\Vert r_{1}^{h}-r_{2}^{h}\right\Vert
_{H_{h,k}^{4}\left( Q_{T_{1}}^{h}\right) },\text{ }\forall
r_{1}^{h},r_{2}^{h}\in B\left( R\right) .  \label{8.7}
\end{equation}

\emph{2. Let }$\lambda _{0}=\lambda _{0}\left( d,\Omega \right) \geq 1$\emph{%
\ be the number chosen in Theorem 6.1. There exists a sufficiently large
number }$\lambda _{3}=\lambda _{3}\left( h_{0},R,T,c_{1},d,\Omega \right)
\geq \lambda _{0}$\emph{\ such that for all }$\lambda \geq \lambda _{3}$%
\emph{\ the functional }$J_{\lambda ,\gamma }\left( r^{h}\right) $\emph{\ is
strongly convex on the set }$B\left( R\right) ,$\emph{\ i.e. the following
inequality holds}%
\begin{equation}
\left. 
\begin{array}{c}
J_{\lambda ,\gamma }\left( r_{2}^{h}\right) -J_{\lambda ,\gamma }\left(
r_{1}^{h}\right) -\left( J_{\lambda ,\gamma }^{\prime }\left(
r_{1}^{h}\right) ,r_{2}^{h}-r_{1}^{h}\right) \geq  \\ 
\geq C_{3}\exp \left[ -2\lambda \left( \left( a-d\right) ^{2}-\left(
B-d\right) ^{2}\right) \right] \left\Vert r_{2}^{h}-r_{1}^{h}\right\Vert
_{H_{h,k}^{1}\left( Q_{T_{1}}^{h}\right) }^{2}+ \\ 
+\gamma \left\Vert r_{2}^{h}-r_{1}^{h}\right\Vert _{H_{h,k}^{4}\left(
Q_{T_{1}}^{h}\right) }^{2},\text{ }\forall r_{1}^{h},r_{2}^{h}\in B\left(
R\right) ,\text{ }\forall \lambda \geq \lambda _{3},%
\end{array}%
\right.   \label{8.8}
\end{equation}%
\emph{where }$\left( ,\right) $\emph{\ denotes the scalar product in }$%
H_{h,k}^{4}\left( Q_{T_{1}}^{h}\right) ,$ \emph{and the number }

$C_{3}=C_{3}\left( h_{0},R,T,c_{1},d,\Omega \right) >0$\emph{\ depends only
on listed parameters.}

\emph{3. For every }$\lambda \geq \lambda _{3}$\emph{\ there exists unique
minimizer }$r_{\min ,\lambda }^{h}\in B\left( R\right) $\emph{\ of the
functional }$J_{\lambda ,\gamma }\left( r^{h}\right) $\emph{\ on the set }$%
B\left( R\right) .$\emph{\ Furthermore, the following inequality holds}%
\begin{equation}
\left( J_{\lambda ,\gamma }^{\prime }\left( r_{\min ,\lambda }^{h}\right)
,r^{h}-r_{\min ,\lambda }^{h}\right) \geq 0,\text{ }\forall r^{h}\in B\left(
R\right) .  \label{8.9}
\end{equation}

\textbf{Proof.} Everywhere below $C_{3}=C_{3}\left( h_{0},R,T,c_{1},d,\Omega
\right) >0$ denotes different numbers depending only on listed parameters.
Let $r_{1}^{h}\left( \mathbf{x},t\right) \in B\left( R\right) $ and $%
r_{2}^{h}\left( \mathbf{x},t\right) \in B\left( R\right) $ be two arbitrary
points of the set $B\left( R\right) .$ Denote $\alpha ^{h}\left( \mathbf{x}%
,t\right) =r_{2}^{h}\left( \mathbf{x},t\right) -r_{1}^{h}\left( \mathbf{x}%
,t\right) .$ Then triangle inequality, (\ref{8.4}) and (\ref{8.6}) imply%
\begin{equation}
r_{2}^{h}\left( \mathbf{x},t\right) =\alpha ^{h}\left( \mathbf{x},t\right)
+r_{1}^{h}\left( \mathbf{x},t\right) ,  \label{8.10}
\end{equation}%
\begin{equation}
\alpha ^{h}\left( \mathbf{x},t\right) \in B_{0}\left( 2R\right) =\left\{
\beta ^{h}\left( \mathbf{x},t\right) \in H_{h,k}^{4,0}\left(
Q_{T_{1}}^{h}\right) :\left\Vert \beta ^{h}\right\Vert _{H_{h,k}^{4}\left(
Q_{T_{1}}^{h}\right) }\leq 2R\right\} .  \label{8.11}
\end{equation}

Next, using (\ref{4.25})\ and (\ref{8.10}), we obtain%
\begin{equation}
\left. 
\begin{array}{c}
L_{j}\left( r_{2}^{h}\right) \left( \mathbf{x},t\right) =L_{j}\left(
r_{1}^{h}+\alpha ^{h}\right) \left( \mathbf{x},t\right) = \\ 
=L_{j,lin}^{r_{1}^{h}}\left( \alpha ^{h}\right) \left( \mathbf{x},t\right)
+L_{j,nonlin}^{r_{1}^{h}}\left( \alpha ^{h}\right) \left( \mathbf{x}%
,t\right) +L_{j}\left( r_{1}^{h}\right) \left( \mathbf{x},t\right) .%
\end{array}%
\right.   \label{8.12}
\end{equation}%
Hence, 
\begin{equation}
\left. 
\begin{array}{c}
\left( L_{j}\left( r_{2}^{h}\right) \left( \mathbf{x},t\right) \right)
^{2}-\left( L_{j}\left( r_{1}^{h}\right) \left( \mathbf{x},t\right) \right)
^{2}= \\ 
=2L_{j,lin}^{r_{1}^{h}}\left( \alpha ^{h}\right) \left( \mathbf{x},t\right)
L_{j}\left( r_{1}^{h}\right) \left( \mathbf{x},t\right) + \\ 
+2L_{j,nonlin}^{r_{1}^{h}}\left( \alpha ^{h}\right) \left( \mathbf{x}%
,t\right) L_{j}\left( r_{1}^{h}\right) \left( \mathbf{x},t\right) + \\ 
+\left( L_{j,lin}^{r_{1}^{h}}\left( \alpha ^{h}\right) \left( \mathbf{x}%
,t\right) \right) ^{2}+\left( L_{j,nonlin}^{r_{1}^{h}}\left( \alpha
^{h}\right) \left( \mathbf{x},t\right) \right) ^{2}.%
\end{array}%
\right.   \label{8.13}
\end{equation}%
In (\ref{8.12})\ and (\ref{8.13}) operators $L_{j,lin}^{r_{1}^{h}}\left(
\alpha ^{h}\right) \left( \mathbf{x},t\right) $ and $%
L_{j,nonlin}^{r_{1}^{h}}\left( \alpha ^{h}\right) \left( \mathbf{x},t\right) 
$ depend on $\alpha ^{h}$ linearly and nonlinearly respectively. More,
precisely, $L_{j,lin}^{r_{1}^{h}}\left( \alpha ^{h}\right) \left( \mathbf{x}%
,t\right) $ is:%
\begin{equation}
\left. 
\begin{array}{c}
L_{j,lin}^{r_{1}^{h}}\left( \alpha ^{h}\right) \left( \mathbf{x},t\right)
=\Delta \alpha \left( \mathbf{x},t_{j}\right) + \\ 
+2\nabla \alpha \left( \mathbf{x},t_{j}\right) \left[ \left( \displaystyle%
\int \limits_{0}^{t_{j}}\nabla r_{1}^{h}\left( \mathbf{x},\tau \right) d\tau
\right) _{h}-\nabla \ln \left( \left\vert \mathbf{x}\right\vert \right) %
\right] + \\ 
+2r_{1}\left( \mathbf{x},t_{j}\right) \left( \displaystyle\int \limits%
_{0}^{t_{j}}\nabla \alpha ^{h}\left( \mathbf{x},\tau \right) d\tau \right)
_{h}- \\ 
-2\displaystyle\sum \limits_{i=1}^{3}\left( x_{i}/\left\vert \mathbf{x}%
\right\vert \right) \partial _{t}^{h}\alpha _{x_{i}}\left( \mathbf{x}%
,t_{j}\right) -2r_{1}\left( \mathbf{x},t_{j}\right) \displaystyle\sum \limits%
_{i=1}^{3}\left( x_{i}/\left\vert \mathbf{x}\right\vert \right) \partial
_{x_{i}}\alpha \left( \mathbf{x},t_{j}\right) - \\ 
-2\alpha \left( \mathbf{x},t_{j}\right) \displaystyle\sum \limits%
_{i=1}^{3}\left( x_{i}/\left\vert \mathbf{x}\right\vert \right) \partial
_{x_{i}}r_{1}\left( \mathbf{x},t_{j}\right) - \\ 
-2\partial _{t}^{h}\alpha \left( \mathbf{x},t_{j}\right) \displaystyle\sum %
\limits_{i=1}^{3}\left( x_{i}/\left\vert \mathbf{x}\right\vert \right)
\left( \displaystyle\int \limits_{0}^{t_{j}}\partial _{x_{i}}r_{1}\left( 
\mathbf{x},\tau \right) d\tau \right) _{h}- \\ 
-2\partial _{t}^{h}r_{1}\left( \mathbf{x},t_{j}\right) \displaystyle\sum %
\limits_{i=1}^{3}\left( x_{i}/\left\vert \mathbf{x}\right\vert \right)
\left( \displaystyle\int \limits_{0}^{t_{j}}\partial _{x_{i}}\alpha \left( 
\mathbf{x},\tau \right) d\tau \right) _{h}- \\ 
-2\partial _{t}^{h}\alpha \left( \mathbf{x},t_{j}\right) \left( \left(
1/\left\vert \mathbf{x}\right\vert \right) -1\right) ,\text{ }\mathbf{x}\in
\Omega ,\text{ }j=0,...,k.%
\end{array}%
\right.   \label{8.14}
\end{equation}

Next, the operator $L_{j,nonlin}^{r_{1}^{h}}\left( \alpha ^{h}\right) \left( 
\mathbf{x},t\right) $ is:%
\begin{equation}
\left. 
\begin{array}{c}
L_{j,nonlin}^{r_{1}^{h}}\left( \alpha ^{h}\right) \left( \mathbf{x},t\right)
=2\nabla \alpha \left( \mathbf{x},t_{j}\right) \left( \displaystyle\int %
\limits_{0}^{t_{j}}\nabla \alpha ^{h}\left( \mathbf{x},\tau \right) d\tau
\right) _{h}- \\ 
-2\alpha \left( \mathbf{x},t_{j}\right) \displaystyle\sum \limits%
_{i=1}^{3}\left( x_{i}/\left\vert \mathbf{x}\right\vert \right) \alpha
_{x_{i}}\left( \mathbf{x},t_{j}\right) - \\ 
-2\partial _{t}^{h}\alpha \left( \mathbf{x},t_{j}\right)
\sum\limits_{i=1}^{3}\left( x_{i}/\left\vert \mathbf{x}\right\vert \right)
\left( \int\limits_{0}^{t_{j}}\alpha _{x_{i}}^{h}\left( \mathbf{x},\tau
\right) d\tau \right) _{h},\text{ }\mathbf{x}\in \Omega ,\text{ }j=0,...,k.%
\end{array}%
\right.   \label{8.15}
\end{equation}%
In particular, it follows from (\ref{4.25}), (\ref{8.40}), (\ref{8.14}) and (%
\ref{8.15}) that the sum of terms in the third and fourth lines of (\ref%
{8.13}) can be estimated from the below as%
\begin{equation}
\left. 
\begin{array}{c}
\left( L_{j,lin}^{r_{1}^{h}}\left( \alpha ^{h}\right) \left( \mathbf{x}%
,t\right) \right) ^{2}+ \\ 
+2L_{j,nonlin}^{r_{1}^{h}}\left( \alpha ^{h}\right) \left( \mathbf{x}%
,t\right) L_{j}\left( r_{1}^{h}\right) \left( \mathbf{x},t\right) +\left(
L_{j,nonlin}^{r_{1}^{h}}\left( \alpha ^{h}\right) \left( \mathbf{x},t\right)
\right) ^{2}\geq  \\ 
\geq C_{3}\left( \Delta \alpha \left( \mathbf{x},t_{j}\right) \right)
^{2}-C_{3}\displaystyle\sum \limits_{i=1}^{k}\left( \left\vert \nabla \alpha
\left( \mathbf{x},t_{i}\right) \right\vert ^{2}+\alpha ^{2}\left( \mathbf{x}%
,t_{i}\right) \right) .%
\end{array}%
\right.   \label{8.16}
\end{equation}

Using (\ref{8.5}), (\ref{8.13}) and (\ref{8.14}), we obtain that the linear,
with respect to $\alpha ^{h}\left( \mathbf{x},t\right) ,$ part of the
expression $J_{\lambda ,\gamma }\left( r_{1}^{h}+\alpha ^{h}\right)
-J_{\lambda ,\gamma }\left( r_{1}^{h}\right) $ is%
\begin{equation}
\widehat{J}\left( \alpha ^{h}\right) =e^{-2\lambda \left( a-d\right)
^{2}}\sum\limits_{j=0}^{k}\int\limits_{\Omega }L_{j,lin}^{r_{2}^{h}}\left(
\alpha ^{h}\right) \varphi _{\lambda }^{2}\left( z\right) d\mathbf{x+}%
2\left( \alpha ^{h},r_{1}^{h}\right) .  \label{8.17}
\end{equation}%
It follows from (\ref{8.11}), (\ref{8.14}) and (\ref{8.17}) that we can
consider $\widehat{J}\left( \alpha ^{h}\right) $ as a bounded linear
functional $\widehat{J}\left( \alpha ^{h}\right) :H_{h,k}^{4,0}\left(
Q_{T_{1}}^{h}\right) \rightarrow \mathbb{R}.$ Hence, by Riesz theorem there
exists a unique element $\widetilde{J}\in H_{h,k}^{4,0}\left(
Q_{T_{1}}^{h}\right) $ such that 
\begin{equation}
\widehat{J}\left( \alpha ^{h}\right) =\left( \widetilde{J},\alpha
^{h}\right) ,\text{ }\forall \alpha ^{h}\in H_{h,k}^{4,0}\left(
Q_{T_{1}}^{h}\right) .  \label{8.18}
\end{equation}%
On the other hand, it is clear from (\ref{8.13}), (\ref{8.17}) and (\ref%
{8.18}) that%
\begin{equation*}
\lim_{\left\Vert \alpha ^{h}\right\Vert _{H_{h,k}^{4}\left(
Q_{T_{1}}^{h}\right) }\rightarrow 0}\left( \frac{J_{\lambda ,\gamma }\left(
r_{1}^{h}+\alpha ^{h}\right) -J_{\lambda ,\gamma }\left( r_{1}^{h}\right)
-\left( \widetilde{J},\alpha ^{h}\right) }{\left\Vert \alpha ^{h}\right\Vert
_{H_{h,k}^{4}\left( Q_{T_{1}}^{h}\right) }}\right) =0.
\end{equation*}%
Hence, $\widetilde{J}=J_{\lambda ,\gamma }^{\prime }\left( r_{1}^{h}\right)
\in H_{h,k}^{4,0}\left( Q_{T_{1}}^{h}\right) $ is the Fr\'{e}chet derivative
of the functional $J_{\lambda ,\gamma }\left( r^{h}\right) $ at the point $%
r_{1}^{h}.$ We omit the proof of the Lipschitz continuity property (\ref{8.7}%
) of the functional $J_{\lambda ,\gamma }^{\prime }$ since this proof is
completely similar with the proof of the corresponding part of the proof of
Theorem 3.1 of \cite{Bak} as well as of Theorem 5.3.1 of \cite{KL}.

Estimate now from the below the expression $J_{\lambda ,\gamma }\left(
r_{2}^{h}\right) -J_{\lambda ,\gamma }\left( r_{1}^{h}\right) -\left(
J_{\lambda ,\gamma }^{\prime }\left( r_{1}^{h}\right)
,r_{2}^{h}-r_{1}^{h}\right) .$ Using (\ref{8.5}), (\ref{8.10}), (\ref{8.13}%
), (\ref{8.16}), (\ref{8.17}) and (\ref{8.18}) we obtain%
\begin{equation}
\left. 
\begin{array}{c}
J_{\lambda ,\gamma }\left( r_{2}^{h}\right) -J_{\lambda ,\gamma }\left(
r_{1}^{h}\right) -\left( J_{\lambda ,\gamma }^{\prime }\left(
r_{1}^{h}\right) ,r_{2}^{h}-r_{1}^{h}\right) \geq  \\ 
\geq e^{-2\lambda \left( a-d\right) ^{2}}C_{3}\displaystyle\sum \limits%
_{j=0}^{k}\displaystyle\int \limits_{\Omega }\left( \Delta \alpha \left( 
\mathbf{x},t_{j}\right) \right) ^{2}\varphi _{\lambda }^{2}\left( z\right) d%
\mathbf{x}- \\ 
-C_{3}e^{-2\lambda \left( a-d\right) ^{2}}\displaystyle\sum \limits_{i=1}^{k}%
\displaystyle\int \limits_{\Omega }\left( \left\vert \nabla \alpha \left( 
\mathbf{x},t_{i}\right) \right\vert ^{2}+\alpha ^{2}\left( \mathbf{x}%
,t_{i}\right) \right) \varphi _{\lambda }^{2}\left( z\right) d\mathbf{x}%
+\left\Vert \alpha ^{h}\right\Vert _{H_{h,k}^{4}\left( Q_{T_{1}}^{h}\right)
}^{2}.%
\end{array}%
\right.   \label{8.19}
\end{equation}%
By (\ref{6.40}), (\ref{8.6}) and (\ref{8.11}) 
\begin{equation}
\left. 
\begin{array}{c}
e^{-2\lambda \left( a-d\right) ^{2}}\displaystyle\sum \limits_{j=0}^{k}%
\displaystyle\int \limits_{\Omega }\left( \Delta \alpha \left( \mathbf{x}%
,t_{j}\right) \right) ^{2}\varphi _{\lambda }^{2}\left( z\right) d\mathbf{%
x\geq } \\ 
\geq Ce^{-2\lambda \left( a-d\right) ^{2}}\displaystyle\sum \limits%
_{j=0}^{k}\lambda \displaystyle\int \limits_{\Omega }\left( \left\vert
\nabla \alpha \left( \mathbf{x},t_{j}\right) \right\vert \right) ^{2}\varphi
_{\lambda }^{2}\left( z\right) d\mathbf{x}+ \\ 
+Ce^{-2\lambda \left( a-d\right) ^{2}}\displaystyle\sum \limits%
_{j=0}^{k}\lambda ^{3}\displaystyle\int \limits_{\Omega }\alpha ^{2}\left( 
\mathbf{x},t_{j}\right) \varphi _{\lambda }^{2}\left( z\right) d\mathbf{x,}%
\text{ }\forall \lambda \geq \lambda _{0}.%
\end{array}%
\right.   \label{8.20}
\end{equation}%
Hence, it follows from (\ref{8.19}) and (\ref{8.20}) that there exists a
sufficiently large number $\lambda _{3}=\lambda _{3}\left(
h_{0},R,T,c_{1},d,\Omega \right) \geq \lambda _{0}$ such that 
\begin{equation}
\left. 
\begin{array}{c}
J_{\lambda ,\gamma }\left( r_{2}^{h}\right) -J_{\lambda ,\gamma }\left(
r_{1}^{h}\right) -\left( J_{\lambda ,\gamma }^{\prime }\left(
r_{1}^{h}\right) ,r_{2}^{h}-r_{1}^{h}\right) \geq  \\ 
\geq C_{3}e^{-2\lambda \left( a-d\right) ^{2}}\displaystyle\sum \limits%
_{i=1}^{k}\displaystyle\int \limits_{\Omega }\left( \left\vert \nabla \alpha
\left( \mathbf{x},t_{i}\right) \right\vert ^{2}+\alpha ^{2}\left( \mathbf{x}%
,t_{i}\right) \right) \varphi _{\lambda }^{2}\left( z\right) d\mathbf{x}%
+\left\Vert \alpha ^{h}\right\Vert _{H_{h,k}^{4}\left( Q_{T_{1}}^{h}\right)
}^{2}, \\ 
\forall \lambda \geq \lambda _{3}.%
\end{array}%
\right.   \label{8.21}
\end{equation}%
Finally, (\ref{3.7}), (\ref{7.150}) and (\ref{8.21}) imply (\ref{8.8}),
which is the target estimate of this theorem.

Given (\ref{8.8}), the existence and uniqueness of the minimizer $r_{\min
,\lambda }^{h}\in B\left( R\right) $ of the functional $J_{\lambda ,\gamma
}\left( r^{h}\right) $ on the set $B\left( R\right) $ as well as inequality (%
\ref{8.9}) follow immediately from a combination of Lemma 2.1 with Theorem
2.1 of \cite{Bak}. $\square $

\subsection{Accuracy estimates }

\label{sec:8.3}

We estimate in this subsection the accuracy of the minimizer $r_{\min
,\lambda }^{h}$ of the functional $J_{\lambda ,\gamma }\left( r^{h}\right) $
on the set $B\left( R\right) $ in two cases: noiseless and noisy input data (%
\ref{4.26}).

\subsubsection{An accuracy estimate of the minimizer of the functional $J_{%
\protect\lambda ,\protect\gamma }\left( r^{h}\right) $ for the case of the
noiseless boundary data (\protect\ref{4.26})}

\label{sec:8.3.1}

Following one of the main concepts of the regularization theory, we assume
that there exists the true solution $r_{\ast }^{h}$ of system (\ref{4.25})
such that 
\begin{equation}
r_{\ast }^{h}\left( \mathbf{x},t\right) \in B^{\ast }\left( R\right)
=\left\{ 
\begin{array}{c}
r^{h}\left( \mathbf{x},t\right) \in H_{h,k}^{4}\left( Q_{T_{1}}^{h}\right)
:\left\Vert r^{h}\right\Vert _{H_{h,k}^{4}\left( Q_{T_{1}}^{h}\right) }\leq
R, \\ 
r^{h}\left( \mathbf{x},t\right) \mid _{S_{T_{1}}^{h}}=\overline{g}_{0,\ast
}^{h}\left( \mathbf{x},t\right) ,\text{ }r_{z}^{h}\left( \mathbf{x},t\right)
\mid _{\Gamma _{0,T_{1}}^{h}}=\overline{g}_{1,\ast }^{h}\left( \mathbf{x}%
,t\right) 
\end{array}%
\right\}   \label{8.22}
\end{equation}%
where $\overline{g}_{0,\ast }^{h}\left( \mathbf{x},t\right) $ and $\overline{%
g}_{1,\ast }^{h}\left( \mathbf{x},t\right) $ are noiseless boundary data.
Theorem 8.2 provides an accuracy estimate for this case.

\textbf{Theorem 8.2 (}accuracy estimate for the minimizer of $J_{\lambda
,\gamma }\left( r^{h}\right) $ in the case of noiseless boundary data). 
\emph{Let the vector function }$r_{\ast }^{h}\left( \mathbf{x},t\right) \in
B^{\ast }\left( R\right) $\emph{\ satisfies equations (\ref{4.25}). Let the
functions }$\overline{g}_{0}^{h}\left( \mathbf{x},t\right) $ \emph{and} $%
\overline{g}_{1}^{h}\left( \mathbf{x},t\right) $\emph{\ be the boundary data
in (\ref{4.26}) and let }%
\begin{equation}
\overline{g}_{0}^{h}\left( \mathbf{x},t\right) \equiv \overline{g}_{0,\ast
}^{h}\left( \mathbf{x},t\right) ,\text{ }\overline{g}_{1}^{h}\left( \mathbf{x%
},t\right) \equiv \overline{g}_{1,\ast }^{h}\left( \mathbf{x},t\right) .
\label{8.23}
\end{equation}%
\emph{Let }$\lambda _{3}=\lambda _{3}\left( h_{0},R,T,c_{1},d,\Omega \right)
\geq \lambda _{0}$\emph{\ be the number, which was found in Theorem 8.1. For 
}$\lambda \geq \lambda _{3}$\emph{, let }$r_{\min ,\lambda }^{h}\in B\left(
R\right) $\emph{\ be the minimizer of the functional }$J_{\lambda ,\gamma
}\left( r^{h}\right) ,$\emph{\ which was found in Theorem 8.1. Then the
following accuracy estimate holds}%
\begin{equation}
\left\Vert r_{\min ,\lambda }^{h}-r_{\ast }^{h}\right\Vert
_{H_{h,k}^{1}\left( Q_{T_{1}}^{h}\right) }\leq C_{3}\exp \left[ \lambda
\left( a-d\right) ^{2}\right] \sqrt{\gamma }.  \label{8.24}
\end{equation}

\textbf{Proof.} Since by (\ref{8.4}), (\ref{8.22}) and (\ref{8.23})\emph{\ }$%
B\left( R\right) =B^{\ast }\left( R\right) ,$ then both vector functions $%
r_{\min ,\lambda }^{h},r_{\ast }^{h}\in B\left( R\right) .$ Hence, we can
apply Theorem 8.1. More precisely, taking $\lambda \geq \lambda _{3}$ and
using (\ref{8.8}), we obtain 
\begin{equation}
\left. 
\begin{array}{c}
J_{\lambda ,\gamma }\left( r_{\ast }^{h}\right) -J_{\lambda ,\gamma }\left(
r_{\min ,\lambda }^{h}\right) -\left[ J_{\lambda ,\gamma }^{\prime }\left(
r_{1}^{h}\right) ,r_{\ast }^{h}-r_{\min ,\lambda }^{h}\right] \geq  \\ 
\geq C_{3}\exp \left[ -2\lambda \left( \left( a-d\right) ^{2}-\left(
B-d\right) ^{2}\right) \right] \left\Vert r_{\ast }^{h}-r_{\min ,\lambda
}^{h}\right\Vert _{H_{h,k}^{1}\left( Q_{T_{1}}^{h}\right) }^{2}.%
\end{array}%
\right.   \label{8.240}
\end{equation}%
By (\ref{8.9}) $-\left( J_{\lambda ,\gamma }^{\prime }\left(
r_{1}^{h}\right) ,r_{\ast }^{h}-r_{\min ,\lambda }^{h}\right) \leq 0.$ Since 
$-J_{\lambda ,\gamma }\left( r_{\min ,\lambda }^{h}\right) $ as well, then (%
\ref{8.240}) implies%
\begin{equation}
C_{3}\exp \left[ -2\lambda \left( a-d\right) ^{2}\right] \left\Vert r_{\ast
}^{h}-r_{\min ,\lambda }^{h}\right\Vert _{H_{h,k}^{1}\left(
Q_{T_{1}}^{h}\right) }^{2}\leq J_{\lambda ,\gamma }\left( r_{\ast
}^{h}\right) .  \label{8.241}
\end{equation}%
Now, since the vector function $r_{\ast }^{h}\left( \mathbf{x},t\right) $%
\emph{\ }satisfies equations (\ref{4.25}), then $L_{j}\left( r_{\ast
}^{h}\right) \left( \mathbf{x},t_{j}\right) =0$ for $j=0,...,k.$ Hence, by (%
\ref{8.5}) $J_{\lambda ,\gamma }\left( r_{\ast }^{h}\right) =\gamma
\left\Vert r_{\ast }^{h}\right\Vert _{H_{h,k}^{3}\left( Q_{T_{1}}^{h}\right)
}^{2}.$ Hence, (\ref{8.241}) implies (\ref{8.24}). $\square $

In the theory of Ill-Posed Problems, the regularization parameter is always
sufficiently small \cite{T}. Thus, choose $\lambda =\lambda \left( \gamma
\right) $ such that 
\begin{equation}
\exp \left[ \lambda \left( \gamma \right) \left( a-d\right) ^{2}\right] 
\sqrt{\gamma }=\gamma ^{1/4}.  \label{8.25}
\end{equation}%
Hence,%
\begin{equation}
\lambda \left( \gamma \right) =\ln \left( \gamma ^{-1/\left( 4\omega \right)
}\right) ,\text{ }\omega =\left( a-d\right) ^{2}.  \label{8.26}
\end{equation}%
Since we should have $\lambda \left( \gamma \right) \geq \lambda _{3},$ then
(\ref{8.25}) and (\ref{8.26}) imply%
\begin{equation}
\gamma \in \left( 0,e^{-4\omega \lambda _{3}}\right) .  \label{8.27}
\end{equation}%
Thus, we arrive at the following Corollary.

\textbf{Corollary 8.1.} \emph{Assume that conditions of Theorem 8.2 hold.
Suppose that the regularization parameter }$\gamma $\emph{\ satisfies (\ref%
{8.26}). Choose }$\lambda =\lambda \left( \gamma \right) $\emph{\ as in (\ref%
{8.27}). Then accuracy estimate (\ref{8.24}) becomes}%
\begin{equation*}
\left\Vert r_{\min ,\lambda }^{h}-r_{\ast }^{h}\right\Vert
_{H_{h,k}^{1}\left( Q_{T_{1}}^{h}\right) }\leq C_{3}\gamma ^{1/4}.
\end{equation*}

\textbf{Remark 8.1.} \emph{Let }$q_{h,\min ,\lambda }\left( \mathbf{x}%
\right) $\emph{\ be the function (\ref{4.28}) in the case when the vector
function }$r^{h}$\emph{\ in (\ref{4.27}) is replaced with }$r_{\min ,\lambda
}^{h},$\emph{\ and then the so obtained vector function }$p_{\min ,\lambda
}^{h}\left( \mathbf{x},t\right) $\emph{\ is used in (\ref{4.28}). We are not
providing here an accuracy estimate for the function }$q_{h,\min ,\lambda
}\left( \mathbf{x}\right) $\emph{\ because this would require an analog of
estimate (\ref{8.24}) for a stronger norm }$\left\Vert r_{\min ,\lambda
}^{h}-r_{\ast }^{h}\right\Vert _{H_{h,k}^{2}\left( Q_{T_{1}}^{h}\right) }.$%
\emph{\ This, in turn would require a modification of Carleman estimate (\ref%
{6.40}) via an incorporation of a term with }$\left( \Delta u\right) ^{2}$%
\emph{\ in its right hand side. This can be done, see, e.g. \cite[Theorem 2.5%
]{Ksurvey} for the parabolic operator. The case of the elliptic operator }$%
\Delta u$\emph{\ easily follows from \cite{Ksurvey} if assuming the
independence on }$t$\emph{. However, such an extension of (\ref{6.40}) is
rather technical, which is why we are not doing this here. The same
considerations are in place for similar results below.}

\subsubsection{An accuracy estimate of the minimizer of the functional $J_{%
\protect\lambda ,\protect\gamma }\left( r^{h}\right) $ for the case of the
noisy boundary data (\protect\ref{4.26})}

\label{sec:8.3.2}

Just as above, let $r_{\ast }^{h}$ be the true solution of system (\ref{4.25}%
) with the noiseless boundary data (\ref{8.22}). Let $\sigma \in \left(
0,1\right) $ be the level of the noise in boundary data (\ref{4.26}). We
assume the existence of two vector functions $G^{h}\left( \mathbf{x}%
,t\right) ,G_{\ast }^{h}\left( \mathbf{x},t\right) $ such that%
\begin{equation}
\left. 
\begin{array}{c}
G^{h}\in B\left( R\right) ,\text{ }G_{\ast }^{h}\mathbf{\in }B^{\ast }\left(
R\right) , \\ 
G^{h}\left( \mathbf{x},t\right) \mid _{S_{T_{1}}^{h}}=\overline{g}%
_{0}^{h}\left( \mathbf{x},t\right) ,\text{ }G_{z}^{h}\left( \mathbf{x}%
,t\right) \mid _{\Gamma _{0,T_{1}}^{h}}=\overline{g}_{1}^{h}\left( \mathbf{x}%
,t\right) , \\ 
G_{\ast }^{h}\left( \mathbf{x},t\right) \mid _{S_{T_{1}}^{h}}=\overline{g}%
_{0,\ast }^{h}\left( \mathbf{x},t\right) ,\text{ }\partial _{z}G_{\ast
}^{h}\left( \mathbf{x},t\right) \mid _{\Gamma _{0,T_{1}}^{h}}=\overline{g}%
_{1,\ast }^{h}\left( \mathbf{x},t\right) , \\ 
\left\Vert G^{h}-G_{\ast }^{h}\right\Vert _{H_{h}^{2}\left( Q_{T}^{h}\right)
}\leq \sigma .%
\end{array}%
\right.   \label{8.28}
\end{equation}%
Denote%
\begin{equation}
\left. 
\begin{array}{c}
P^{h}\left( \mathbf{x},t\right) =r^{h}\left( \mathbf{x},t\right)
-G^{h}\left( \mathbf{x},t\right) ,\text{ }\forall r^{h}\in B\left( R\right) ,
\\ 
P_{\ast }^{h}\left( \mathbf{x},t\right) =r_{\ast }^{h}\left( \mathbf{x}%
,t\right) -G_{\ast }\left( \mathbf{x},t\right) .%
\end{array}%
\right.   \label{8.29}
\end{equation}%
It follows from (\ref{8.11}), (\ref{8.28}) and (\ref{8.29}) that 
\begin{equation}
\left. 
\begin{array}{c}
P^{h}\in B_{0}\left( 2R\right) ,\text{ }\forall r^{h}\in B\left( R\right) ,
\\ 
P_{\ast }^{h}\in B_{0}\left( 2R\right) .%
\end{array}%
\right.   \label{8.30}
\end{equation}%
Consider the functional $I_{\lambda ,\gamma },$%
\begin{equation}
\left. 
\begin{array}{c}
I_{\lambda ,\gamma }:B_{0}\left( 2R\right) \rightarrow \mathbb{R}, \\ 
I_{\lambda ,\gamma }\left( P^{h}\right) =J_{\lambda ,\gamma }\left(
P^{h}+G^{h}\right) .%
\end{array}%
\right.   \label{8.31}
\end{equation}%
Obviously, an analog of Theorem 8.1 is valid for the functional $I_{\lambda
,\gamma }.$ We formulate it now. The proof is omitted since it is almost
identical with the proof of Theorem 8.1.

\textbf{Theorem 8.3. }\emph{Assume that condition (\ref{6.6}) holds. The
functional }$I_{\lambda ,\gamma }\left( P^{h}\right) $\emph{\ has the
following properties:}

\emph{1. For any pair }$\lambda ,\gamma >0$\emph{\ and for any }$P^{h}\in
B_{0}\left( 2R\right) $\emph{\ there exists Fr\'{e}chet derivative }$%
I_{\lambda ,\gamma }^{\prime }\left( P^{h}\right) \in H_{h,k}^{4,0}\left(
Q_{T_{1}}^{h}\right) $\emph{\ of the functional }$I_{\lambda ,\gamma }\left(
P^{h}\right) $\emph{\ at the point }$P^{h}.$\emph{\ This derivative is
Lipschitz continuous on }$B_{0}\left( 2R\right) .$\emph{\ }

\emph{2. Let }$\lambda _{0}=\lambda _{0}\left( d,\Omega \right) \geq 1$\emph{%
\ be the number chosen in Theorem 6.1 and let }$\lambda _{3}=\lambda
_{3}\left( h_{0},R,T_{1},d,\Omega \right) \geq \lambda _{0}$\emph{\ be the
number chosen in Theorem 8.1. Denote }

$\overline{\lambda }_{3}=\lambda _{3}\left( h_{0},2R,T,c_{1},d,\Omega
\right) .$\emph{\ Then the functional }$I_{\lambda ,\gamma }\left(
P^{h}\right) $ \emph{is strongly convex on the set }$B_{0}\left( 2R\right) $ 
\emph{for all} $\lambda \geq \overline{\lambda }_{3},$\emph{\ i.e. the
following inequality holds}%
\begin{equation}
\left. 
\begin{array}{c}
I_{\lambda ,\gamma }\left( P_{2}^{h}\right) -I_{\lambda ,\gamma }\left(
P_{1}^{h}\right) -\left( I_{\lambda ,\gamma }^{\prime }\left(
P_{1}^{h}\right) ,P_{2}^{h}-P_{1}^{h}\right) \geq  \\ 
\geq C_{3}\exp \left[ -2\lambda \left( \left( a-d\right) ^{2}-\left(
B-d\right) ^{2}\right) \right] \left\Vert P_{2}^{h}-P_{1}^{h}\right\Vert
_{H_{h,k}^{1}\left( Q_{T_{1}}^{h}\right) }^{2}+ \\ 
+\gamma \left\Vert P_{2}^{h}-P_{1}^{h}\right\Vert _{H_{h,k}^{4}\left(
Q_{T_{1}}^{h}\right) }^{2}, \\ 
\text{ }\forall P_{1}^{h},P_{2}^{h}\in B_{0}\left( 2R\right) ,\text{ }%
\forall \lambda \geq \overline{\lambda }_{3}.%
\end{array}%
\right.   \label{8.32}
\end{equation}%
\emph{3. For every }$\lambda \geq \overline{\lambda }_{3}$\emph{\ there
exists unique minimizer }$P_{\min ,\lambda }^{h}\in B_{0}\left( 2R\right) $%
\emph{\ of the functional }$I_{\lambda ,\gamma }\left( P^{h}\right) $\emph{\
on the set }$B_{0}\left( 2R\right) .$\emph{\ Furthermore, the following
inequality holds}%
\begin{equation}
\left( I_{\lambda ,\gamma }^{\prime }\left( P_{\min ,\lambda }^{h}\right)
,P^{h}-P_{\min ,\lambda }^{h}\right) \geq 0,\text{ }\forall P^{h}\in
B_{0}\left( 2R\right) .  \label{8.33}
\end{equation}

This theorem allows us to estimate the accuracy of the minimizer $P_{\min
,\lambda }^{h}$ for the case of noisy data. Denote 
\begin{equation}
B_{1}\left( R\right) =\left\{ 
\begin{array}{c}
r^{h}\left( \mathbf{x},t\right) \in H_{h,k}^{4}\left( Q_{T_{1}}^{h}\right)
:\left\Vert r^{h}\right\Vert _{H_{h,k}^{1}\left( Q_{T_{1}}^{h}\right) }\leq
R, \\ 
r^{h}\left( \mathbf{x},t\right) \mid _{S_{T_{1}}^{h}}=\overline{g}%
_{0}^{h}\left( \mathbf{x},t\right) ,\text{ }r_{z}^{h}\left( \mathbf{x}%
,t\right) \mid _{\Gamma _{0,T_{1}}^{h}}=\overline{g}_{1}^{h}\left( \mathbf{x}%
,t\right)%
\end{array}%
\right\} .  \label{8.34}
\end{equation}%
Since the norm in the space $H_{h,k}^{1}\left( Q_{T_{1}}^{h}\right) $ is
weaker than the norm in the space $H_{h,k}^{4}\left( Q_{T_{1}}^{h}\right) ,$
then, by (\ref{8.4}) and (\ref{8.34}), the set $B_{1}\left( R\right) $ is a
\textquotedblleft weaker" analog of the set $B\left( R\right) ,$ and $%
B\left( R\right) \subset B_{1}\left( R\right) .$

\textbf{Theorem 8.4 (}accuracy estimate for the minimizer of $J_{\lambda
,\gamma }\left( r^{h}\right) $ in the case of noisy boundary data). \emph{%
Assume that conditions (\ref{8.28}) and (\ref{8.29}) hold. Choose }$\lambda =%
\overline{\lambda }_{3},$\emph{\ where the number }$\overline{\lambda }_{3}$%
\emph{\ was found in Theorem 8.3. For }$\lambda \geq \overline{\lambda }_{3},
$\emph{\ let }

$P_{\min ,\lambda }^{h}\in B_{0}\left( 2R\right) $\emph{\ be the minimizer
of the functional }$I_{\lambda ,\gamma }\left( P^{h}\right) ,$\emph{\ which
was found in Theorem 8.3. Following (\ref{8.29}), denote \ }%
\begin{equation}
\overline{r}_{\min ,\lambda }^{h}=P_{\min ,\lambda }^{h}+G^{h}\emph{.\ }
\label{8}
\end{equation}%
\emph{Choose the regularization parameter }%
\begin{equation}
\gamma =\sigma ^{2}.  \label{8.340}
\end{equation}%
\emph{Then the following accuracy estimate holds} 
\begin{equation}
\left\Vert \ \overline{r}_{\min ,\overline{\lambda }_{3}}^{h}-r_{\ast
}^{h}\right\Vert _{H_{h,k}^{1}\left( Q_{T_{1}}^{h}\right) }\leq C_{3}\sigma .
\label{8.341}
\end{equation}%
\emph{Hence, if the noise level }$\sigma $\emph{\ is so small that }$%
R-C_{3}\sigma >0$\emph{\ and also if }$r_{\ast }^{h}\in B_{1}\left(
R-C_{3}\sigma \right) ,$\emph{\ then} 
\begin{equation}
\ \overline{r}_{\min ,\overline{\lambda }_{3}}^{h}\in B_{1}\left( R\right) .
\label{8.342}
\end{equation}

\textbf{Remark 8.2.} \emph{The choice of the dependence of the
regularization parameter }$\gamma $\emph{\ on the noise level }$\sigma $%
\emph{\ as in (\ref{8.340}) is a conventional one in the theory of Ill-Posed
problems \cite{T}. It follows from (\ref{8.34}), (\ref{8.341}) and (\ref%
{8.342}) that, in the case of noisy data, even though we cannot estimate the
accuracy of the minimizer of the functional }$J_{\lambda ,\gamma }\left(
r^{h}\right) $\emph{\ on the set }$B\left( R\right) ,$\emph{\ unlike the
noiseless case of Theorem 8.2, we still find an element }

$\ \overline{r}_{\min ,\overline{\lambda }_{3}}^{h}\in B_{1}\left( R\right) $%
\emph{\ in a weaker analog }$B_{1}\left( R\right) $\emph{\ of }$B\left(
R\right) $\emph{\ such that its distance from the true solution }$r_{\ast
}^{h}$\emph{\ is sufficiently small.}

\textbf{Proof of Theorem 8.4.} By the second line of (\ref{8.30}) and (\ref%
{8.33}) $-\left( I_{\lambda ,\gamma }^{\prime }\left( P_{\min ,\lambda
}^{h}\right) ,P_{\ast }^{h}-P_{\min ,\lambda }^{h}\right) \leq 0.$ Also,
obviously $-I_{\lambda ,\gamma }\left( P_{1}^{h}\right) \leq 0.$ Hence, (\ref%
{8.32}) implies%
\begin{equation}
\left. 
\begin{array}{c}
I_{\lambda ,\gamma }\left( P_{\ast }^{h}\right) \geq C_{3}\exp \left[
-2\lambda \left( \left( a-d\right) ^{2}-\left( B-d\right) ^{2}\right) \right]
\left\Vert P_{\ast }^{h}-P_{\min ,\lambda }^{h}\right\Vert
_{H_{h,k}^{1}\left( Q_{T_{1}}^{h}\right) }^{2}, \\ 
\forall \lambda \geq \overline{\lambda }_{3}.%
\end{array}%
\right.   \label{8.35}
\end{equation}%
By (\ref{8.29}) and (\ref{8.31}) 
\begin{equation*}
I_{\lambda ,\gamma }\left( P_{\ast }^{h}\right) =J_{\lambda ,\gamma }\left(
P_{\ast }^{h}+G^{h}\right) =J_{\lambda ,\gamma }\left( r_{\ast }^{h}+\left(
G^{h}-G_{\ast }^{h}\right) \right) .
\end{equation*}%
Hence, by (\ref{6.3}), (\ref{8.5}) and (\ref{8.28}) 
\begin{equation}
I_{\lambda ,\gamma }\left( P_{\ast }^{h}\right) \leq \displaystyle\sum %
\limits_{j=0}^{k}\displaystyle\int \limits_{\Omega }\left[ L_{j}\left(
r_{\ast }^{h}+\left( G^{h}-G_{\ast }^{h}\right) \right) \left( \mathbf{x}%
,t_{j}\right) \right] ^{2}\varphi _{\lambda }^{2}\left( z\right) d\mathbf{x}%
+C_{3}\gamma .  \label{8.36}
\end{equation}%
Using (\ref{4.25}), the fourth line of (\ref{8.28}) and Cauchy-Schwarz
inequality, we obtain 
\begin{equation*}
\left[ L_{j}\left( r_{\ast }^{h}+\left( G^{h}-G_{\ast }^{h}\right) \right)
\left( \mathbf{x},t_{j}\right) \right] ^{2}\leq 2\left( L_{j}\left( r_{\ast
}^{h}\right) \right) ^{2}+C_{3}\sigma ^{2}=C_{3}\sigma ^{2}.
\end{equation*}%
Combining this with (\ref{8.340}), (\ref{8.35}) and (\ref{8.36}) and taking $%
\lambda =\overline{\lambda }_{3}$, we obtain%
\begin{equation*}
\left\Vert P_{\ast }^{h}-P_{\min ,\overline{\lambda }_{3}}^{h}\right\Vert
_{H_{h,k}^{1}\left( Q_{T_{1}}^{h}\right) }\leq C_{3}\exp \left[ \overline{%
\lambda }_{3}\left( a-d\right) ^{2}\right] \sigma .
\end{equation*}%
Next, recalling (\ref{8.29}) and (\ref{8}) as well as using triangle
inequality, we obtain%
\begin{equation*}
\left. 
\begin{array}{c}
C_{3}\sigma \geq \left\Vert P_{\ast }^{h}-P_{\min ,\overline{\lambda }%
_{3}}^{h}\right\Vert _{H_{h,k}^{1}\left( Q_{T_{1}}^{h}\right) }\geq  \\ 
\geq \left\Vert \ \overline{r}_{\min ,\overline{\lambda }_{3}}^{h}-r_{\ast
}^{h}\right\Vert _{H_{h,k}^{1}\left( Q_{T_{1}}^{h}\right) }-\left\Vert
G^{h}-G_{\ast }^{h}\right\Vert _{H_{h,k}^{1}\left( Q_{T_{1}}^{h}\right)
}\geq  \\ 
\geq \left\Vert \ \overline{r}_{\min ,\overline{\lambda }_{3}}^{h}-r_{\ast
}^{h}\right\Vert _{H_{h,k}^{1}\left( Q_{T_{1}}^{h}\right) }-\sigma ,%
\end{array}%
\right. 
\end{equation*}%
which implies (\ref{8.341}). Then (\ref{8.342}) is self evident. $\square $

\subsection{Global convergence of the gradient descent method}

\label{sec:8.4}

We now assume instead of (\ref{8.22}) that 
\begin{equation}
r_{\ast }^{h}\in B^{\ast }\left( \frac{R}{3}\right) .  \label{8.37}
\end{equation}

\subsubsection{The noiseless data}

\label{sec:8.4.1}

In the case of the noiseless data, (\ref{8.23}) holds. Hence, by (\ref{8.4})
and (\ref{8.22})%
\begin{equation}
B\left( \frac{R}{3}\right) =B^{\ast }\left( \frac{R}{3}\right) .
\label{8.38}
\end{equation}%
Consider an arbitrary vector function 
\begin{equation}
r_{0}^{h}\in B\left( \frac{R}{3}\right) .  \label{8.380}
\end{equation}%
Let $\xi \in \left( 0,1\right) $ be an arbitrary number. We construct the
gradient descent method of the minimization of the functional $J_{\lambda
,\gamma }\left( r^{h}\right) $ as follows:%
\begin{equation}
r_{n}^{h}=r_{n-1}^{h}-\xi J_{\lambda ,\gamma }^{\prime }\left(
r_{n-1}^{h}\right) ,\text{ }n=1,2,...  \label{8.39}
\end{equation}%
Note that since by Theorem 8.1 $J_{\lambda ,\gamma }^{\prime }\left(
r^{h}\right) \in H_{h,k}^{4,0}\left( Q_{T_{1}}^{h}\right) $, then boundary
conditions $\overline{g}_{0}^{h}\left( \mathbf{x},t\right) $ and $\overline{g%
}_{1}^{h}\left( \mathbf{x},t\right) $ are kept unchanged for all terms $%
r_{n}^{h}$ in (\ref{8.39}).

\textbf{Theorem 8.5} (the global convergence of the gradient descent method (%
\ref{8.39}) for noiseless data). \emph{Assume that }$\sigma =0$\emph{\ in (%
\ref{8.28}). Let (\ref{8.37}) and (\ref{8.380}) hold. Let }$\lambda =\lambda
\left( \gamma \right) $\emph{\ be as in (\ref{8.26}), and let (\ref{8.27})
be valid. Let }$r_{\min ,\lambda \left( \gamma \right) }^{h}\in B\left(
R\right) $\emph{\ be the minimizer of the functional }$J_{\lambda \left(
\gamma \right) ,\gamma }\left( r^{h}\right) $\emph{\ on the set }$B\left(
R\right) ,$\emph{\ which was found in Theorem 8.1. Assume that }%
\begin{equation}
r_{\min ,\lambda \left( \gamma \right) }^{h}\in B\left( \frac{R}{3}\right) .
\label{8.41}
\end{equation}%
\emph{Then there exists a sufficiently small number }$\xi _{0}\in \left(
0,1\right) $\emph{\ such that  }%
\begin{equation}
\left\{ r_{n}^{h}\right\} _{n=1}^{\infty }\subset B\left( R\right) ,\text{ }%
\forall \xi \in \left( 0,\xi _{0}\right) .  \label{8.42}
\end{equation}%
\emph{Furthermore, for every }$\xi \in \left( 0,\xi _{0}\right) $\emph{\
there exists a number }$\theta =\theta \left( \xi \right) \in \left(
0,1\right) $\emph{\ such that }%
\begin{equation}
\left\Vert r_{\min ,\lambda \left( \gamma \right) }^{h}-r_{n}^{h}\right\Vert
_{H_{h,k}^{4}\left( Q_{T_{1}}^{h}\right) }\leq \theta ^{n}\left\Vert r_{\min
,\lambda \left( \gamma \right) }^{h}-r_{0}^{h}\right\Vert
_{H_{h,k}^{4}\left( Q_{T_{1}}^{h}\right) },\text{ }n=1,2,...  \label{8.43}
\end{equation}%
\emph{Finally, sequence (\ref{8.39}) converges to the true solution }$%
r_{\ast }^{h}$\emph{\ of system (\ref{4.25}) with the boundary data (\ref%
{8.22}) in the norm of the space }$H_{h,k}^{1}\left( Q_{T_{1}}^{h}\right) ,$%
\emph{\ and the following estimate of the convergence rate holds:} 
\begin{equation}
\left\Vert r_{n}^{h}-r_{\ast }^{h}\right\Vert _{H_{h,k}^{1}\left(
Q_{T_{1}}^{h}\right) }\leq C_{3}\gamma ^{1/4}+\theta ^{n}\left\Vert r_{\min
,\lambda \left( \gamma \right) }^{h}-r_{0}^{h}\right\Vert
_{H_{h,k}^{4}\left( Q_{T_{1}}^{h}\right) },\text{ }n=1,2,...  \label{8.44}
\end{equation}%
\textbf{Proof.} First, we note that assumption (\ref{8.41}) is a natural
one. Indeed, by Corollary 8.1 the vector function $r_{\min ,\lambda \left(
\gamma \right) }^{h}$ is close to the vector function $r_{\ast }^{h}$ in
terms of the norm $H_{h,k}^{1}\left( Q_{T_{1}}^{h}\right) ,$ and also
conditions (\ref{8.37}) and (\ref{8.38}) hold.

Next, given (\ref{8.380}), (\ref{8.39}) and (\ref{8.41}), relations (\ref%
{8.42}) and (\ref{8.43}) follow immediately from Theorem 6 of \cite{Kind}%
\emph{. }We now prove estimate (\ref{8.44}). Using Corollary 8.1\textbf{, }(%
\ref{8.43}), triangle inequality and the fact that the norm in the space $%
H_{h,k}^{1}\left( Q_{T_{1}}^{h}\right) $ is weaker than the norm in the
space $H_{h,k}^{4}\left( Q_{T_{1}}^{h}\right) ,$ we obtain%
\begin{equation*}
\left. 
\begin{array}{c}
\left\Vert r_{n}^{h}-r_{\ast }^{h}\right\Vert _{H_{h,k}^{1}\left(
Q_{T_{1}}^{h}\right) }\leq \left\Vert r_{\min ,\lambda \left( \gamma \right)
}^{h}-r_{\ast }^{h}\right\Vert _{H_{h,k}^{1}\left( Q_{T_{1}}^{h}\right)
}+\left\Vert r_{\min ,\lambda \left( \gamma \right)
}^{h}-r_{n}^{h}\right\Vert _{H_{h,k}^{1}\left( Q_{T_{1}}^{h}\right) }\leq 
\\ 
\leq C_{3}\gamma ^{1/4}+\theta ^{n}\left\Vert r_{\min ,\lambda \left( \gamma
\right) }^{h}-r_{0}^{h}\right\Vert _{H_{h,k}^{4}\left( Q_{T_{1}}^{h}\right)
}.\text{ }\square 
\end{array}%
\right. 
\end{equation*}

\subsubsection{The noisy data}

\label{sec:8.4.2}

Let $\overline{r}_{\min ,\overline{\lambda }_{3}}^{h}$ be the vector
function defined in (\ref{8}) at $\lambda =\overline{\lambda }_{3}$. And let 
$r_{\min ,\overline{\lambda }_{3}}^{h}$ be the minimizer of the functional $%
J_{\overline{\lambda }_{3},\gamma }\left( r^{h}\right) $ on the set $B\left(
R\right) ,$ which was found in Theorem 8.1. Based on (\ref{8.341}), (\ref%
{8.342}) and (\ref{8.37}), it is reasonable to impose the following analog
of assumption (\ref{8.41}) 
\begin{equation}
r_{\min ,\overline{\lambda }_{3}}^{h}=\overline{r}_{\min ,\overline{\lambda }%
_{3}}^{h}\in B\left( \frac{R}{3}\right) .  \label{8.45}
\end{equation}%
In particular, (\ref{8.341}) and (\ref{8.45}) imply 
\begin{equation}
\left\Vert \ r_{\min ,\overline{\lambda }_{3}}^{h}-r_{\ast }^{h}\right\Vert
_{H_{h,k}^{1}\left( Q_{T_{1}}^{h}\right) }\leq C_{3}\sigma .  \label{8.46}
\end{equation}%
We again consider the gradient descent method (\ref{8.39}).

\textbf{Theorem 8.6} (the global convergence of the gradient descent method (%
\ref{8.39}) for noisy data). \emph{Let (\ref{8.37}), (\ref{8.380}) and (\ref%
{8.45}) hold. Let in (\ref{8.39}) }$J_{\lambda ,\gamma }^{\prime }\left(
r_{n-1}^{h}\right) $\emph{\ be replaced with }$J_{\overline{\lambda }%
_{3},\gamma }^{\prime }\left( r_{n-1}^{h}\right) .$\emph{\ Then there exists
a sufficiently small number }$\xi _{0}\in \left( 0,1\right) $\emph{\ such
that (\ref{8.42}) is valid. Furthermore, for every }$\xi \in \left( 0,\xi
_{0}\right) $\emph{\ there exists a number }$\theta =\theta \left( \xi
\right) \in \left( 0,1\right) $\emph{\ such that the following analog of (%
\ref{8.43}) is valid }%
\begin{equation*}
\left\Vert r_{\min ,\overline{\lambda }_{3}}^{h}-r_{n}^{h}\right\Vert
_{H_{h,k}^{4}\left( Q_{T_{1}}^{h}\right) }\leq \theta ^{n}\left\Vert r_{\min
,\overline{\lambda }_{3}}^{h}-r_{0}^{h}\right\Vert _{H_{h,k}^{4}\left(
Q_{T_{1}}^{h}\right) },\text{ }n=1,2,...
\end{equation*}%
\emph{Next, suppose that (\ref{8.340}) holds. Then sequence (\ref{8.39})
converges to }$r_{\ast }^{h}$\emph{\ in the norm of the space }$%
H_{h,k}^{1}\left( Q_{T_{1}}^{h}\right) $\emph{\ and the following estimate
of the convergence rate holds:}%
\begin{equation}
\left\Vert r_{n}^{h}-r_{\ast }^{h}\right\Vert _{H_{h,k}^{1}\left(
Q_{T_{1}}^{h}\right) }\leq C_{3}\sigma +\theta ^{n}\left\Vert r_{\min ,%
\overline{\lambda }_{3}}^{h}-r_{0}^{h}\right\Vert _{H_{h,k}^{4}\left(
Q_{T_{1}}^{h}\right) },\text{ }n=1,2,...  \label{8.47}
\end{equation}

The proof of Theorem 8.6 is almost identical with the proof of Theorem 8.5
with the only difference that when proving (\ref{8.47}), one should use (\ref%
{8.46}) instead of Corollary 8.1.

\section{Convexification Numerical Method for CIP2 }

\label{sec:9}

All theorems of section 8 for CIP1 can be reformulated for CIP2 in a similar
way. Hence, we formulate in this section only Theorems 9.1 and 9.2, which
are close analogs of Theorems 8.1 and 8.2 respectively, and do not formulate
close analogs of Theorems 8.3-8.6. We omit proofs in this section since they
are quite similar with proofs of corresponding results in section 8.

\subsection{The convexification functional}

\label{sec:9.1}

Again, let $R>0$ be an arbitrary number. Note that since we have the term $%
1/w\left( \mathbf{x},0\right) =1/w\left( \mathbf{x},t_{0}\right) $ in (\ref%
{5.25}) and (\ref{5.27}), then we need to ensure that the inequality in the
last line of (\ref{6.11}) holds, where $\left\Vert \tau \left( \mathbf{x}%
\right) \right\Vert _{C^{3}\left( \overline{\Omega }\right) }\leq M_{3},$
also, see (\ref{5.60}). On the other hand, by embedding theorem 
\begin{equation}
H^{5}\left( \Omega \right) \subset C^{3}\left( \overline{\Omega }\right) ,%
\text{ }\left\Vert f\right\Vert _{C^{3}\left( \overline{\Omega }\right)
}\leq \widehat{C}\left\Vert f\right\Vert _{H^{5}\left( \Omega \right) },%
\text{ }\forall f\in H^{5}\left( \Omega \right) ,  \label{9.0}
\end{equation}%
where the number $\widehat{C}=\widehat{C}\left( \Omega \right) >0$ depends
only on the domain $\Omega .$ Hence, using considerations, which are similar
with (\ref{8.2})-(\ref{8.4}), and taking into account (\ref{5.28}), (\ref%
{5.29}) and (\ref{9.0}), we set%
\begin{equation}
\widetilde{B}\left( R\right) =\left\{ 
\begin{array}{c}
\left( w^{h}\left( \mathbf{x},t\right) ,\tau \left( \mathbf{x}\right)
\right) \in H_{h,k}^{4}\left( Q_{T_{2}}^{h}\right) \times H^{5}\left( \Omega
\right) : \\ 
\left\Vert w^{h}\right\Vert _{H_{h,k}^{4}\left( Q_{T_{2}}^{h}\right) }\leq R,%
\text{ }\left\Vert \tau \right\Vert _{H^{5}\left( \Omega \right) }\leq R, \\ 
w\left( \mathbf{x},t_{0}\right) \geq A_{0}\left( \widehat{C}R,\Omega \right)
, \\ 
w^{h}\mid _{S_{T_{2}}^{h}}=s_{0}^{h}\left( \mathbf{x},t\right) ,\text{ }%
w_{z}^{h}\mid _{\Gamma _{0,T_{2}}^{h}}=s_{1}^{h}\left( \mathbf{x},t\right) ,
\\ 
\tau \mid _{\partial \Omega }=\tau _{0}\left( \mathbf{x}\right) ,\text{ }%
\tau _{z}\mid _{\Gamma _{0}}=\tau _{1}\left( \mathbf{x}\right) 
\end{array}%
\right\} .  \label{9.1}
\end{equation}%
Hence, similarly with (\ref{8.40})%
\begin{equation*}
\left. 
\begin{array}{c}
\widetilde{B}\left( R\right) \subset C^{2}\left( \overline{Q}%
_{T_{2}}^{h}\right) \times C^{3}\left( \overline{\Omega }\right) , \\ 
\left\Vert w^{h}\right\Vert _{C^{2}\left( \overline{Q}_{T_{2}}^{h}\right)
}\leq C^{\left( 1\right) }R,\text{ }\left\Vert \tau \right\Vert
_{C^{3}\left( \overline{\Omega }\right) }\leq \widehat{C}R,\text{ }\forall
\left( w^{h},\tau \right) \in \widetilde{B}\left( R\right) ,%
\end{array}%
\right. 
\end{equation*}%
where the number $C^{\left( 1\right) }=C^{\left( 1\right) }\left(
Q_{T_{2}}^{h},\Omega \right) >0$ depends only on listed items. Denote%
\begin{equation*}
H_{h,k}^{4}\left( Q_{T_{2}}^{h}\right) \times H^{5}\left( \Omega \right) =
\end{equation*}%
\begin{equation}
=\left\{ 
\begin{array}{c}
\left( w^{h},\tau \right) :w^{h}\in H_{h,k}^{4}\left( Q_{T_{2}}^{h}\right) ,%
\text{ }\tau \in H^{5}\left( \Omega \right) , \\ 
\left\Vert \left( w^{h},\tau \right) \right\Vert _{H_{h,k}^{4}\left(
Q_{T_{2}}^{h}\right) \times H^{5}\left( \Omega \right) }^{2}=\left\Vert
w^{h}\right\Vert _{H_{h,k}^{4}\left( Q_{T_{2}}^{h}\right) }^{2}+\left\Vert
\tau \right\Vert _{H^{5}\left( \Omega \right) }^{2}%
\end{array}%
\right\} ,  \label{9.02}
\end{equation}%
\begin{equation*}
H_{h,k}^{4,0}\left( Q_{T_{2}}^{h}\right) \times H_{0}^{5}\left( \Omega
\right) =
\end{equation*}%
\begin{equation}
=\left\{ 
\begin{array}{c}
\left( w^{h},\tau \right) \in H_{h,k}^{4}\left( Q_{T_{2}}^{h}\right) \times
H^{5}\left( \Omega \right) : \\ 
w^{h}\mid _{S_{T_{2}}^{h}}=\partial _{z}w^{h}\mid _{\Gamma
_{0,T_{2}}^{h}}=\tau \mid _{\partial \Omega }=\tau _{z}\mid _{\Gamma _{0}}=0%
\end{array}%
\right\} .  \label{9.2}
\end{equation}%
As to the regularity conditions in (\ref{9.1})-(\ref{9.2}), we refer to
Remark 2.1.

Let $L_{j}^{\left( 1\right) }\left( w^{h},\tau \right) $ and $L^{\left(
2\right) }\left( w^{h}\left( \mathbf{x},t_{0}\right) ,\tau \left( \mathbf{x}%
\right) \right) $ be the Partial Differential Operators in (\ref{5.26}) and (%
\ref{5.27}) respectively, $\lambda \geq 1$ be a parameter, $\gamma \in
\left( 0,1\right) $ be the regularization parameter and $\varphi _{\lambda
}\left( z\right) $ be the CWF (\ref{6.2}). We now construct the following
analog $\widetilde{J}_{\lambda ,\gamma }\left( w^{h},\tau \right) $ of the
functional $J_{\lambda ,\gamma }\left( r^{h}\right) $ in (\ref{8.5}): 
\begin{equation}
\left. 
\begin{array}{c}
\widetilde{J}_{\lambda ,\gamma }\left( w^{h},\tau \right) :\widetilde{B}%
\left( R\right) \rightarrow \mathbb{R}, \\ 
\widetilde{J}_{\lambda ,\gamma }\left( w^{h},\tau \right) =e^{-2\lambda
\left( a-d\right) ^{2}}\displaystyle\sum \limits_{j=0}^{k}\displaystyle\int %
\limits_{\Omega }\left[ L_{j}^{\left( 1\right) }\left( w^{h},\tau \right) %
\right] ^{2}\varphi _{\lambda }^{2}\left( z\right) d\mathbf{x}+ \\ 
+e^{-2\lambda \left( a-d\right) ^{2}}\dint\limits_{\Omega }\left[ L^{\left(
2\right) }\left( w^{h}\left( \mathbf{x},t_{0}\right) ,\tau \left( \mathbf{x}%
\right) \right) \right] ^{2}\varphi _{\lambda }^{2}\left( z\right) d\mathbf{x%
}+\gamma \left\Vert \left( w^{h},\tau \right) \right\Vert
_{H_{h,k}^{4}\left( Q_{T_{2}}^{h}\right) \times H^{5}\left( \Omega \right)
}^{2}.%
\end{array}%
\right.   \label{9.3}
\end{equation}%
The following analog of Theorem 8.1 holds:

\textbf{Theorem 9.1 }(the strong convexity of the functional $\widetilde{J}%
_{\lambda ,\gamma }\left( w^{h},\tau \right) $ on the set $\widetilde{B}%
\left( R\right) $)\textbf{.} \emph{Assume that condition (\ref{6.6}) holds.
Then the functional }$\widetilde{J}_{\lambda ,\gamma }\left( w^{h},\tau
\right) $\emph{\ has the following properties:}

\emph{1. For any pair }$\lambda ,\gamma >0$\emph{\ and for any pair }$\left(
w^{h},\tau \right) \in \widetilde{B}\left( R\right) $\emph{\ there exists Fr%
\'{e}chet derivative }$\widetilde{J}_{\lambda ,\gamma }^{\prime }\left(
w^{h},\tau \right) \in H_{h,k}^{4,0}\left( Q_{T_{2}}^{h}\right) \times
H_{0}^{5}\left( \Omega \right) $\emph{\ of the functional }$\widetilde{J}%
_{\lambda ,\gamma }\left( w^{h},\tau \right) $\emph{\ at the point }$\left(
w^{h},\tau \right) .$\emph{\ This derivative is Lipschitz continuous on the
set }$\widetilde{B}\left( R\right) .$\emph{\ }

\emph{2. Let }$\lambda _{0}=\lambda _{0}\left( d,\Omega \right) \geq 1$\emph{%
\ be the number of Theorem 6.1. There exists a sufficiently large number }$%
\lambda _{4}=\lambda _{4}\left( h_{0},R,T_{2},d,\Omega \right) \geq \lambda
_{0}$\emph{\ depending only on listed parameters such that for all }$\lambda
\geq \lambda _{4}$\emph{\ the functional }$\widetilde{J}_{\lambda ,\gamma
}\left( w^{h},\tau \right) $\emph{\ is strongly convex on the set }$%
\widetilde{B}\left( R\right) ,$\emph{\ i.e. the following inequality holds}%
\begin{equation*}
\left. 
\begin{array}{c}
\widetilde{J}_{\lambda ,\gamma }\left( w_{2}^{h},\tau _{2}\right) -%
\widetilde{J}_{\lambda ,\gamma }\left( w_{1}^{h},\tau _{1}\right) -\left\{ 
\widetilde{J}_{\lambda ,\gamma }^{\prime }\left( w_{1}^{h},\tau _{1}\right)
,\left( w_{2}^{h}-w_{1}^{h},\tau _{2}-\tau _{1}\right) \right\} \geq  \\ 
\geq C_{4}\exp \left[ -2\lambda \left( \left( a-d\right) ^{2}-\left(
B-d\right) ^{2}\right) \right] \times  \\ 
\times \left( \left\Vert w_{2}^{h}-w_{1}^{h}\right\Vert _{H_{h,k}^{1}\left(
Q_{T_{2}}^{h}\right) }^{2}+\left\Vert \tau _{2}-\tau _{1}\right\Vert
_{H^{1}\left( \Omega \right) }^{2}\right) + \\ 
+\gamma \left\Vert \left( w_{2}^{h}-w_{1}^{h},\tau _{2}-\tau _{1}\right)
\right\Vert _{H_{h,k}^{4}\left( Q_{T_{2}}^{h}\right) \times H^{5}\left(
\Omega \right) }^{2},\text{ } \\ 
\forall \left( w_{1}^{h},\tau _{1}\right) ,\left( w_{2}^{h},\tau _{2}\right)
\in \widetilde{B}\left( R\right) ,\text{ }\forall \lambda \geq \lambda _{4},%
\end{array}%
\right. 
\end{equation*}%
\emph{where }$\left\{ ,\right\} $\emph{\ denotes the scalar product in }$%
H_{h,k}^{4}\left( Q_{T_{2}}^{h}\right) \times H^{5}\left( \Omega \right) ,$ 
\emph{and the number }$C_{4}=C_{4}\left( h_{0},R,T_{2},d,\Omega \right) >0$%
\emph{\ depends only on listed parameters.}

\emph{3. For every }$\lambda \geq \lambda _{4}$\emph{\ there exists unique
minimizer }$\left( w_{\min ,\lambda }^{h},\tau _{\min ,\lambda }\right) \in 
\widetilde{B}\left( R\right) $\emph{\ of the functional }$\widetilde{J}%
_{\lambda ,\gamma }\left( w^{h},\tau \right) $\emph{\ on the set }$%
\widetilde{B}\left( R\right) .$\emph{\ Furthermore, the following inequality
holds}%
\begin{equation*}
\left\{ \widetilde{J}_{\lambda ,\gamma }^{\prime }\left( w_{\min ,\lambda
}^{h},\tau _{\min ,\lambda }\right) ,\left( w^{h}-w_{\min ,\lambda
}^{h},\tau -\tau _{\min ,\lambda }\right) \right\} \geq 0,\text{ }\forall
\left( w^{h},\tau \right) \in \widetilde{B}\left( R\right) .
\end{equation*}

Below $C_{4}=C_{4}\left( h_{0},R,T_{2},d,\Omega \right) >0$\ denotes
different numbers depending only on listed parameters.

\subsection{An accuracy estimate of the minimizer of the functional $\protect%
\widetilde{J}_{\protect\lambda ,\protect\gamma }\left( w^{h},\protect\tau %
\right) $ for the case of the noiseless boundary data (\protect\ref{5.28}), (%
\protect\ref{5.29})}

\label{sec:9.2}

Assume that there exists the true solution $\left( w_{\ast }^{h},\tau _{\ast
}\right) $ of problem (\ref{5.26})-(\ref{5.29}) with noiseless boundary data 
$s_{0,\ast }^{h}\left( \mathbf{x},t\right) ,s_{1,\ast }^{h}\left( \mathbf{x}%
,t\right) ,\tau _{0,\ast }\left( \mathbf{x}\right) ,\tau _{1,\ast }\left( 
\mathbf{x}\right) $ such that 
\begin{equation*}
\left( w_{\ast }^{h},\tau _{\ast }\right) \in \widetilde{B}^{\ast }\left(
R\right) =\left\{ 
\begin{array}{c}
\left( w^{h}\left( \mathbf{x},t\right) ,\tau \left( \mathbf{x}\right)
\right) \in H_{h,k}^{4}\left( Q_{T_{2}}^{h}\right) \times H^{5}\left( \Omega
\right) : \\ 
\left\Vert w^{h}\right\Vert _{H_{h,k}^{4}\left( Q_{T_{2}}^{h}\right) }\leq R,%
\text{ }\left\Vert \tau \right\Vert _{H^{5}\left( \Omega \right) }\leq R, \\ 
w\left( \mathbf{x},t_{0}\right) \geq A_{0}\left( \widehat{C}R,\Omega \right)
, \\ 
w^{h}\mid _{S_{T_{2}}^{h}}=s_{0,\ast }^{h}\left( \mathbf{x},t\right) ,\text{ 
}w_{z}^{h}\mid _{\Gamma _{0,T_{2}}^{h}}=s_{1,\ast }^{h}\left( \mathbf{x}%
,t\right) , \\ 
\tau \mid _{\partial \Omega }=\tau _{0,\ast }\left( \mathbf{x}\right) ,\text{
}\tau _{z}\mid _{\Gamma _{0}}=\tau _{1,\ast }\left( \mathbf{x}\right) 
\end{array}%
\right\} .
\end{equation*}%
\textbf{Theorem 9.2 (}accuracy estimate for the minimizer of $\widetilde{J}%
_{\lambda ,\gamma }\left( w^{h},\tau \right) $ in the case of noiseless
boundary data). \emph{Let the vector function }$\left( w_{\ast }^{h},\tau
_{\ast }\right) \in \widetilde{B}^{\ast }\left( R\right) $\emph{\ satisfies
equations (\ref{5.26}), (\ref{5.27}). Let the vector functions }$%
s_{0}^{h}\left( \mathbf{x},t\right) ,s_{1}^{h}\left( \mathbf{x},t\right)
,\tau _{0}\left( \mathbf{x}\right) ,\tau _{1}\left( \mathbf{x}\right) $\emph{%
\ be the boundary conditions (\ref{5.28}), (\ref{5.29}). Assume that}%
\begin{equation*}
\left. 
\begin{array}{c}
s_{0}^{h}\left( \mathbf{x},t\right) =s_{0,\ast }^{h}\left( \mathbf{x}%
,t\right) ,s_{1}^{h}\left( \mathbf{x},t\right) =s_{1,\ast }^{h}\left( 
\mathbf{x},t\right) , \\ 
\tau _{0}\left( \mathbf{x}\right) =\tau _{0,\ast }\left( \mathbf{x}\right)
,\tau _{1}\left( \mathbf{x}\right) =\tau _{1,\ast }\left( \mathbf{x}\right) .%
\end{array}%
\right. 
\end{equation*}%
\emph{\ Let }$\lambda _{4}=\lambda _{4}\left( h_{0},R,T_{2},d,\Omega \right)
\geq \lambda _{0}$\emph{\ be the number which was found in Theorem 9.1. For }%
$\lambda \geq \lambda _{4}$\emph{, let }$\left( w_{\min ,\lambda }^{h},\tau
_{\min ,\lambda }\right) \in \widetilde{B}^{\ast }\left( R\right) $\emph{\
be the minimizer of the functional }$\widetilde{J}_{\lambda ,\gamma }\left(
w^{h},\tau \right) ,$\emph{\ which was found in Theorem 9.1. Then the
following accuracy estimate holds}%
\begin{equation*}
\left. \left\Vert w_{\min ,\lambda }^{h}-w_{\ast }^{h}\right\Vert
_{H_{h,k}^{1}\left( Q_{T_{2}}^{h}\right) }+\left\Vert \tau _{\min ,\lambda
}-\tau _{\ast }\right\Vert _{H^{1}\left( \Omega \right) }\leq C_{4}\exp 
\left[ \lambda \left( a-d\right) ^{2}\right] \sqrt{\gamma },\text{ }\forall
\lambda \geq \lambda _{4}.\right. 
\end{equation*}%
\emph{Furthermore, choose }$\lambda =\lambda \left( \gamma \right) $\emph{\
\ as in (\ref{8.26}) and, similarly with (\ref{8.27}) let }$\gamma \in
\left( 0,e^{-4\omega \lambda _{4}}\right) .$\emph{\ Then}%
\begin{equation*}
\left\Vert w_{\min ,\lambda \left( \gamma \right) }^{h}-w_{\ast
}^{h}\right\Vert _{H_{h,k}^{1}\left( Q_{T_{2}}^{h}\right) }+\left\Vert \tau
_{\min ,\lambda \left( \gamma \right) }-\tau _{\ast }\right\Vert
_{H^{1}\left( \Omega \right) }\leq C_{4}\gamma ^{1/4}.
\end{equation*}

\section{Numerical Studies}
\label{sec:10}
In this section we present results of our numerical studies of the convexification method of section 8 for CIP1. We use the coefficients of 3D letter shapes to demonstrate the performance of the convexity method. The test letters include $L$, $K$, $O$, $S$, and $SZ$, these letters feature relatively thin stroke structures or contain voids. Our primary objective is to verify whether the shapes of these complex targets and the values of the coefficient $q(\mathbf{x})$ can be accurately reconstructed. Blur can be reduced on a later stage via applying one of standard \textquotedblleft image cleaning" procedures. However, the latter is outside of the scope of this paper. It is worth noting that, although our theory requires a sufficient smoothness of the coefficient $q(\mathbf{x})$ (see (\ref{2.7})), our numerical studies show that it can perform well even when $q(\mathbf{x})$ is a piecewise continuous function. 

We specify the domain $\Omega \subset \mathbb{R}^3$. The available boundary data are indicated by (\ref{2.13}) and (\ref{2.14}). We have worked numerically only with the CIP1. First, we need to obtain the boundary data $g_0(\mathbf{x}, t)$ and $g_1(\mathbf{x},t)$ in (\ref{2.13}) and (\ref{2.14}). So, we need to solve the direct problem (\ref{2.1}) and(\ref{2.2}) to obtain $g_0(\mathbf{x}, t)$ and $g_1(\mathbf{x},t)$. In our numerical studies, we approximate the $\delta -$function in {\ref{2.1}} by the $C^{\infty }\left( \mathbb{R}^{n}\right) -$function 
\begin{equation}\label{10.1}
	\delta _{\xi }\left( \mathbf{x}\right) =C_{\xi }\left\{ 
	\begin{array}{c}
		\exp \left( \frac{\left\vert \mathbf{x}\right\vert ^{2}}{\left\vert \mathbf{x%
			}\right\vert ^{2}-\xi ^{2}}\right) ,\text{ }\left\vert \mathbf{x}\right\vert
		<\xi , \\ 
		0,\text{ }\left\vert \mathbf{x}\right\vert \geq \xi ,%
	\end{array}%
	\right. 
\end{equation}%
where the parameter $\xi =0.025$, and the constant $C_{\xi }>0$ is chosen
such that 
\begin{equation}\label{10.2}
	\displaystyle\int\limits_{\left\vert \mathbf{x}\right\vert <\xi }\delta
	_{\xi }\left( \mathbf{x}\right) d\mathbf{x}=1.
\end{equation}%

We employ the finite element method to generate the boundary observation data. The computations are carried out in the spatial domain $[-0.32,0.32]^3$ with a spatial mesh size $0.0025$. The time interval $[0,1]$ is discretized into 1600 uniform time steps, corresponding to a time step size of $h =6.25\times 10^{-4}$. For the inversion, the physical region $(0.08,0.16)^3$ and the retarded-time interval $[0,0.075]$ are mapped to 
\begin{equation*}
	\Omega=(1,2)^3,\qquad [0,T_1]=[0,1],
	\label{eq:ns-domain}
\end{equation*}
respectively. In our data
simulations we have chosen the function $q\left( \mathbf{x}\right) $ as%
\begin{equation}\label{10.3}
	q\left( \mathbf{x}\right) =\left\{ 
	\begin{array}{c}
		M=const.\geq 2\text{ inside of an inclusion,} \\ 
		0\text{ otherwise.}%
	\end{array}%
	\right.
\end{equation}
In our studies we have taken
\begin{equation*}
	M=3,5,10.
\end{equation*}
Random noise is added to the normalized boundary inputs $g_0,g_1$. Specifically, the noisy data are generated by
\begin{equation}\label{10.4}
	g_{0}^{\sigma }(\mathbf{x},t)=g_{0}(\mathbf{x},t)(1+\sigma \zeta ^{(0)}),
\end{equation}%
\begin{equation}\label{10.5}
	g_{1}^{\sigma }(\mathbf{x},t)=g_{1}(\mathbf{x},t)(1+\sigma \zeta ^{(1)}),
\end{equation}%
where $\sigma\in[0,1)$ is the noise level, and the values of $\zeta^{(j)}$ are sampled independently from the uniform distribution on $[-1,1]$. Thus, $\sigma=0.01$ corresponds to a pointwise relative noise bound of $1\%$.

\begin{figure}[htbp]
	\centering
	\includegraphics[width=0.8\textwidth]{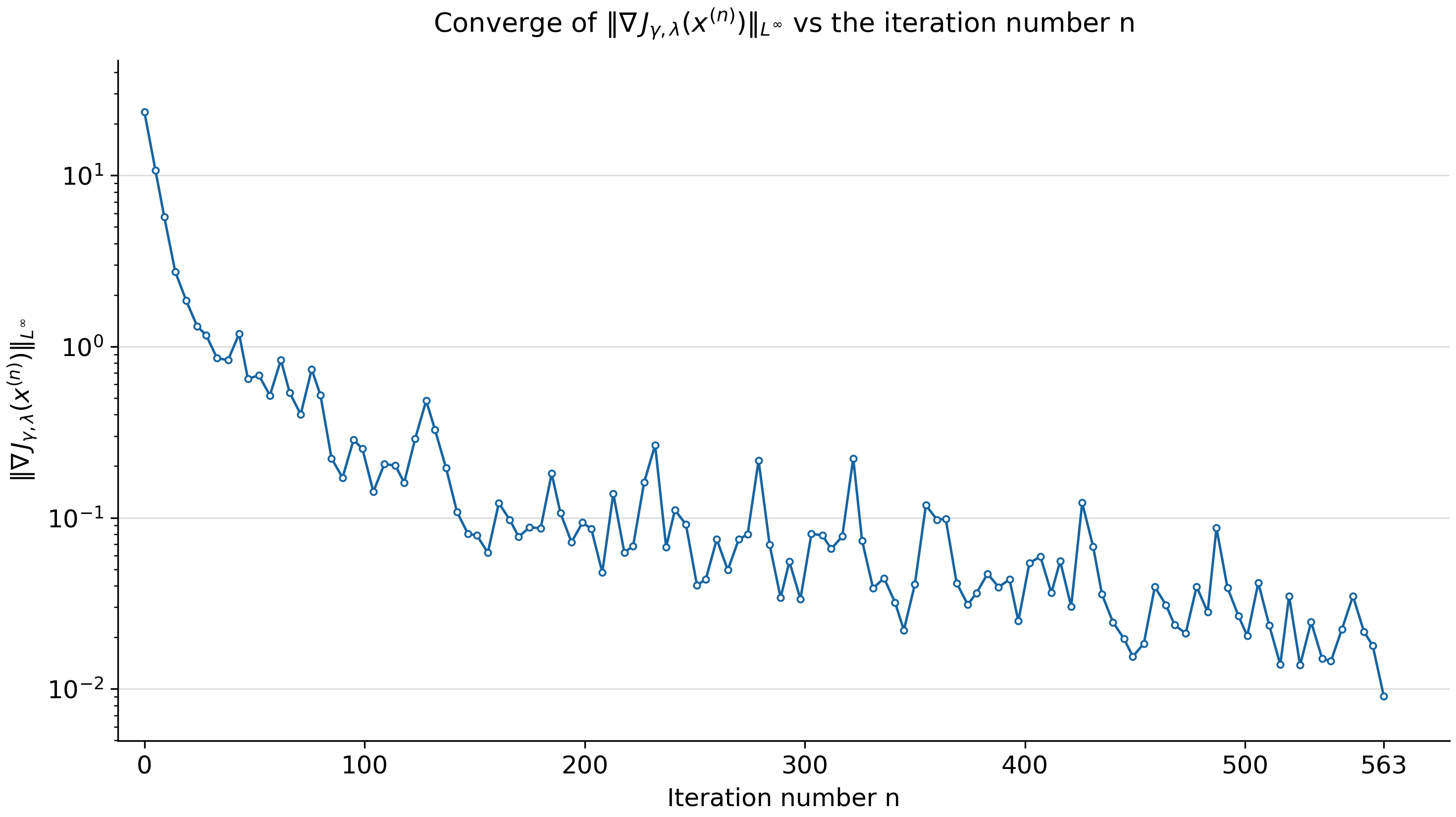}
	\caption{Convergence history of the L-BFGS-B method for the
		$L$-shaped coefficient with $N_t=70$.}
	\label{fig:L-Nt70-gradient}
\end{figure}

For the inverse problem, we use a uniform spatial grid with $19^3$ nodes and set $d=2.2$ and $\gamma=1\times10^{-5}$. The discrete functional is minimized using the L-BFGS-B method, with a stopping tolerance of $10^{-2}$ for the scaled discrete gradient, i.e.,
\begin{equation*}
	\left\|\nabla J_{\gamma,\lambda}
	\bigl(x^{(n)}\bigr)\right\|_{L^\infty}
	\leq 10^{-2}.
\end{equation*} 
Figure~\ref{fig:L-Nt70-gradient} shows the convergence history of the gradient norm for the reconstruction in Test~\ref{ex:timestep} with $N_t=70$.

For the Dirichlet boundary constraints (the first equation in \ref{4.24}), we fix them on the six boundary surfaces of the observation region. These values remain unchanged during the iteration process. We employ an elimination method to deal with the Neumann boundary constraints (the second equation in \ref{4.240}). Using transformed Dirichlet and Neumann boundary data, we rewrite the discretized Neumann boundary constraints as the equality constraints
\begin{equation*}
	\begin{aligned}
		r^h(\mathbf{x}_1,t_j)
		=\frac{1}{48}\Big[&
		25\overline{g}_0^h(\mathbf{x}_0,t_j)
		+12|\Delta z|\,
		\overline{g}_1^h(\mathbf{x}_0,t_j)\\
		&+36r^h(\mathbf{x}_2,t_j)
		-16r^h(\mathbf{x}_3,t_j)
		+3r^h(\mathbf{x}_4,t_j)\Big].
	\end{aligned}
\end{equation*}
where $\mathbf{x}_0\in\Gamma_0$ is a boundary grid point and $\mathbf{x}_\ell=\mathbf{x}_0+\ell\Delta z\,\mathbf{e}_z$,
$\ell=1,\ldots,4$, are the four successive interior grid points along the positive $z$ direction. Here $\Delta z>0$ is the spatial grid spacing and $\mathbf{e}_z=(0,0,1)^T$. During optimization, only the remaining free nodes are updated, there by automatically satisfying the discrete Neumann boundary condition at each iteration. We perform a discrete harmonic extension of the Dirichlet boundary data on the six faces to generate an initial guess. At each temporal node $t_j$, the discrete harmonic extension $\widetilde r_0^h$ is determined by
\begin{equation*}
	\begin{cases}
		\Delta_{\mathrm{FD}}\widetilde r_0^h(\mathbf{x},t_j)=0,
		& \text{at interior grid nodes},\\[2mm]
		\widetilde r_0^h(\mathbf{x},t_j)
		=\overline{g}_0^h(\mathbf{x},t_j),
		& \text{at boundary grid nodes},
	\end{cases}
	\qquad j=0,\ldots,k,
\end{equation*}
where $\Delta_{\mathrm{FD}}$ denotes the standard second-point discrete Laplace operator on the uniform spatial grid..

\begin{figure}[htbp]
	\centering
	
	\begin{subfigure}[t]{0.15\textwidth}
		\centering
		\includegraphics[width=\linewidth]{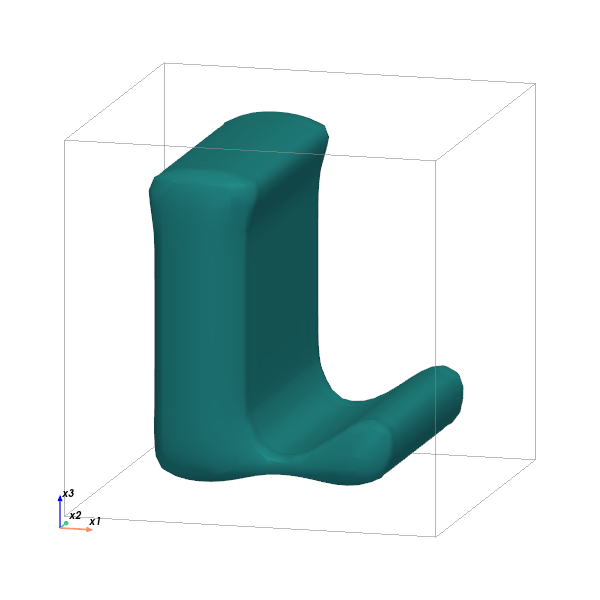}
		\caption{True coefficient}
	\end{subfigure}\hfill
	\begin{subfigure}[t]{0.15\textwidth}
		\centering
		\includegraphics[width=\linewidth]{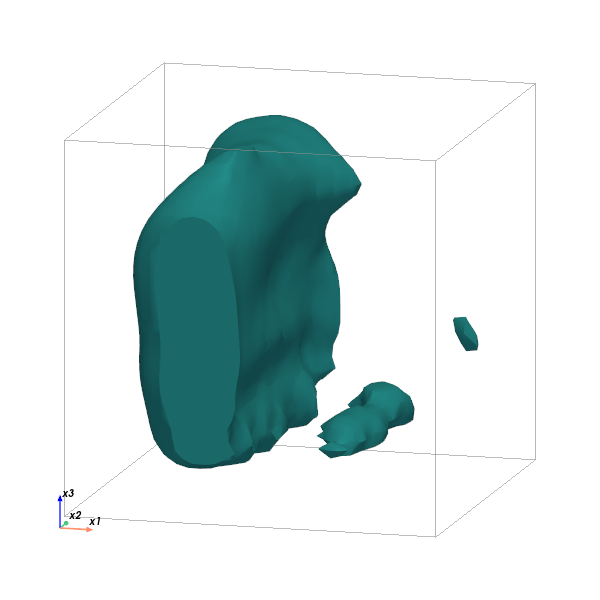}
		\caption{$N_t=10$}
	\end{subfigure}\hfill
	\begin{subfigure}[t]{0.15\textwidth}
		\centering
		\includegraphics[width=\linewidth]{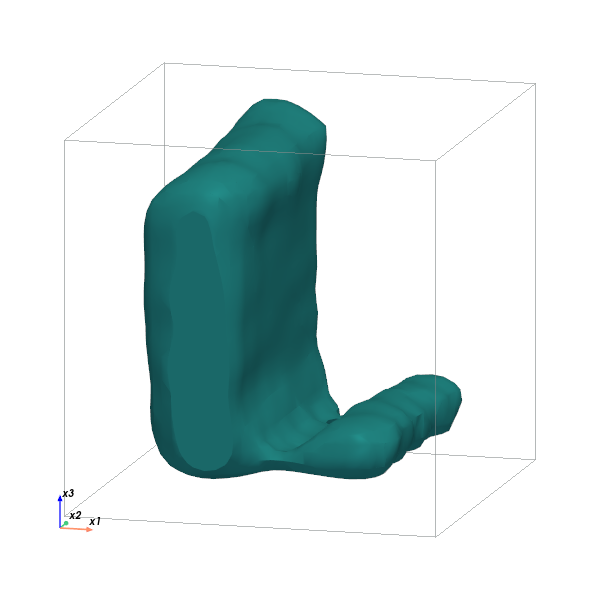}
		\caption{$N_t=30$}
	\end{subfigure}
	\begin{subfigure}[t]{0.15\textwidth}
		\centering
		\includegraphics[width=\linewidth]{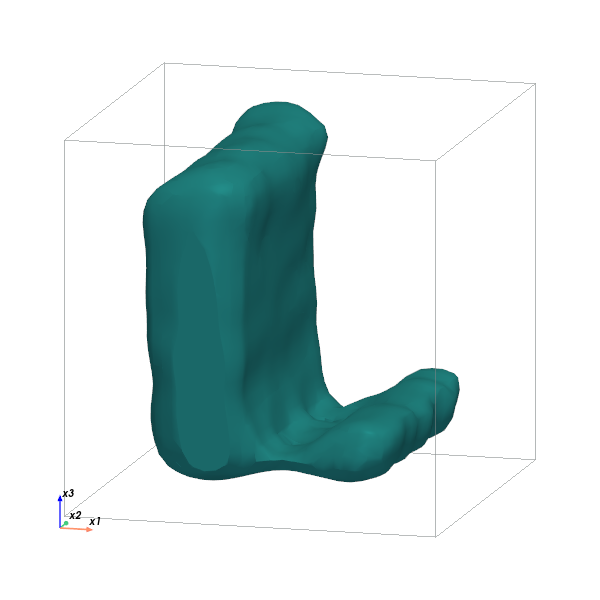}
		\caption{$N_t=50$}
	\end{subfigure}\hfill
	\begin{subfigure}[t]{0.15\textwidth}
		\centering
		\includegraphics[width=\linewidth]{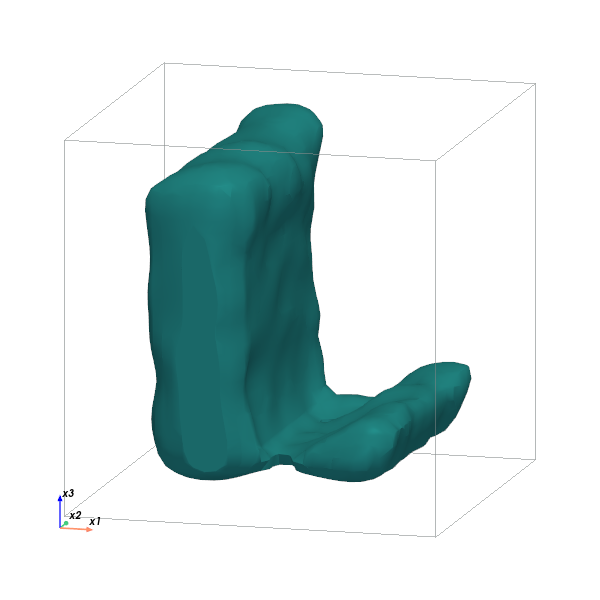}
		\caption{$N_t=70$}
	\end{subfigure}\hfill
	\begin{subfigure}[t]{0.15\textwidth}
		\centering
		\includegraphics[width=\linewidth]{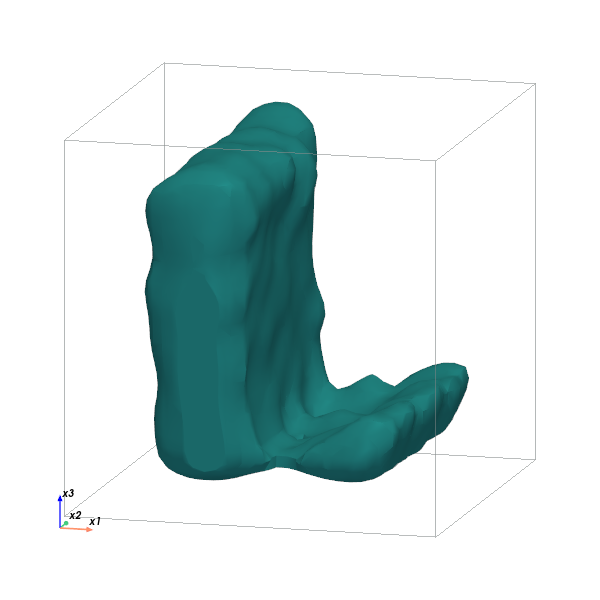}
		\caption{$N_t=100$}
	\end{subfigure}
	
	\medskip
	
	\begin{subfigure}[t]{0.15\textwidth}
		\centering
		\includegraphics[width=\linewidth]{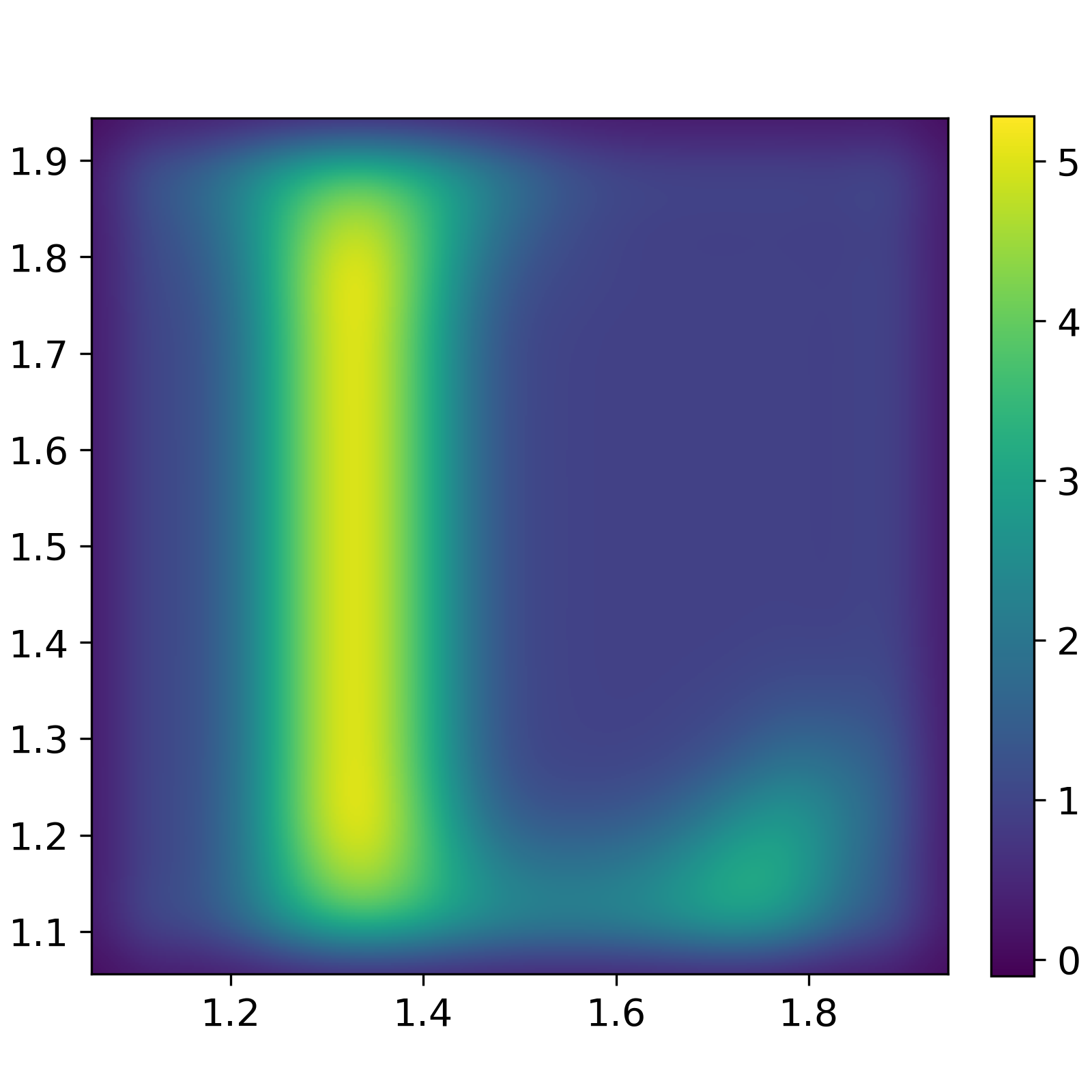}
		\caption{True coefficient}
	\end{subfigure}\hfill
	\begin{subfigure}[t]{0.15\textwidth}
		\centering
		\includegraphics[width=\linewidth]{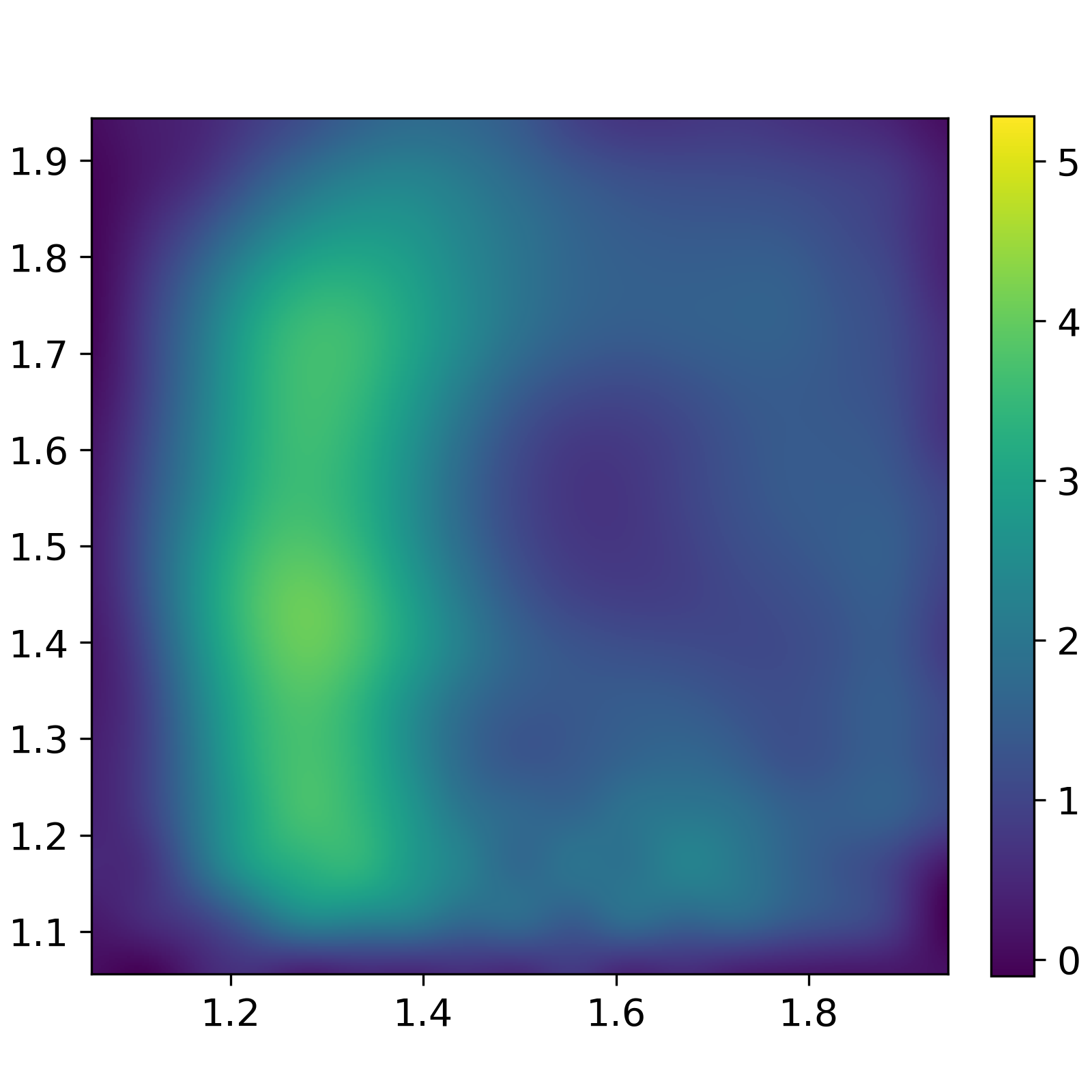}
		\caption{$N_t=10$}
	\end{subfigure}\hfill
	\begin{subfigure}[t]{0.15\textwidth}
		\centering
		\includegraphics[width=\linewidth]{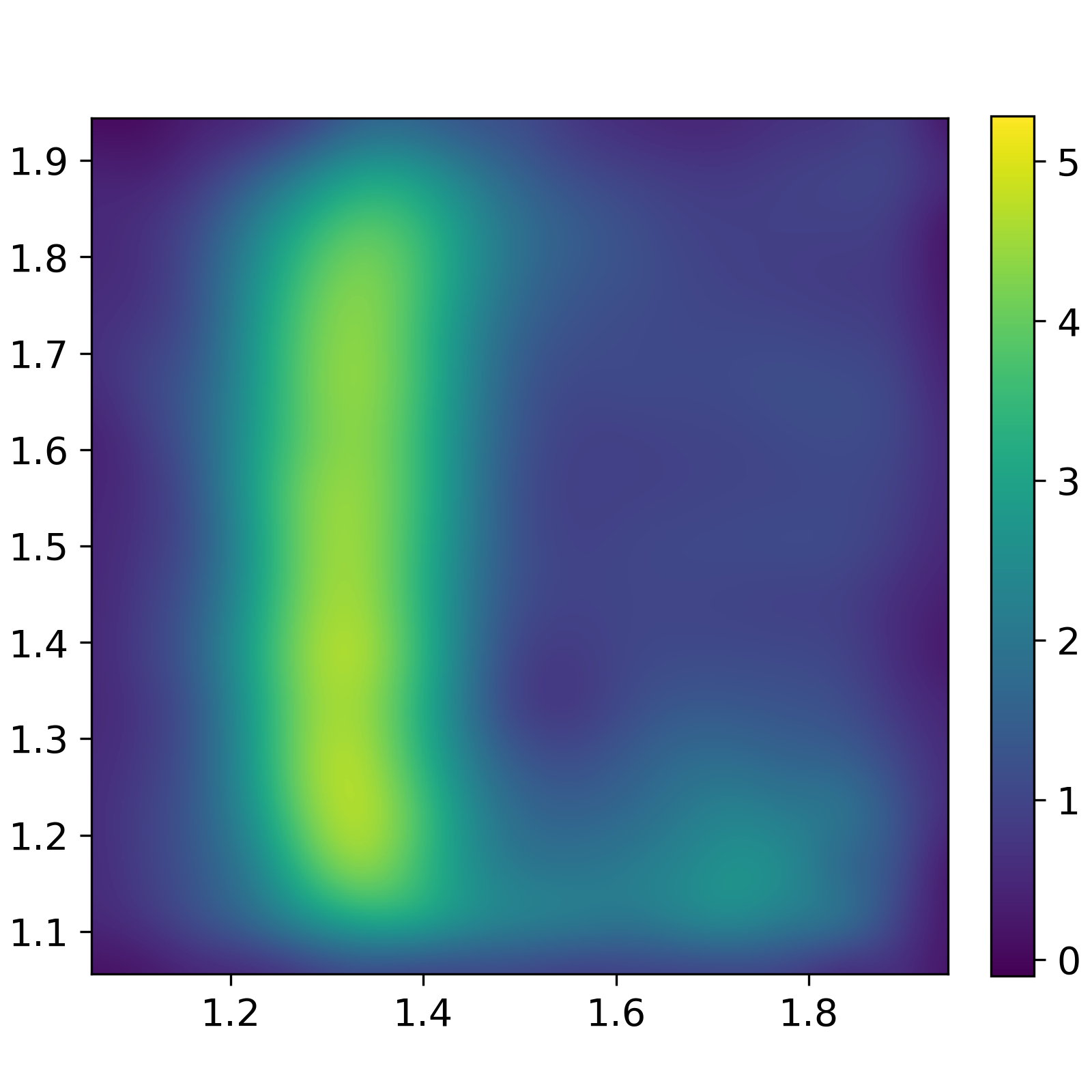}
		\caption{$N_t=30$}
	\end{subfigure}
	\begin{subfigure}[t]{0.15\textwidth}
		\centering
		\includegraphics[width=\linewidth]{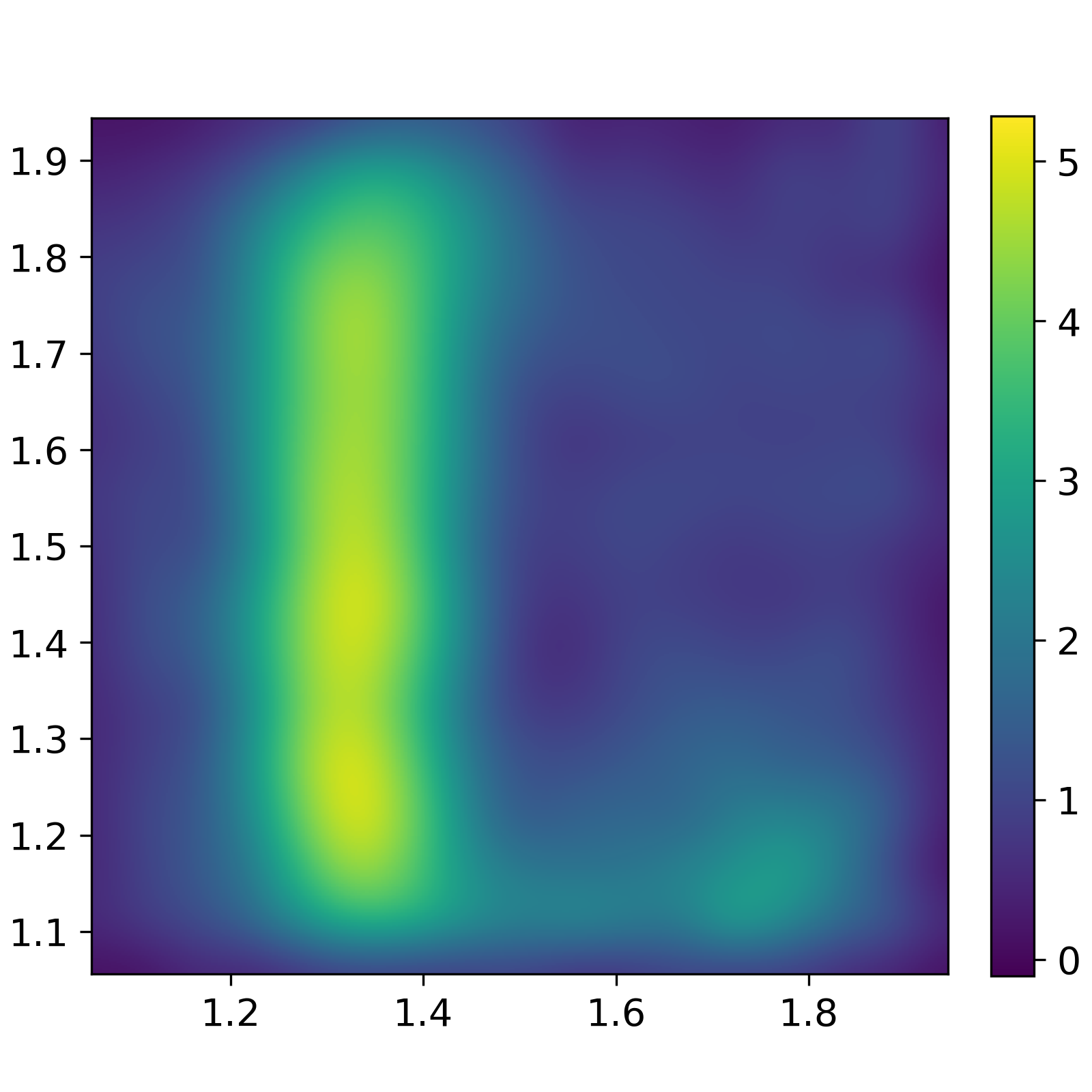}
		\caption{$N_t=50$}
	\end{subfigure}\hfill
	\begin{subfigure}[t]{0.15\textwidth}
		\centering
		\includegraphics[width=\linewidth]{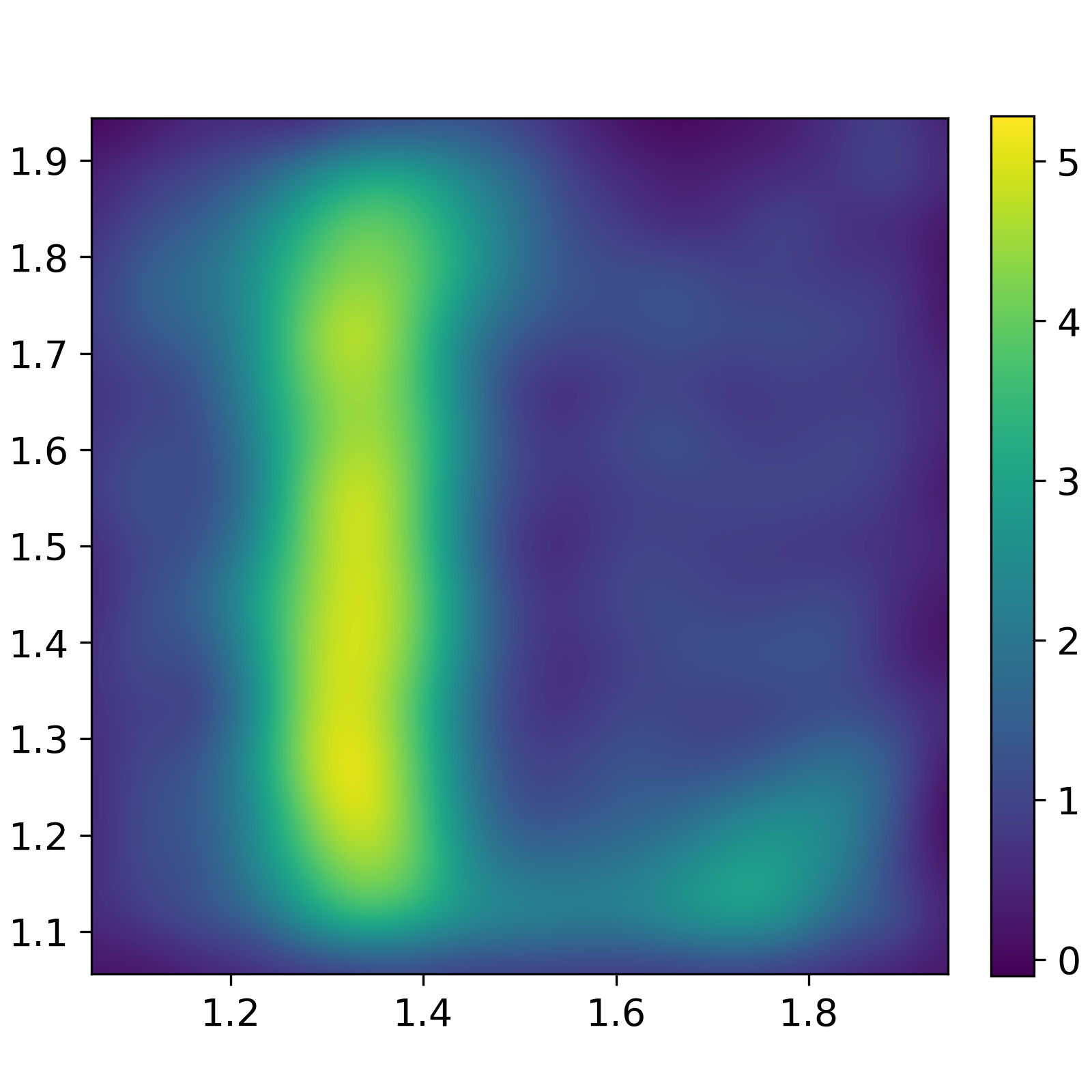}
		\caption{$N_t=70$}
	\end{subfigure}\hfill
	\begin{subfigure}[t]{0.15\textwidth}
		\centering
		\includegraphics[width=\linewidth]{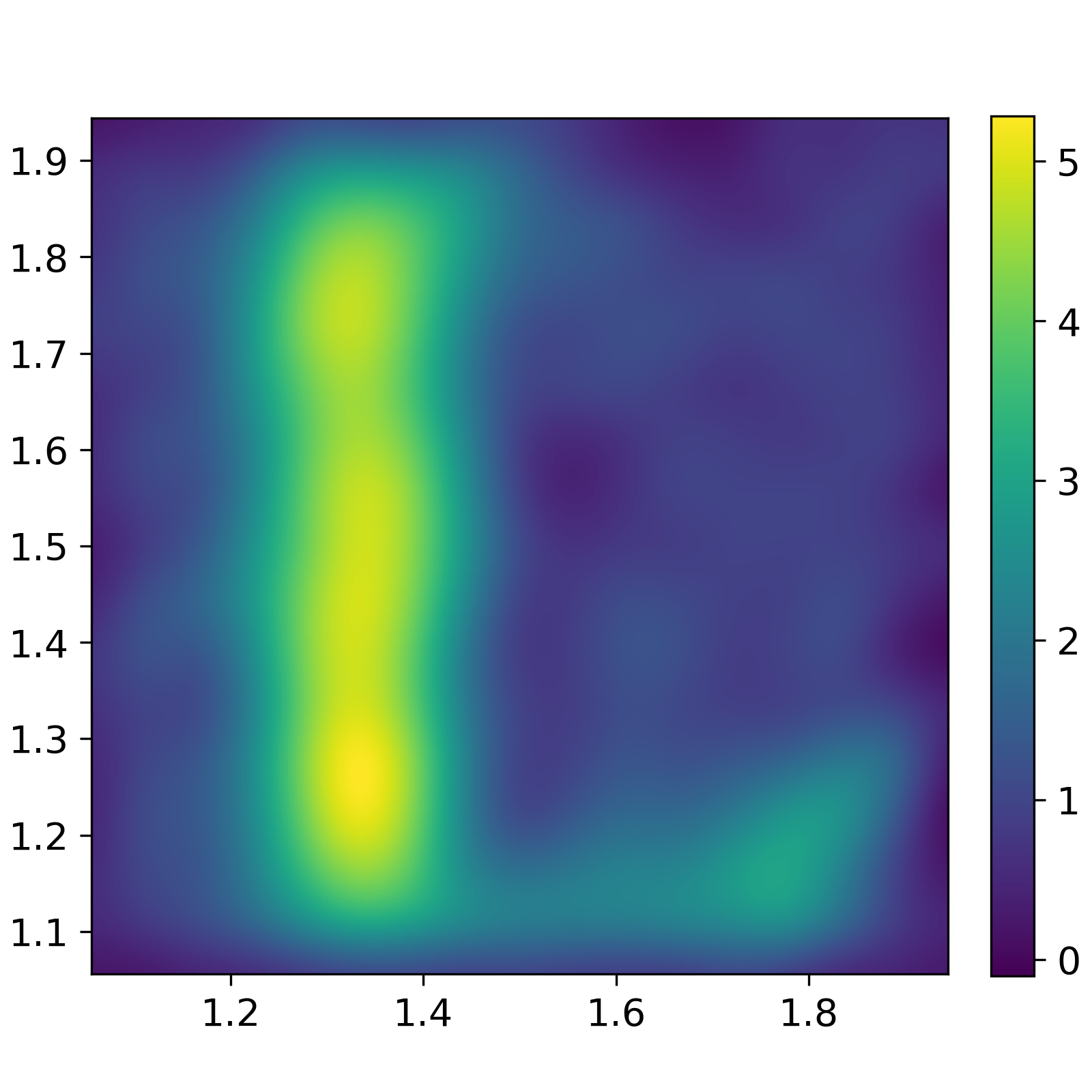}
		\caption{$N_t=100$}
	\end{subfigure}
	\caption{Test \ref{ex:timestep}:Three-dimensional and cross-sectional views of the true and reconstructed $L$-shaped coefficients for different temporal resolutions. At the tested resolutions, $N_t=70$ yielded the best reconstruction results, and was therefore used in subsequent experiments.}
	\label{fig:timestep_test}
\end{figure}

\begin{test}\label{ex:timestep}
	We investigate the impact of the temporal discretization by recovering an inclusion in the shape of the letter `$L$', for which the true coefficient is $q(\mathbf{x})=5$ inside of this inclusion and $q(\mathbf{x}) = 0$ outside of it, see (\ref{10.3}). In this test, $\lambda =2$, and the Tikhonov regularization parameter is set to $\gamma =1\times 10^{-5}$. To simulate a realistic scenario, we add $1\%$ noise to the measurement data, i.e. $\sigma =0.01$ in (\ref{10.4}), (\ref{10.5}). We perform reconstructions using five different numbers of time steps: $N_{t}\in \{10,30,50,70,100\}$. The results are
	displayed on Figure~\ref{fig:timestep_test}.
\end{test}

\begin{figure}[htbp]
	\centering
	
	\begin{subfigure}[t]{0.15\textwidth}
		\centering
		\includegraphics[width=\linewidth]{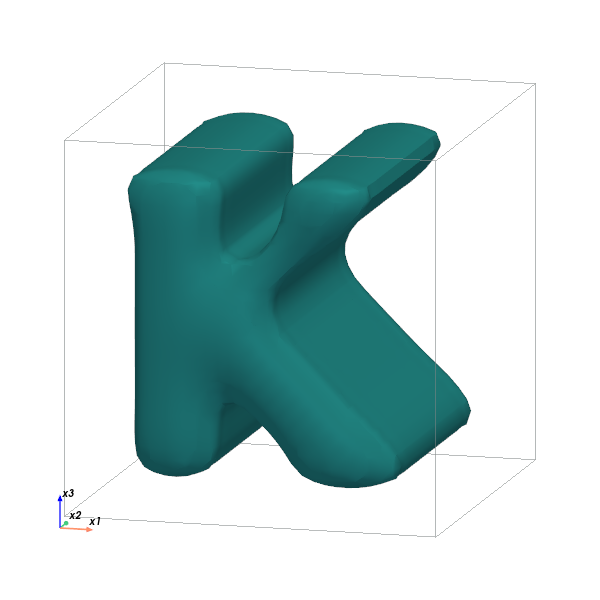}
		\caption{True coefficient}
	\end{subfigure}\hfill
	\begin{subfigure}[t]{0.15\textwidth}
		\centering
		\includegraphics[width=\linewidth]{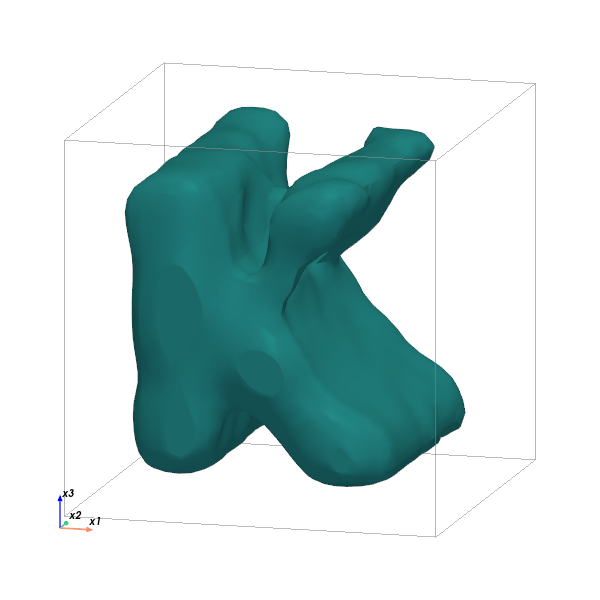}
		\caption{$\lambda=1$}
	\end{subfigure}\hfill
	\begin{subfigure}[t]{0.15\textwidth}
		\centering
		\includegraphics[width=\linewidth]{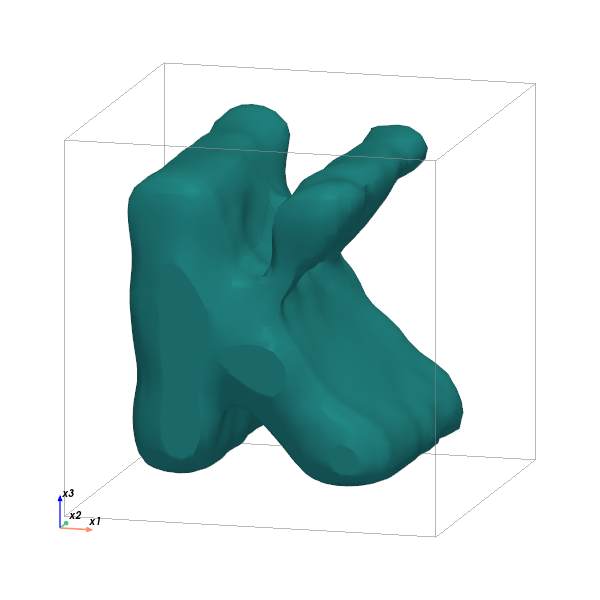}
		\caption{$\lambda=2$}
	\end{subfigure}
	\begin{subfigure}[t]{0.15\textwidth}
		\centering
		\includegraphics[width=\linewidth]{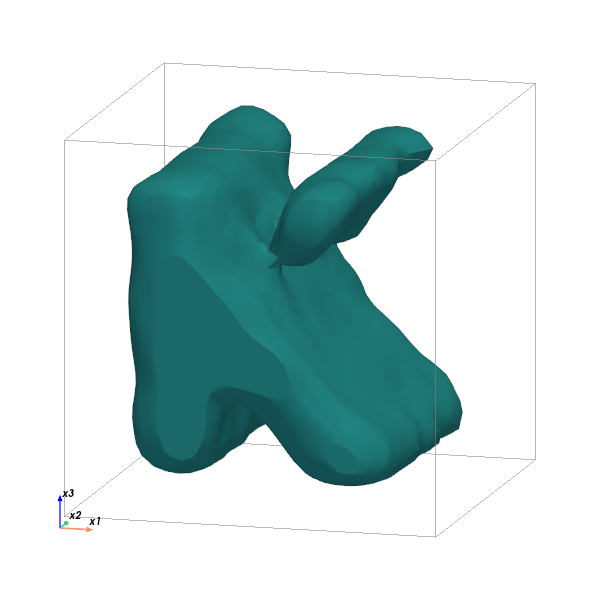}
		\caption{$\lambda=3$}
	\end{subfigure}\hfill
	\begin{subfigure}[t]{0.15\textwidth}
		\centering
		\includegraphics[width=\linewidth]{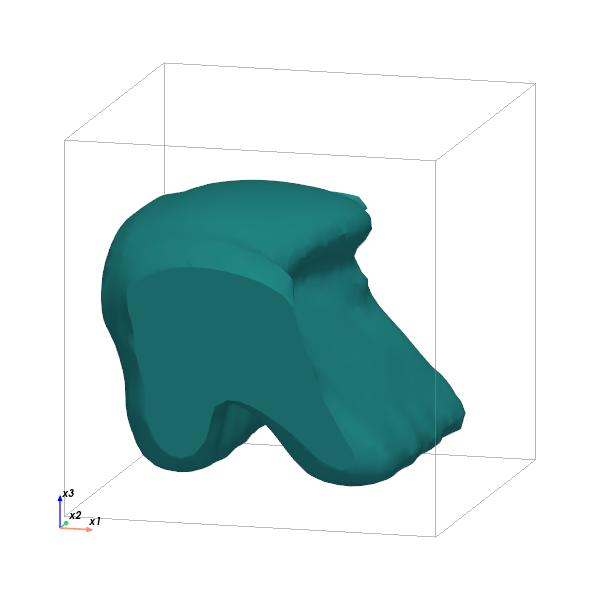}
		\caption{$\lambda=5$}
	\end{subfigure}\hfill
	\begin{subfigure}[t]{0.15\textwidth}
		\centering
		\includegraphics[width=\linewidth]{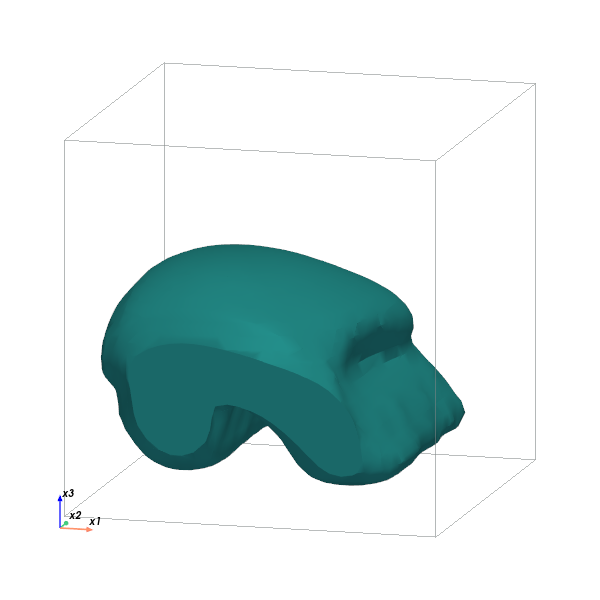}
		\caption{$\lambda=7$}
	\end{subfigure}
	
	\medskip
	
	\begin{subfigure}[t]{0.15\textwidth}
		\centering
		\includegraphics[width=\linewidth]{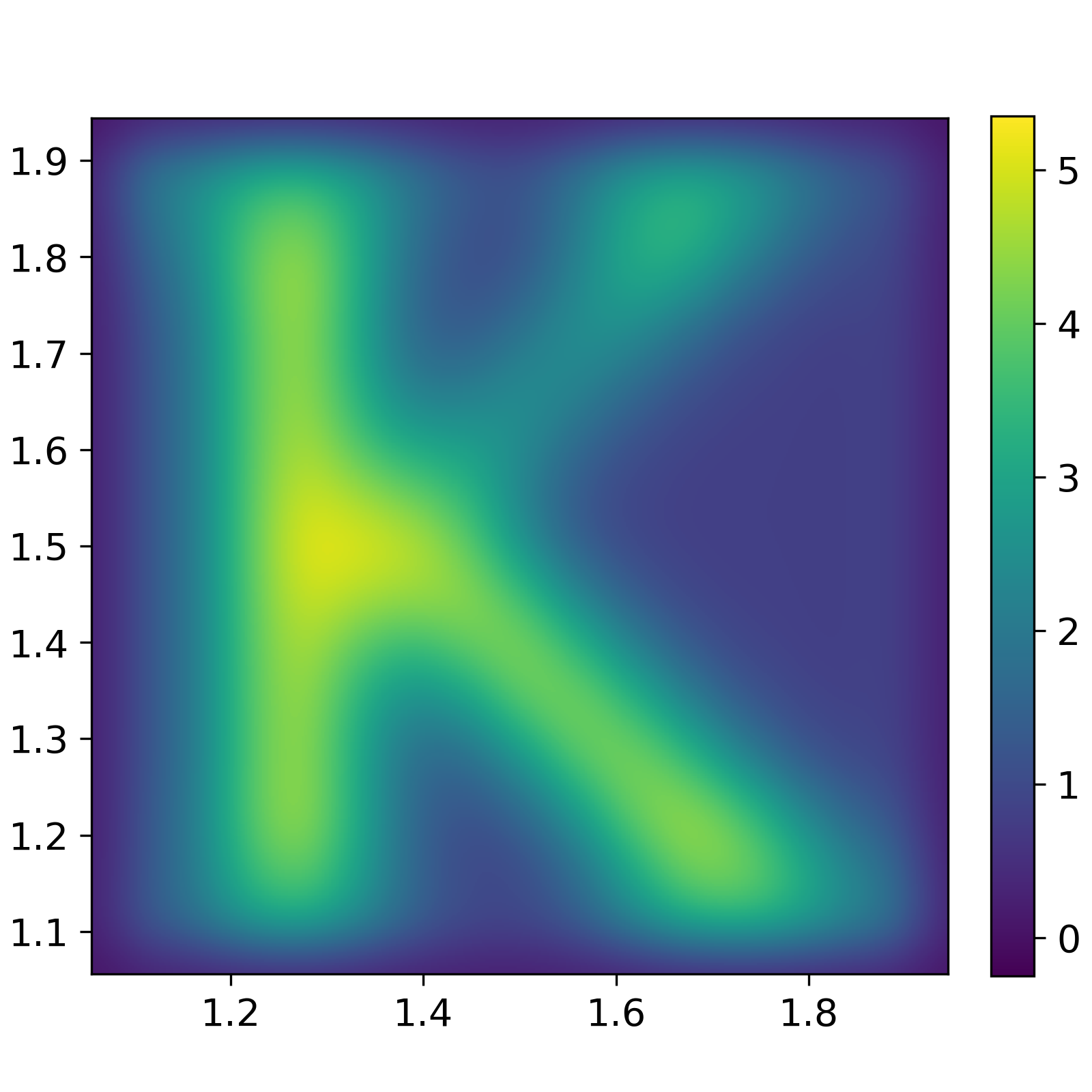}
		\caption{True coefficient}
	\end{subfigure}\hfill
	\begin{subfigure}[t]{0.15\textwidth}
		\centering
		\includegraphics[width=\linewidth]{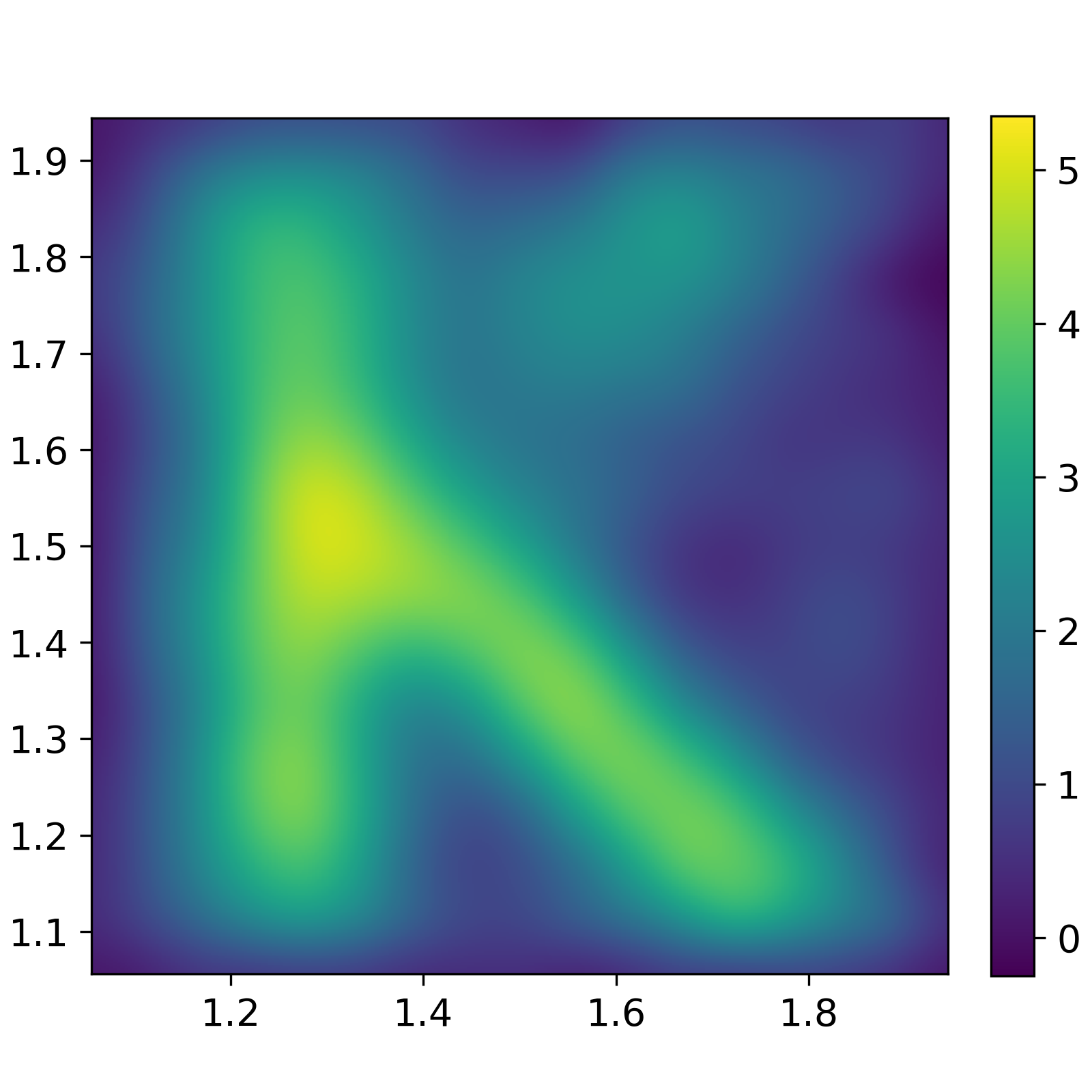}
		\caption{$\lambda=1$}
	\end{subfigure}\hfill
	\begin{subfigure}[t]{0.15\textwidth}
		\centering
		\includegraphics[width=\linewidth]{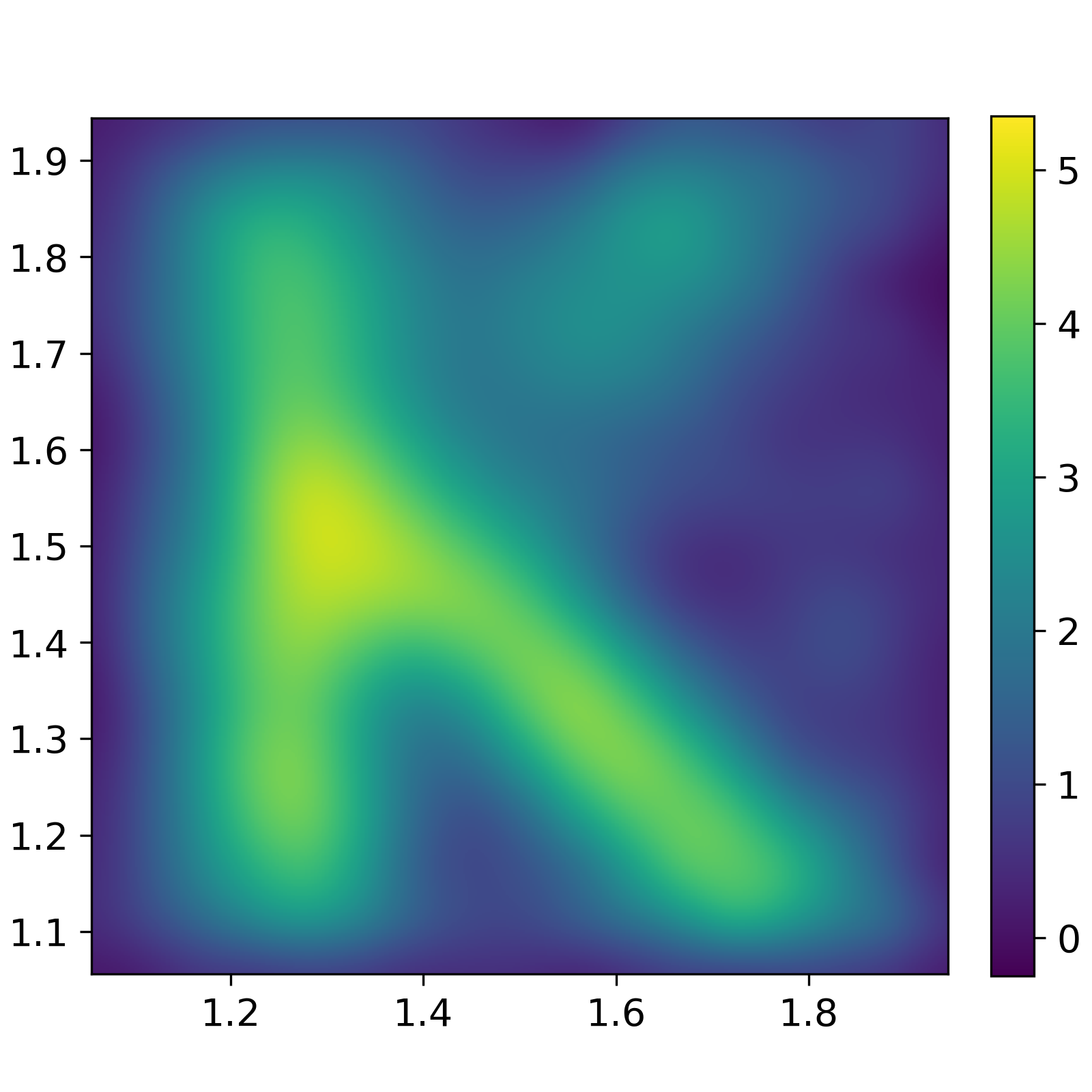}
		\caption{$\lambda=2$}
	\end{subfigure}
	\begin{subfigure}[t]{0.15\textwidth}
		\centering
		\includegraphics[width=\linewidth]{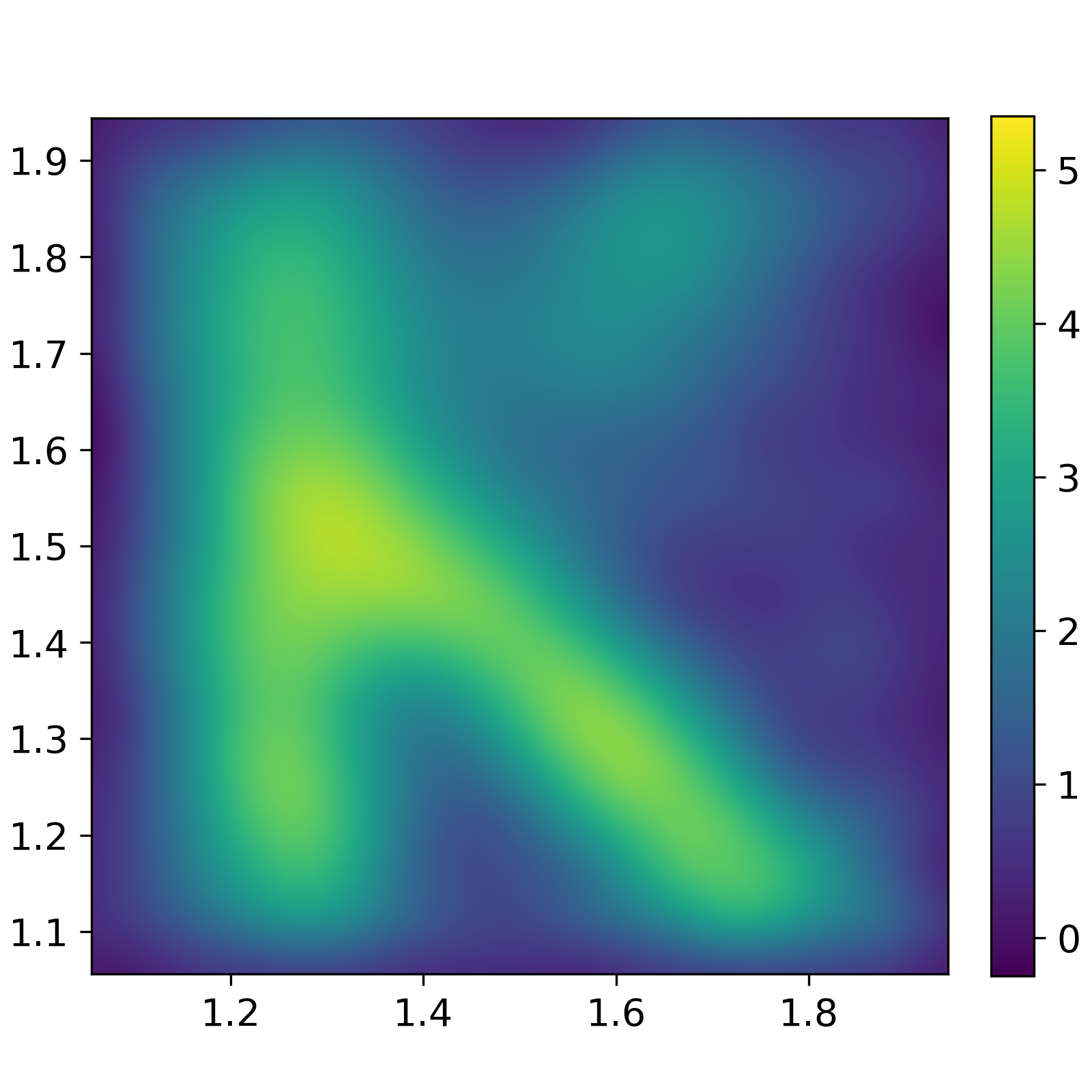}
		\caption{$\lambda=3$}
	\end{subfigure}\hfill
	\begin{subfigure}[t]{0.15\textwidth}
		\centering
		\includegraphics[width=\linewidth]{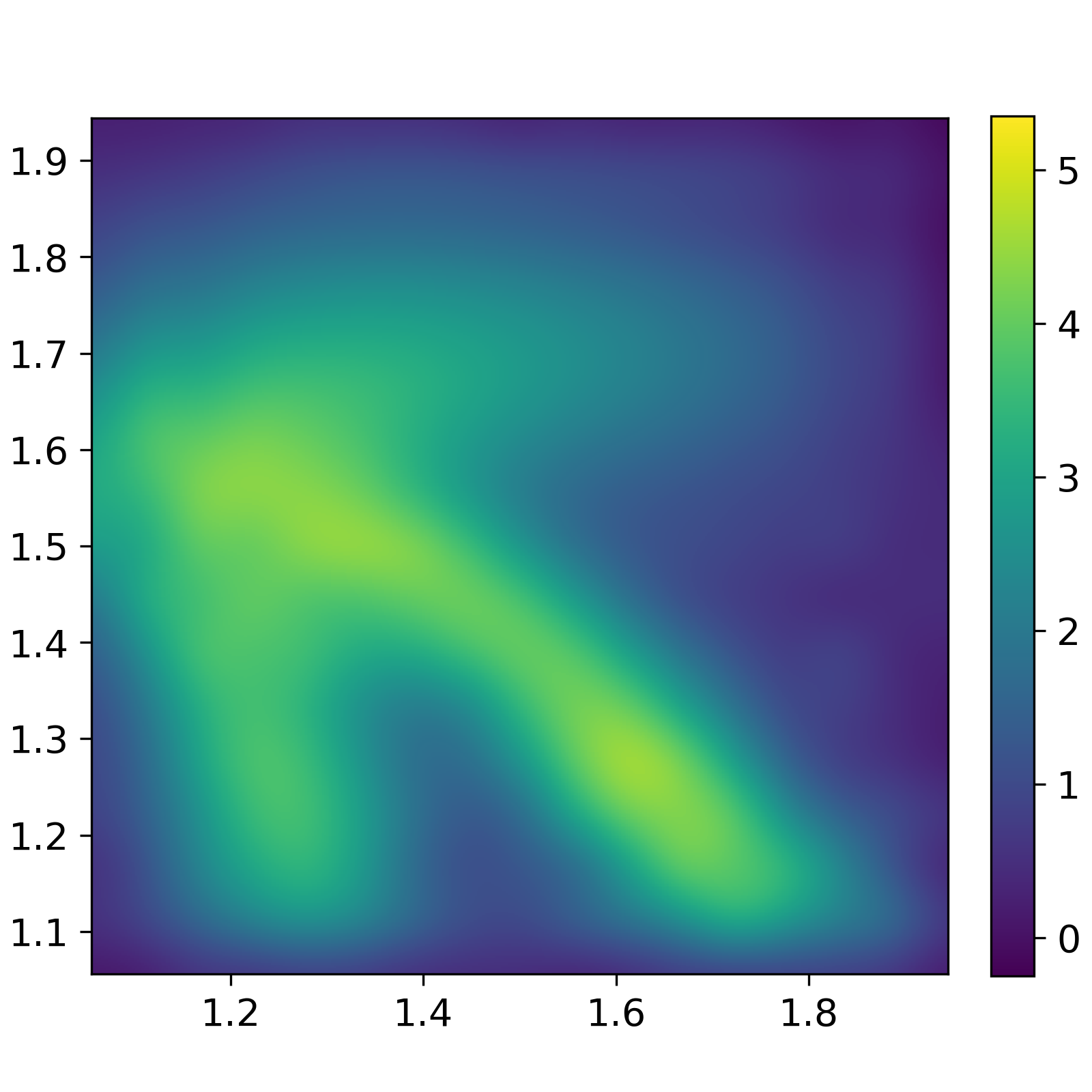}
		\caption{$\lambda=5$}
	\end{subfigure}\hfill
	\begin{subfigure}[t]{0.15\textwidth}
		\centering
		\includegraphics[width=\linewidth]{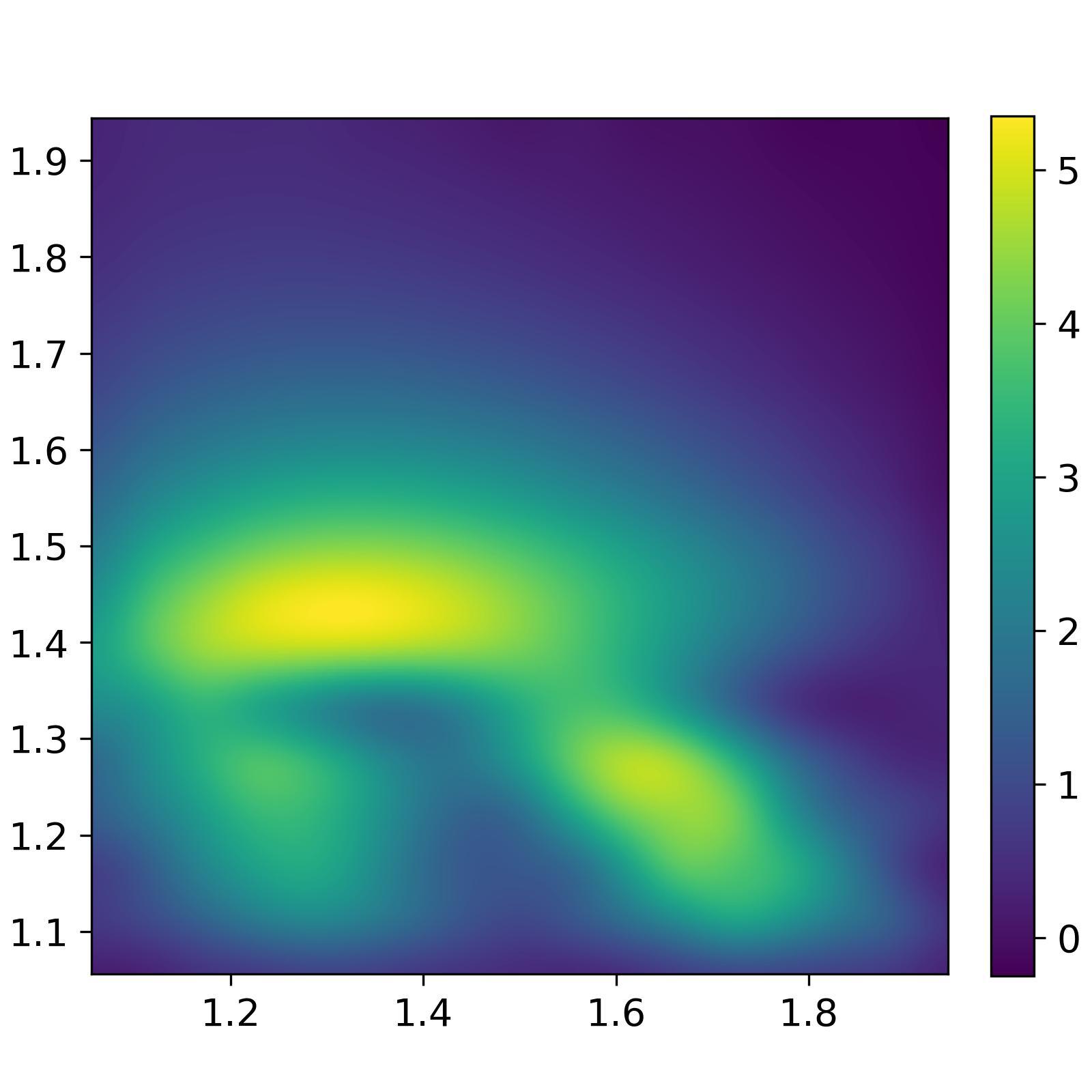}
		\caption{$\lambda=7$}
	\end{subfigure}
	
	\caption{Test \ref{ex:lambda}: Three-dimensional and cross-sectional views of the true and reconstructed $K$-shaped coefficients for different values of $\lambda$. In these tests, the reconstruction with $\lambda = 2$ performed best, and we set $\lambda = 2$ in subsequent experiments.}
	\label{fig:weight_test}
\end{figure}

\begin{test}\label{ex:lambda}
	We investigate the impact of the parameter $\lambda$ in CWF (\ref{6.2}) by recovering an inclusion in the shape of the letter `$K$', for which the true coefficient is $q(\mathbf{x})=5$ inside of this inclusion and $q(\mathbf{x}) = 0$ outside of it, see (\ref{10.3}). In this test, $N_t = 70$ and the Tikhonov regularization parameter is set to $\gamma =1\times 10^{-5}$. To simulate a realistic scenario, we add $1\%$ noise to the measurement data, i.e. $\sigma =0.01$ in (\ref{10.4}), (\ref{10.5}). We perform reconstructions using five different numbers of weight: $\lambda \in \{1,2,3,5,7\}$. The results are
	displayed on Figure~\ref{fig:weight_test}.
\end{test}

\begin{figure}[htbp]
	\centering
	
	\begin{subfigure}[t]{0.24\textwidth}
		\centering
		\includegraphics[width=\linewidth]{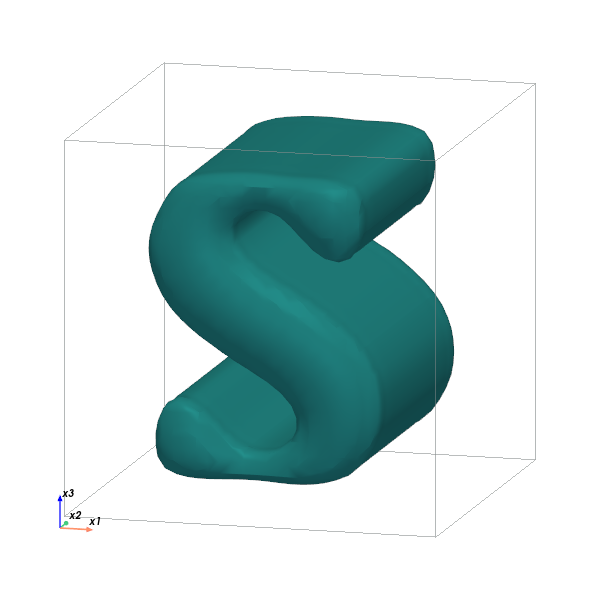}
		\caption{True profile}
	\end{subfigure}\hfill
	\begin{subfigure}[t]{0.24\textwidth}
		\centering
		\includegraphics[width=\linewidth]{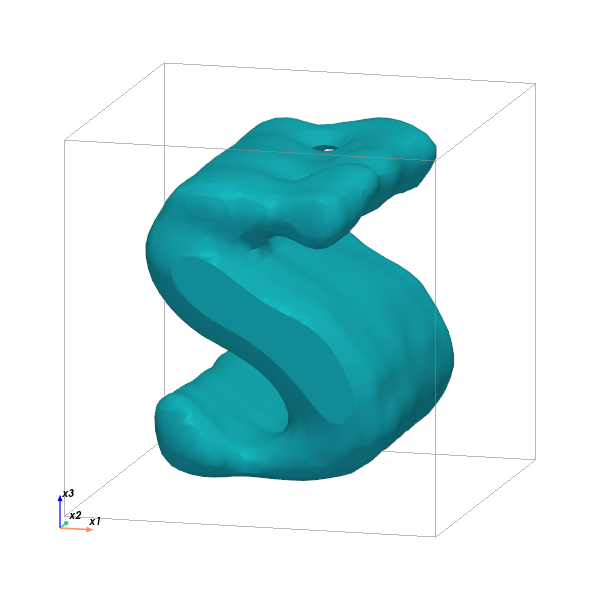}
		\caption{$M=3$}
	\end{subfigure}\hfill
	\begin{subfigure}[t]{0.24\textwidth}
		\centering
		\includegraphics[width=\linewidth]{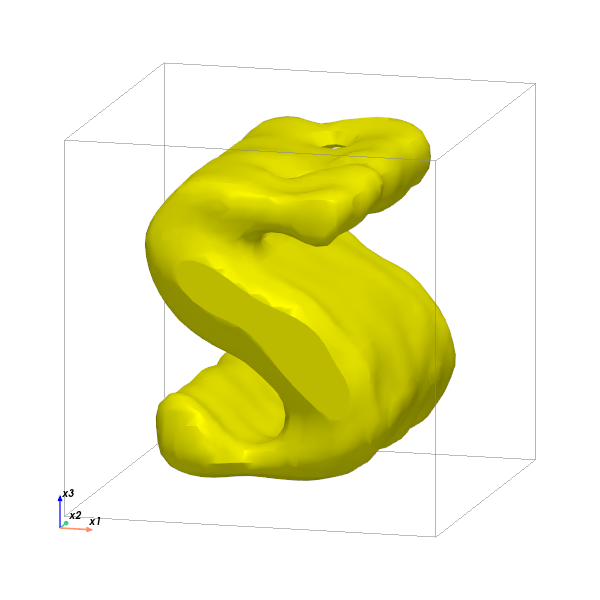}
		\caption{$M=5$}
	\end{subfigure}\hfill
	\begin{subfigure}[t]{0.24\textwidth}
		\centering
		\includegraphics[width=\linewidth]{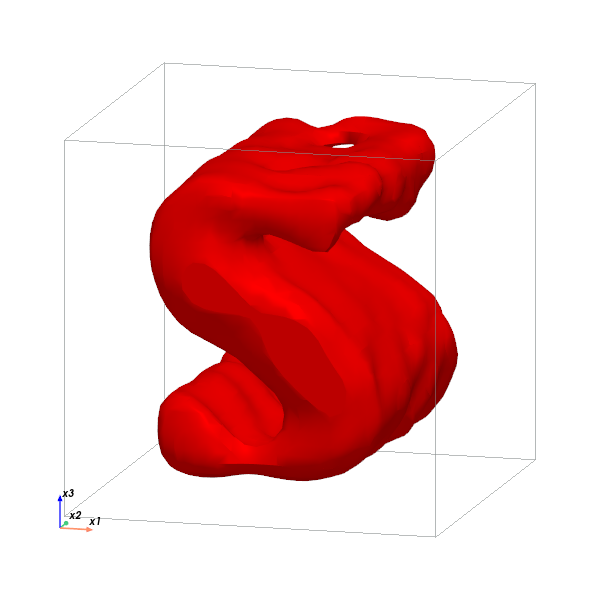}
		\caption{$M=10$}
	\end{subfigure}
	
	\par\medskip
	
	\begin{subfigure}[t]{0.24\textwidth}
		\centering
		\includegraphics[width=\linewidth]{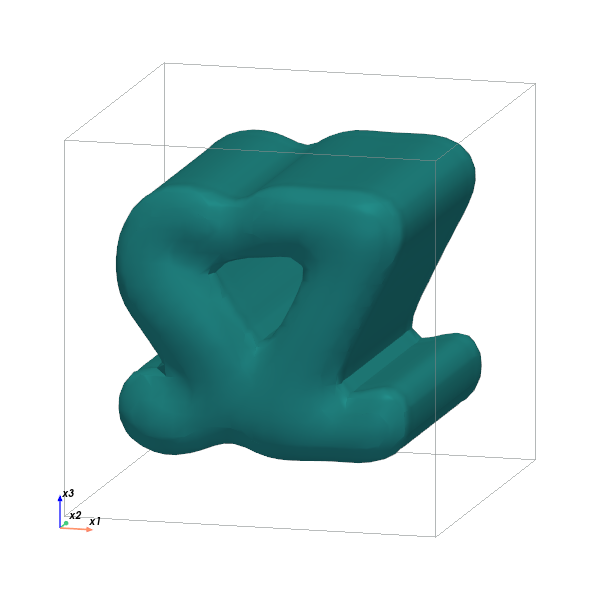}
		\caption{True profile}
	\end{subfigure}\hfill
	\begin{subfigure}[t]{0.24\textwidth}
		\centering
		\includegraphics[width=\linewidth]{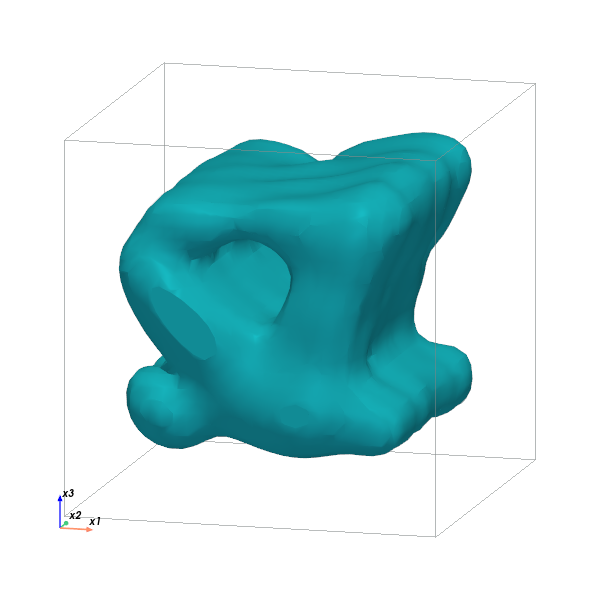}
		\caption{$M=3$}
	\end{subfigure}\hfill
	\begin{subfigure}[t]{0.24\textwidth}
		\centering
		\includegraphics[width=\linewidth]{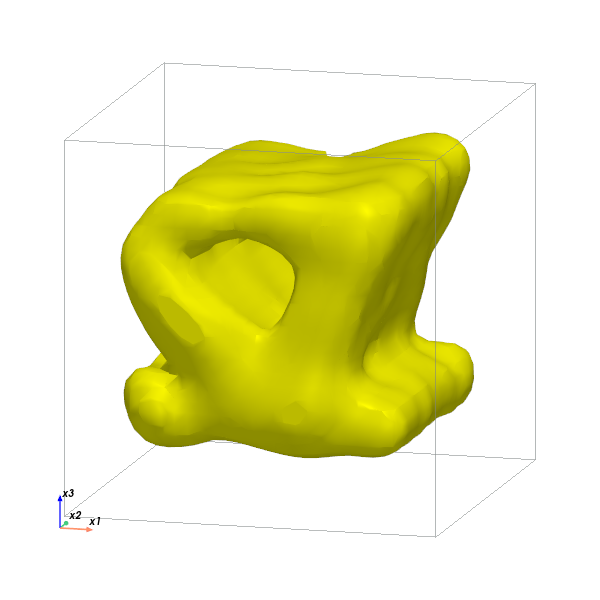}
		\caption{$M=5$}
	\end{subfigure}\hfill
	\begin{subfigure}[t]{0.24\textwidth}
		\centering
		\includegraphics[width=\linewidth]{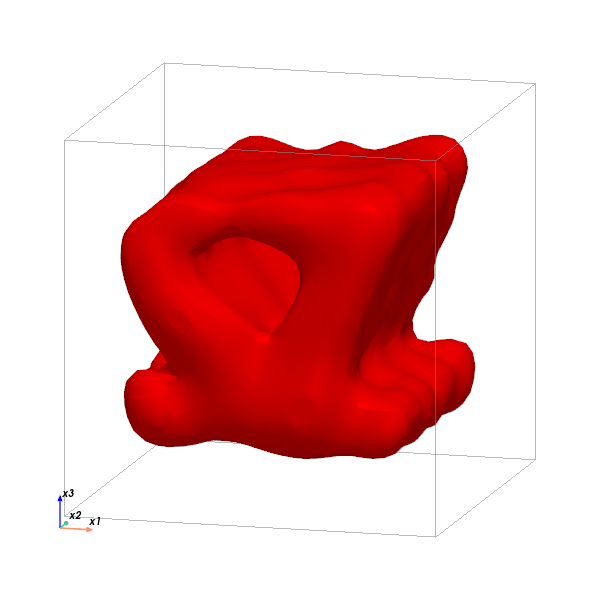}
		\caption{$M=10$}
	\end{subfigure}
	
	\caption{Test \ref{ex:values}: Three-dimensional views of the true and reconstructed
		profiles for different coefficient amplitudes.
		Top row: $S$; bottom row: $SZ$.}
	\label{fig:value_test}
\end{figure}

\begin{figure}[htbp]
	\centering
	
	\begin{subfigure}[t]{0.24\textwidth}
		\centering
		\includegraphics[width=\linewidth]{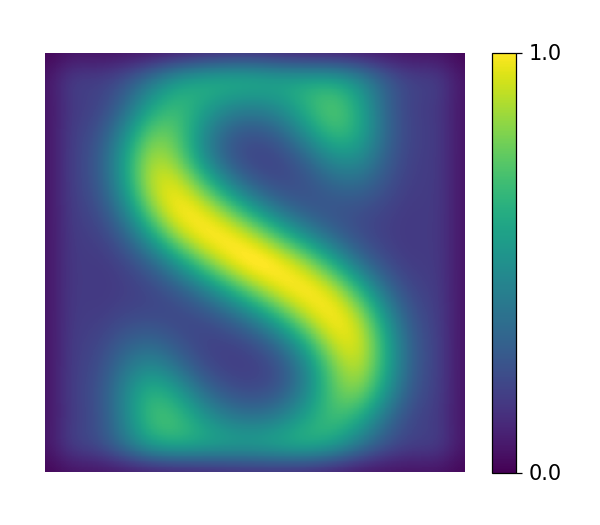}
		\caption{True profile}
	\end{subfigure}\hfill
	\begin{subfigure}[t]{0.24\textwidth}
		\centering
		\includegraphics[width=\linewidth]{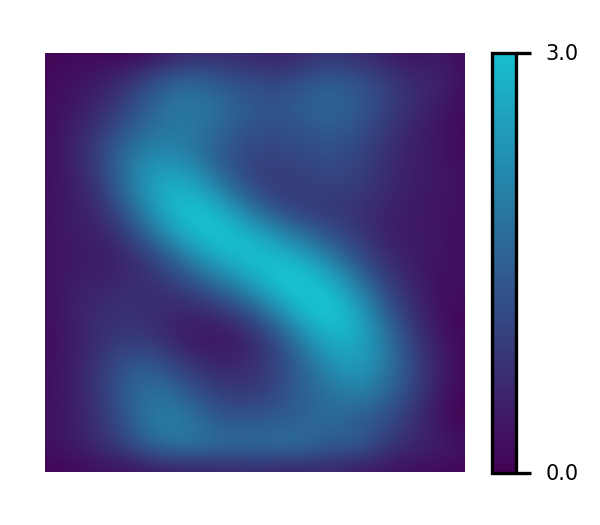}
		\caption{$M=3$}
	\end{subfigure}\hfill
	\begin{subfigure}[t]{0.24\textwidth}
		\centering
		\includegraphics[width=\linewidth]{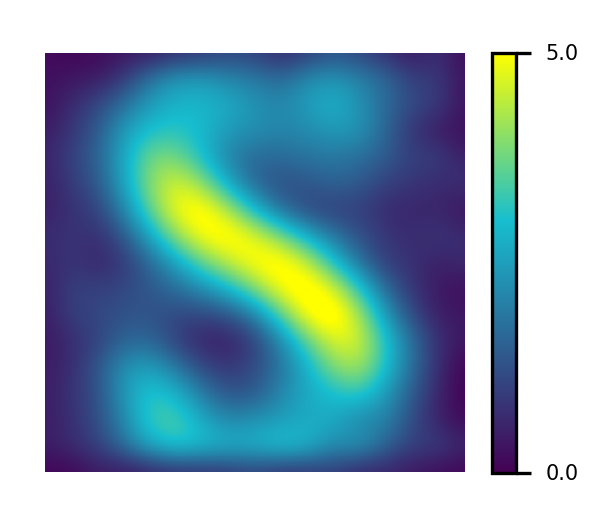}
		\caption{$M=5$}
	\end{subfigure}\hfill
	\begin{subfigure}[t]{0.24\textwidth}
		\centering
		\includegraphics[width=\linewidth]{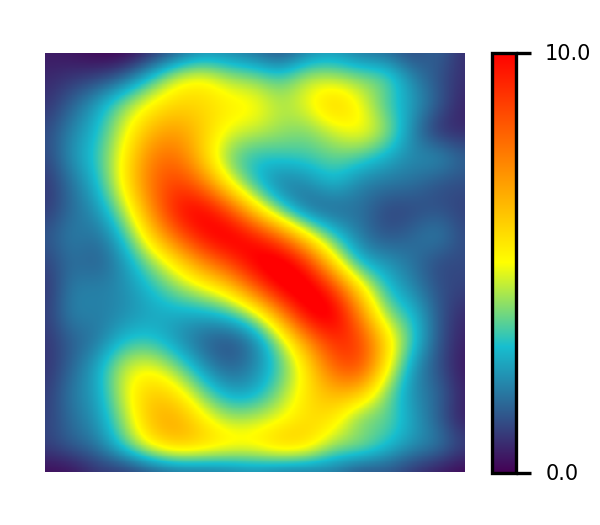}
		\caption{$M=10$}
	\end{subfigure}
	
	\par\medskip
	
	\begin{subfigure}[t]{0.24\textwidth}
		\centering
		\includegraphics[width=\linewidth]{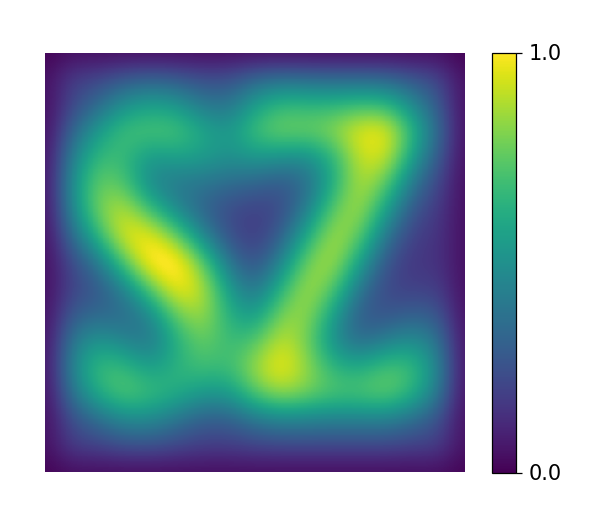}
		\caption{True profile}
	\end{subfigure}\hfill
	\begin{subfigure}[t]{0.24\textwidth}
		\centering
		\includegraphics[width=\linewidth]{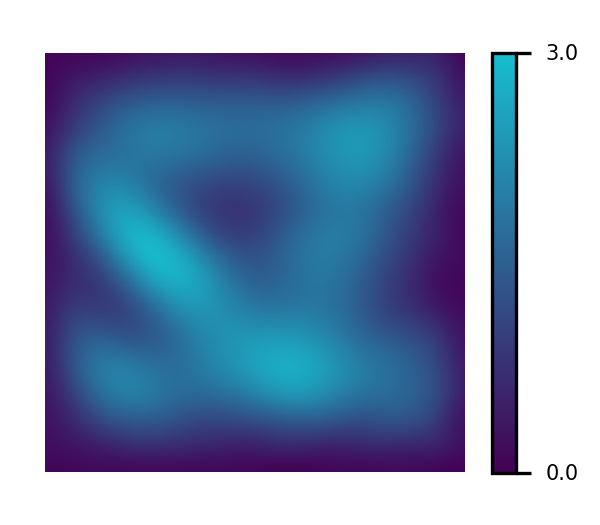}
		\caption{$M=3$}
	\end{subfigure}\hfill
	\begin{subfigure}[t]{0.24\textwidth}
		\centering
		\includegraphics[width=\linewidth]{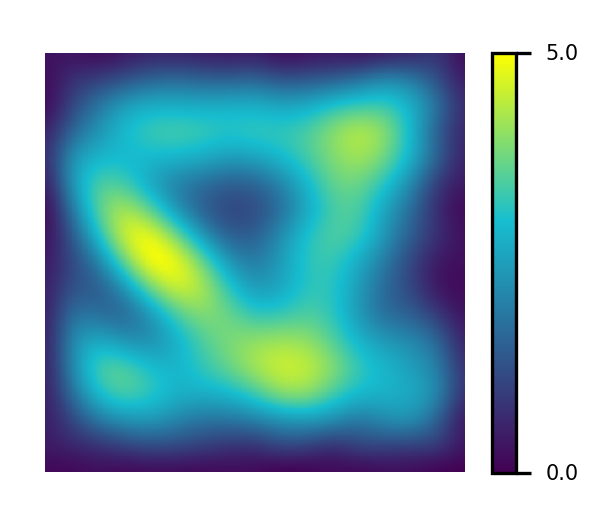}
		\caption{$M=5$}
	\end{subfigure}\hfill
	\begin{subfigure}[t]{0.24\textwidth}
		\centering
		\includegraphics[width=\linewidth]{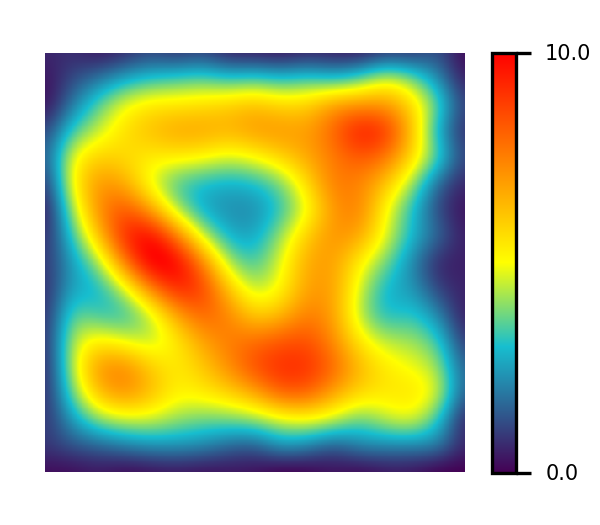}
		\caption{$M=10$}
	\end{subfigure}
	
	\caption{Test \ref{ex:values}:Cross-sectional views of the true and reconstructed coefficients $q$ for different amplitudes $M$. The true coefficients are divided by their corresponding amplitudes $M$ and presented as normalized profiles with a value range of $[0,1]$. Top row: $S$; bottom row: $SZ$.}
	\label{fig:value_test_2D}
\end{figure}

\begin{test}\label{ex:values}
	This experiment investigates the performance of our method for different values of the true coefficient $M$ in (\ref{10.3}). The shapes of inclusions are the letter '$S$' and the letters '$SZ$'. We set the true coefficient to be $M=3$, $5$, $10$, with $q(\mathbf{x}) = M$ inside of this inclusion and $q(\mathbf{x}) = 0$ outside of it.  In all cases, $1\%$ noise is added to the measurements. The reconstructions are compared to assess the sensitivity of the method to the magnitude of the unknown coefficient. The results are
	displayed on Figure~\ref{fig:value_test} and Figure~\ref{fig:value_test_2D}. In Figure~\ref{fig:value_test_2D}, the true coefficients are divided by their corresponding amplitudes $M$ and presented as normalized profiles with a value range of $[0,1]$.
\end{test}

\begin{figure}[htbp]
	\centering
	
	\begin{subfigure}[t]{0.24\textwidth}
		\centering
		\includegraphics[width=\linewidth]{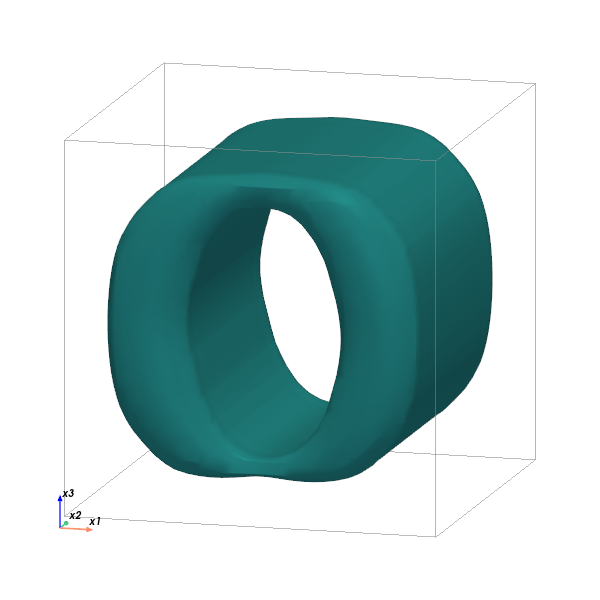}
		\caption{True coefficient}
	\end{subfigure}
	\begin{subfigure}[t]{0.24\textwidth}
		\centering
		\includegraphics[width=\linewidth]{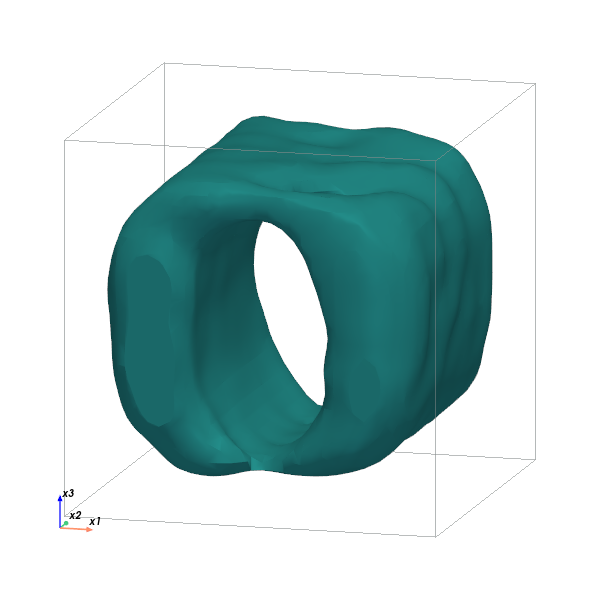}
		\caption{$\sigma=1\%$}
	\end{subfigure}
	\begin{subfigure}[t]{0.24\textwidth}
		\centering
		\includegraphics[width=\linewidth]{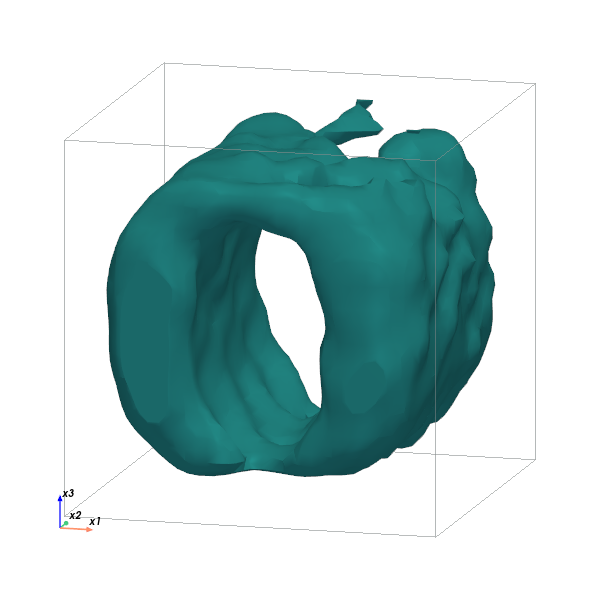}
		\caption{$\sigma=3\%$}
	\end{subfigure}\hfill
	\begin{subfigure}[t]{0.24\textwidth}
		\centering
		\includegraphics[width=\linewidth]{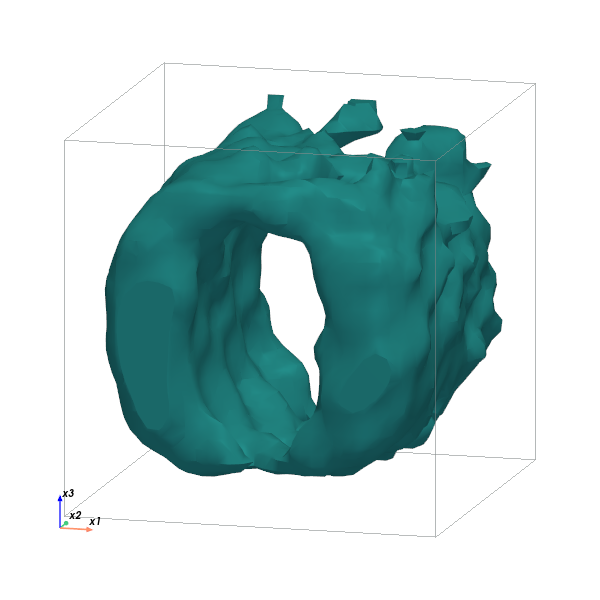}
		\caption{$\sigma=5\%$}
	\end{subfigure}
	
	\medskip
	
	\begin{subfigure}[t]{0.24\textwidth}
		\centering
		\includegraphics[width=\linewidth]{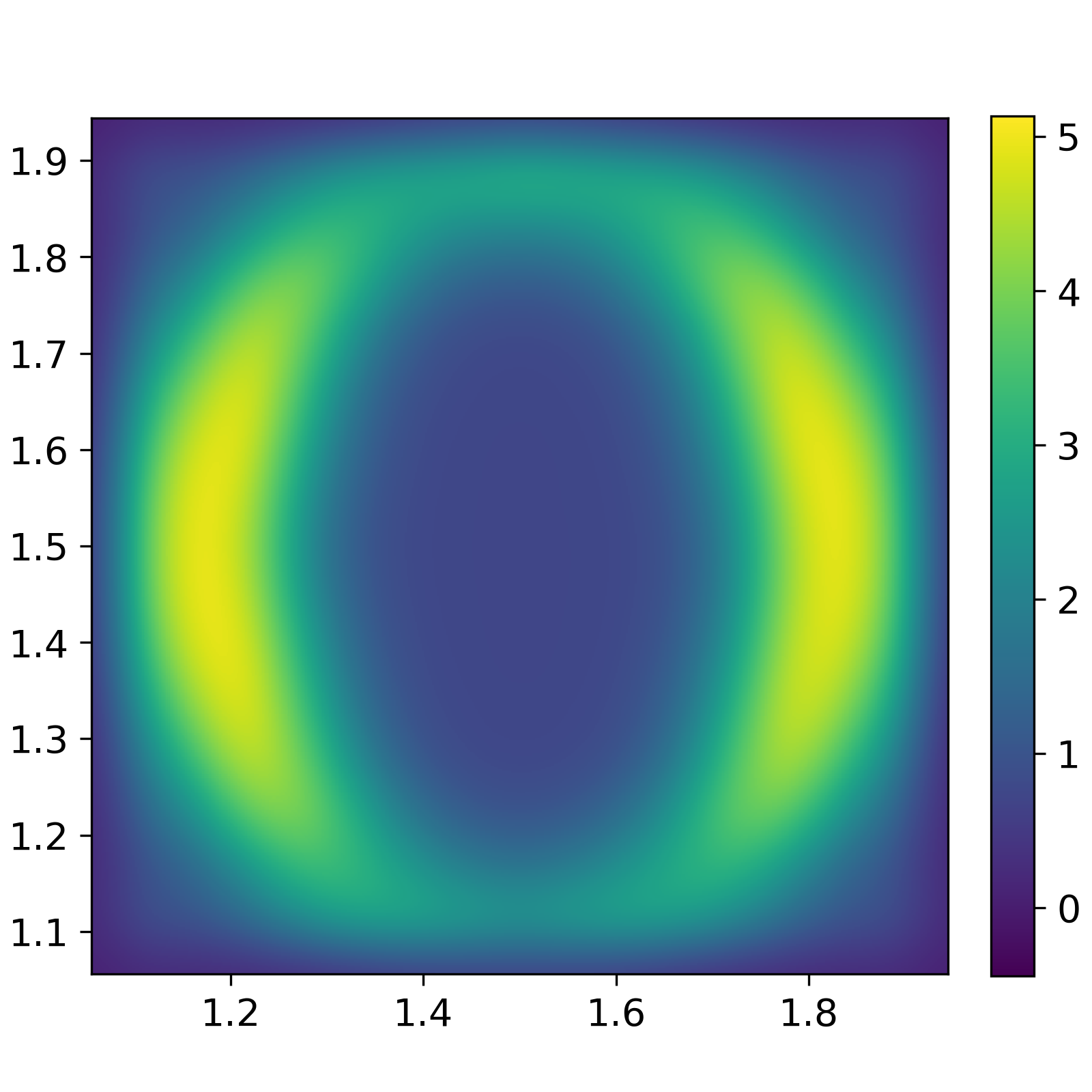}
		\caption{True coefficient}
	\end{subfigure}\hfill
	\begin{subfigure}[t]{0.24\textwidth}
		\centering
		\includegraphics[width=\linewidth]{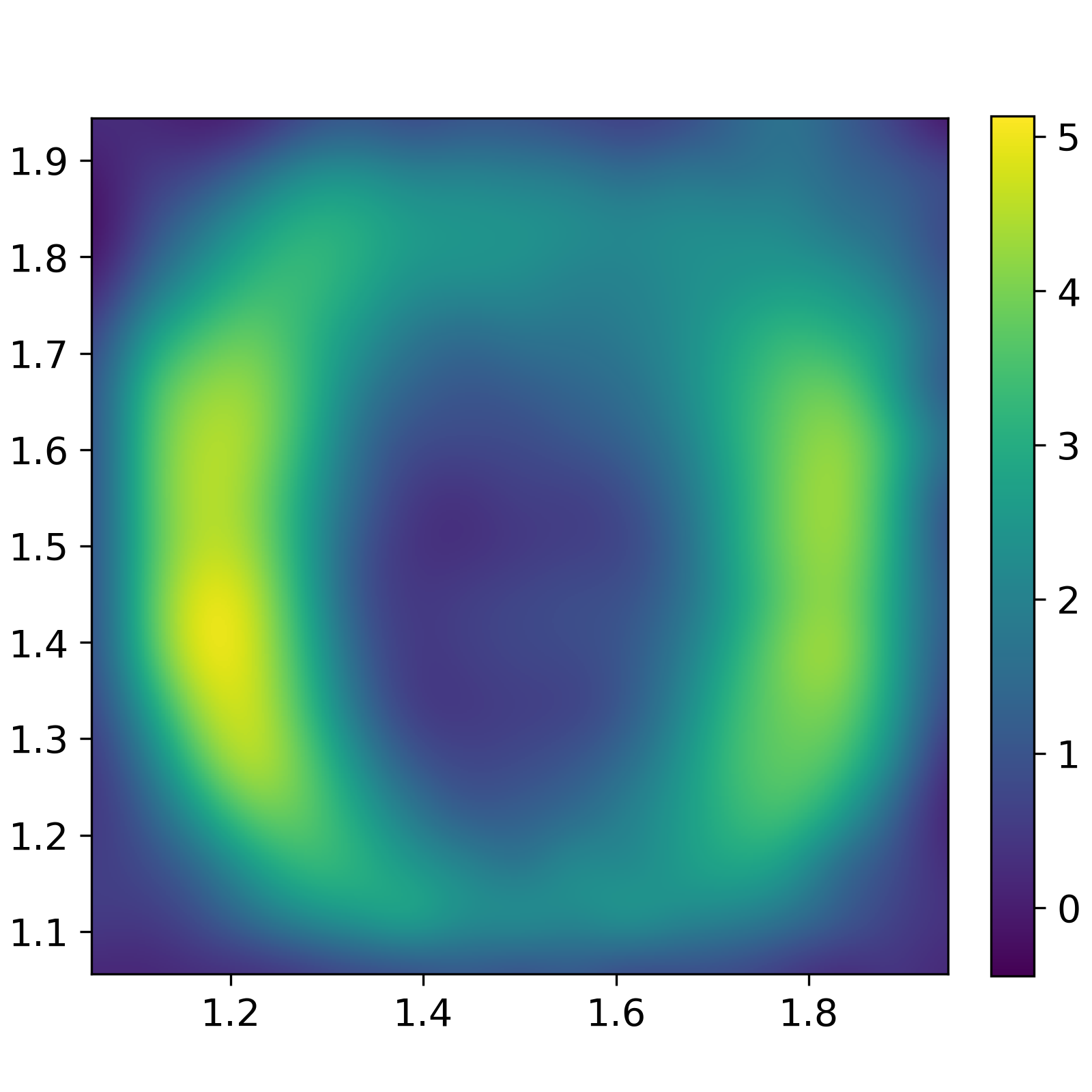}
		\caption{$\sigma=1\%$}
	\end{subfigure}
	\begin{subfigure}[t]{0.24\textwidth}
		\centering
		\includegraphics[width=\linewidth]{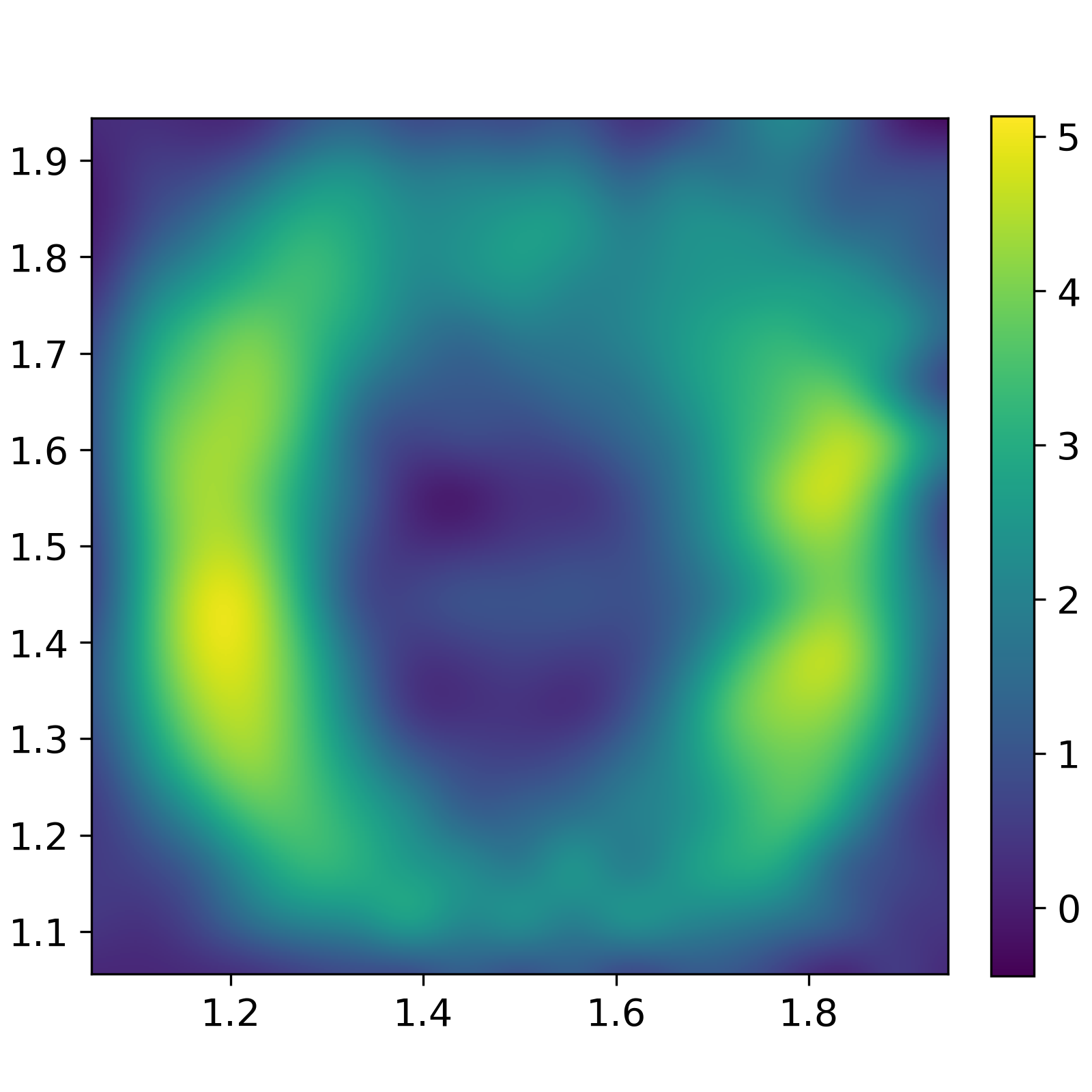}
		\caption{$\sigma=3\%$}
	\end{subfigure}\hfill
	\begin{subfigure}[t]{0.24\textwidth}
		\centering
		\includegraphics[width=\linewidth]{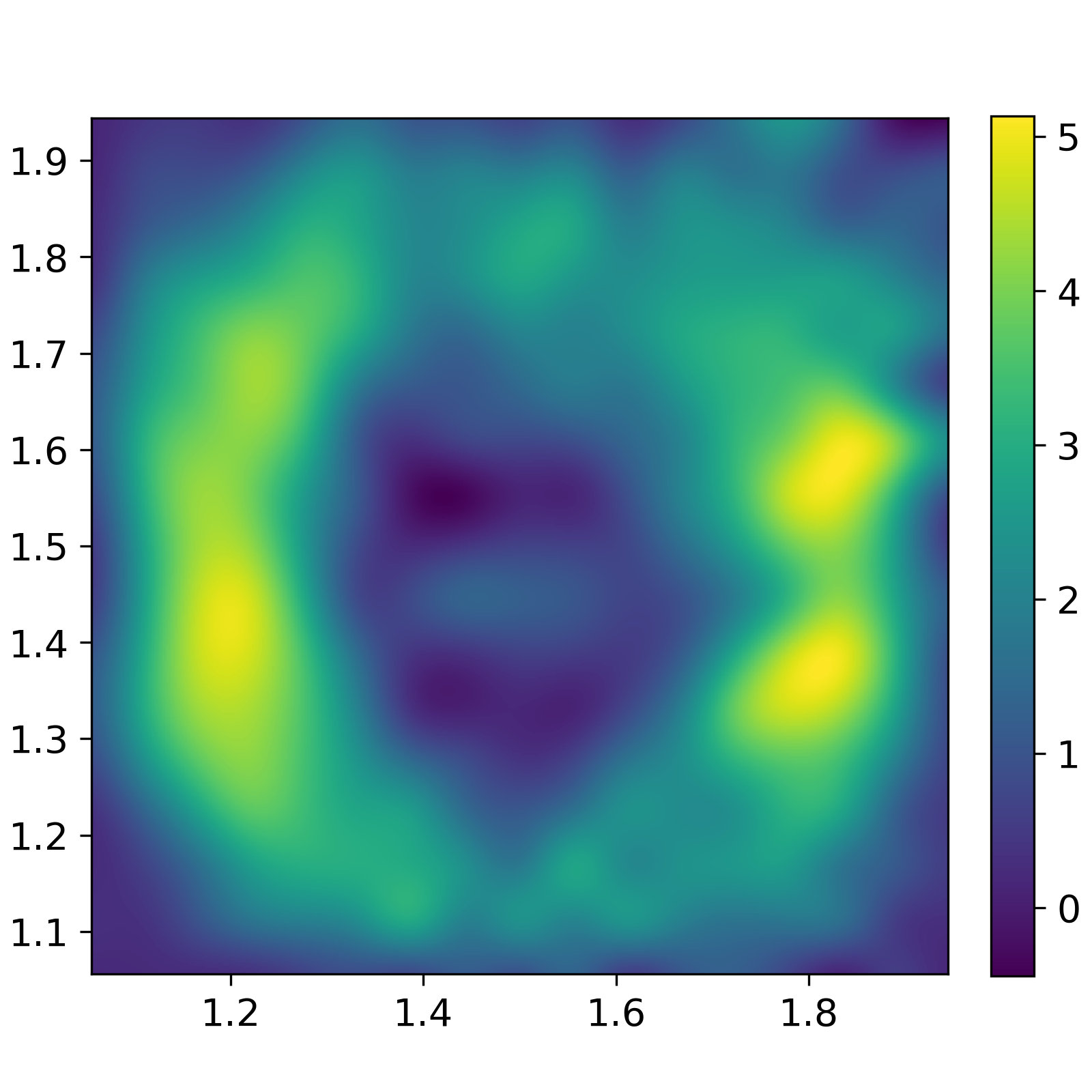}
		\caption{$\sigma=5\%$}
	\end{subfigure}
	\caption{Test \ref{ex:noise}: Three-dimensional and cross-sectional views of the true and reconstructed
		$O$-shaped coefficients for different noise levels.}
	\label{fig:noise_test}
\end{figure}

\begin{test}\label{ex:noise}
	We investigate the impact of the noise level $\sigma $ in the boundary data (%
	\ref{10.4}), (\ref{10.5}) via imaging an inclusion of the shape of the
	letters '$O$'. The true coefficient is $q(\mathbf{x})=5$ inside of this
	inclusion and $q(\mathbf{x})=0$ outside of it.  Noise levels of $1\%$%
	, $3\%$, $5\%$ are added to the measurements, i.e. $\sigma =0.01,$ $0.03,$ $%
	0.05$. The results are displayed on Figure~\ref{fig:noise_test}.
\end{test}

\bigskip \textbf{Data Availability:} Data sharing is not applicable to this
article as no datasets were generated or analysed during the current study.

\textbf{Declaration:} The authors declare that artificial intelligence was
not involved in any part of this work.

\textbf{Conflict of interest:} There are no conflicts of interest regarding
the publication of this paper.

\textsc{Department of Mathematics and Statistics, University of North
Carolina at Charlotte, Charlotte, North Carolina, 28223, United States}

\emph{Email address:} mklibanv@charlotte.edu

\textsc{Department of Mathematics, Southern University of Science and
Technology, Shenzhen 518055, China}

\emph{Email address:} li.jz@sustech.edu.cn

\textsc{Department of Mathematics, Southern University of Science and
Technology, Shenzhen 518055, China}

\emph{Email address}: niut@sustech.edu.cn


\begin{thebibliography}{99}
\bibitem{Bak} A. B. Bakushinskii, M. V. Klibanov and N. A. Koshev, Carleman
weight functions for a globally convergent numerical method for ill-posed
Cauchy problems for some quasilinear PDEs, \emph{Nonlinear Analysis: Real
World Applications}, 34, 201-224, 2017.

\bibitem{Beilina} L. Beilina, Domain decomposition finite element/finite
difference method for the conductivity reconstruction in a hyperbolic
equation, \emph{Communications in Nonlinear Science and Numerical Simulation}%
, 37, 222-237, 2016.

\bibitem{Beilina2} L. Beilina and E. Lindstr\"{o}m, An adaptive finite
element/finite difference domain decomposition method for applications in
microwave imaging, \emph{Electronics}, 11, 2079-9292, 2022.

\bibitem{Bell} M. Bellassoued and M. Choulli, Logarithmic stability in the
dynamical inverse problem for the Schr\"{o}dinger equation by arbitrary
boundary observation, \emph{J. de Math\'{e}matiques Pures et Appliqu\'{e}es}%
, 91, 233--255, 2009.

\bibitem{Ber} Yu. M. Berezansky, On the uniqueness theorem in the inverse
problem of spectral analysis for the Schr\"{o}dinger equation,\emph{\
Transl., Ser. 2, American Mathematical Society}, 35, 167--235, 1964.

\bibitem{Boul} M. Boulakia, M. de Buhan, T. Delaunay, S. Imperiale and P.
Moireau, Solving inverse source wave problem - from Carleman estimates to
observer design, \emph{Mathematical Control and Related Fields}, 16, 40-81,
2026.

\bibitem{BukhKlib} A.~L. Bukhgeim and M.~V. Klibanov, Uniqueness in the
large of a class of multidimensional inverse problems, \emph{Soviet
Mathematics Doklady}, 17, 244--247, 1981.

\bibitem{Chavent} G. Chavent, \emph{Nonlinear Least Squares for Inverse
Problems: Theoretical Foundations and Step-by-Step Guide for Applications},
Springer, Dordrecht, 2009.

\bibitem{CK} D. Colton and R. Kress, \emph{Inverse Acoustic and
Electromagnetic Scattering Theory}, Springer-Verlag, Berlin, 1992.

\bibitem{Hamr} M. Hamrouni, M. Khenissi and E. Soccorsi, Stability analysis
of inverse problems for coupled magnetic Schr\"{o}dinger equations, \emph{%
Inverse Problems}, 41, 035009, 2025.

\bibitem{Isakov} V. Isakov, \emph{Inverse Problems for Partial Differential
Equations}, Springer, New York, 2006.

\bibitem{Khoa} V.A. Khoa, M.V. Klibanov and L.H. Nguyen, Convexification for
a three-dimensional inverse scattering problem with the moving point source, 
\emph{SIAM J. Imaging Sciences}, 13, 871-904, 2020.

\bibitem{Klib92} M.~V. Klibanov. \newblock Inverse problems and {C}arleman
estimates. \newblock {\em Inverse Probl.}, 8, 575--596, 1992.

\bibitem{Klib95} M.V. Klibanov and O.V. Ioussoupova, Uniform strict
convexity of a cost functional for three dimensional inverse scattering
problem, \emph{SIAM\ J.\ Mathematical Analysis}, 26, 147-179, 1995.

\bibitem{Klib97} M.V. Klibanov, Global convexity in a three-dimensional
inverse acoustic problem, \emph{SIAM J. Mathematical Analysis,} 28,
1371--1388, 1997.

\bibitem{Ksurvey} M.V. Klibanov, Carleman estimates for global uniqueness,
stability and numerical methods for coefficient inverse problems,\emph{\ J.
of Inverse and Ill-Posed Problems}, 21, 477-510, 2013.

\bibitem{KR1} M.V. Klibanov and V.G. Romanov, The first solution of a long
standing problem: Reconstruction formula for a 3-d phaseless inverse
scattering problem for the Schr\"{o}dinger equation, \emph{J. Inverse and
Ill-Posed Problems}, 23, 415--428, 2015.

\bibitem{KL} M.V. Klibanov and J. Li, \emph{Inverse Problems and Carleman
Estimates: Global Uniqueness, Global Convergence and Experimental Data}, De
Gruyter, 2021.

\bibitem{Kind} M. V. Klibanov, V. A. Khoa, A. V. Smirnov, Loc H. Nguyen, G.
W. Bidney, Lam H. Nguyen, A. Sullivan and V. N. Astratov, Convexification
inversion method for nonlinear SAR imaging with experimentally collected
data, \emph{J. Applied and Industrial Mathematics}, 15, 413--436, 2021.

\bibitem{KR} M.V. Klibanov and V.G. Romanov, A H\"{o}lder stability estimate
for a coefficient inverse problem for the wave equation with a point source, 
\emph{Eurasian J. of Mathematical and Computer Applications}, 10, Issue 2,
11 --25, 2022.

\bibitem{TFD} M.V. Klibanov, A new type of ill-posed and inverse problems
for parabolic equations, \emph{Communications on Analysis and Computation},
2, 367-398, 2024.

\bibitem{MFG} M.V. Klibanov and J. Li, \emph{Carleman Estimates in Mean
Field Games}, De Gruyter, 2025.

\bibitem{KT} M.V. Klibanov and A.A. Timonov, Acoustic imaging via a
viscosity approximation of an elliptic system generated by the Lavrent'ev
integral operator, \emph{SIAM J. Imaging Sciences}, 18, 1002-1027, 2025.

\bibitem{moving} M.V. Klibanov, J. Li, V.G. Romanov and Z. Yang,
Convexification numerical method for imaging of moving targets, \emph{%
Mathematische Annalen}, 396:25, 2026.

\bibitem{LRS} M.M. Lavrent'ev, V.G. Romanov and S.P. Shishatskii, \emph{%
Ill-Posed Problems of Mathematical Physics and Analysis}, AMS, Providence,
RI, 1986.

\bibitem{Liu} S. Liu, A. Pierrottet and S. Scruggs, Stability analysis of
inverse problems for coupled magnetic Schr\"{o}dinger equations, \emph{%
Inverse Problems}, 39, 105009, 2023.

\bibitem{Nov} R. G. Novikov, The $\overline{\partial }-$approach to
approximate inverse scattering at fixed energy in three dimensions, \emph{%
Int. Math. Res. Papers,} 6, 287--349, 2005.

\bibitem{QY} H. Qiao and G. Yuan, Convexification of a complex
Ginzburg--Landau equation with reversed time using a new Carleman estimate, 
\emph{Mathematics}, 14, 3078, 2026.

\bibitem{Rom} V.G. Romanov, \emph{Inverse Problems of Mathematical Physics},
VNU Press, 1986.

\bibitem{Rom1} V.G. Romanov, \emph{Investigation Methods for Inverse
Problems,} VSP, Utrecht, 2002.

\bibitem{Sam} A.A. Samarskii, \emph{The Theory of Difference Schemes}, Boca
Raton, 2001.

\bibitem{Scales} J.A. Scales, M.L. Smith, T.L. Fisher, Global optimization
methods for multimodal inverse problems, \emph{J. Computational Physics},
103, 258-268, 1992.

\bibitem{T} A. N. Tikhonov and V. Y. Arsenin, \textit{Solutions of Ill-Posed
Problems}, Winston, Washington, DC, 1977.
\end{thebibliography}
\end{document}